%% file: Siam_article.tex
\documentclass[onefignum,onetabnum]{siamonline250211}
\usepackage{comment}
\usepackage{physics}
\usepackage{graphicx}
\usepackage{listings}
\usepackage{enumitem}
\usepackage{url}
\usepackage{kantlipsum}
\usepackage{tabularx}
\usepackage[utf8]{inputenc}
\usepackage{csquotes}
\usepackage[nottoc]{tocbibind}
\usepackage{xcolor}
\usepackage{makecell}
\usepackage{float}

\DeclareMathOperator{\nullspace}{null}
\DeclareMathOperator*{\argmin}{arg\,min}

\input{ex_shared}

\ifpdf
\hypersetup{
  pdftitle={Forced-Term Modeling in the Koopman Framework for Linear Operator Learning in Nonlinear Dynamical Systems},
  pdfauthor={Paolo Climaco, Jochen Garcke and Xenia F. Gerloff}
}
\fi

\makeatletter
\providecommand{\@openbib@code}{}
\makeatother

\begin{document}

\maketitle

\begin{abstract}
We study the problem of extracting reliable and interpretable linear models from nonlinear dynamical systems with periodic behavior. Dynamic Mode Decomposition (DMD) is a widely used data-driven method for approximating dynamics by a linear system and is commonly interpreted through Koopman operator theory. However, because standard DMD assumes linear evolution in the chosen observables, it can produce misleading results on nonlinear systems. Recent extensions address this limitation by modeling the dynamics as a linear system driven by a nonlinear forcing term, but their relationship to Koopman theory and their reliability remain unclear. We observe that the Koopman theoretical framweork allows for explicit, finite-dimensional, forced linear representations of dynamical systems, independently of the chosen observable space. The main question is therefore how to choose a useful, computable, and interpretable model for the forcing term, since forced linear representations are generally nonunique. In the framework considered here, this modeling question is related to the choice of a defect space: a complementary space that closes the action of the Koopman operator on the chosen observables. Existing forced-DMD approaches can then be related to different modeling choices for this defect space. Importantly, the structure of the learned linear operator depends on the choice of the defect space. As an example of defect-space modeling for algorithm design, we introduce Correlation-Basis enhanced DMD (CB-DMD), a data-driven method whose goal is to learn an effective autonomous linear operator for periodic nonlinear dynamics. The forcing term is used to absorb nonlinear components that cannot be represented by the learned linear operator. We further analyze two existing DMD variants for forced linear modeling and show that they often produce poor operator estimates and spurious modes. Experiments on three nonlinear periodic systems demonstrate that CB-DMD yields more reliable linear models than standard DMD and related variants in the tested examples.
\end{abstract}



\section{Introduction}
~\label{koopmanforsystemidentification}

In this work we consider the analysis of data arising from nonlinear dynamical systems evolving on a limit cycle. That is, nonlinear systems that exhibit a periodic behavior. Analyzing data from such nonlinear dynamical systems is important due to its wide applicability in domains such as neuroscience (e.g., neural oscillations), physics (e.g., superconductivity), and engineering (e.g., fluid dynamics). These systems exhibit inherently periodic but nonlinear behavior. While nonlinear models can capture complex dynamics, they often lack interpretability. Linear models, on the other hand, are easier to analyze and deploy, but struggle to represent periodic nonlinear behavior. Therefore, developing methods to extract interpretable linear representations from nonlinear, limit-cycle systems can enhance their scientific understanding. 

Dynamic Mode Decomposition (DMD)~\cite{schmidt,dmd222} is a data-driven algorithm for constructing linear models from time-series data, with applications in fields such as epidemiology \cite{infectDMD} and engineering \cite{Climaco2023}. DMD characterizes the evolution of finite dynamical systems as a linear combination of spatio-temporal structures, called DMD-modes, which can be interpreted as the eigenvectors of the linear operator evolving the system. 

DMD relies on the underlying assumption of linear evolution of the analyzed data measurements, which may not be true and may affect the reliability and interpretability of the computed linear models. For instance, one notable issue of DMD is the creation of spurious modes \cite{resDMD}, that is, the computation of DMD-modes that may not have any physical meaning or theoretical interpretation. To address this, several DMD-based approaches have been developed, including: coordinate transformations to linearize the system (e.g., time-delay embedding \cite{time-delay}, neural network-based observables \cite{Lusch2018, Mardt2018});  decomposition into a linear component, whose evolution is governed by a linear operator, and a nonlinear component, capturing the remaining nonlinearities \cite{chao,PCT,Khodkar2019};  or learning linear operators under specific constraints \cite{Baddoo2023}. 

In this work, we focus on DMD-based approaches aiming at characterizing the dynamics evolution as a forced linear model, i.e., a system that is governed by linear equations and is under an external force. 
The aim is the interpretability of the linear operator as learned by means of Koopman theory.
Koopman operator theory, introduced in 1931 \cite{koopm}, allows for linear representations of nonlinear systems by considering the dynamics on the observable space, which is a function space defined on the state variables of the system. In the observable space, the dynamics are ruled by a linear but infinite-dimensional operator, the Koopman operator, which encodes all the information related to the dynamical system. Thus, the central focus of several numerical techniques to analyse dynamical systems consists in learning the Koopman operator from data. DMD can provide finite-dimensional approximations of the Koopman operator by computing the best linear fit of measurements representing the temporal evolution of a dynamical system~\cite{onconv}.

This Koopman-theoretic perspective also clarifies how forced linear models arise. Since the Koopman operator acts on an infinite-dimensional observable space, a finite-dimensional approximation requires restricting attention to a prescribed finite-dimensional space of observables. In general, however, the Koopman image of this space is not contained within the space itself. We therefore decompose the Koopman action into two components: one that remains in the chosen observable space, and one that does not. The first component is represented by a finite-dimensional linear operator, while the second is represented through a finite-dimensional space, which we call the defect space, and appears as the forcing term in the model. Since the defect space is not unique, its choice is part of the modeling procedure. Consequently, the resulting linear operator and forcing term are meaningful only relative to the chosen observable and defect spaces. Within this framework, existing forced linear representation methods, including those in~\cite{PCT} and~\cite{Khodkar2019}, can be interpreted as DMD-based techniques corresponding to particular choices of the defect space.

Based on this viewpoint and on dynamic mode decomposition with control (DMDc)~\cite{DMDc}, we develop ``Correlation-Basis enhanced DMD'' (CB-DMD), a data-driven approach for computing forced linear decompositions of nonlinear systems. The underlying idea is that, once an observable space has been fixed, the unresolved part of the Koopman action can be modeled through a suitable finite-dimensional defect space. In CB-DMD, the choice of this defect space is tied to the eigenfunctions of an associated correlation operator, which are used to identify directions separated from the dominant snapshot row space. This construction yields a decomposition of the observed dynamics into an autonomous linear component on the chosen observables and a data-dependent residual component determined by the defect space. In this sense, CB-DMD provides a Koopman-consistent forced linear representation in which the modeling role of the observable and defect spaces is explicit, while mitigating spurious modes relative to classical DMD and existing forced-linear DMD variants.

CB-DMD does not require an a priori model for the control or forcing term, and it does not involve an iterative learning process. 
The procedure can be applied directly in the coordinates in which the data are measured, including the original state coordinates, and therefore does not require a time-delay embedding or another coordinate transformation.

\section{Related work}

Dynamic Mode Decomposition (DMD) is a data-driven technique for analyzing the dynamics of complex systems. Originally developed in fluid dynamics, DMD has since been applied in fields such as neuroscience, finance, video processing, and epidemiology. The main idea behind DMD is to build linear models of dynamical systems from time-series data, without needing to know the underlying equations. One of the advantages of DMD is its interpretability. The linear model computed by DMD can be interpreted by means of Koopman operator theory. However, DMD assumes that the system's behavior can be described by a linear operator, but this assumption often fails for nonlinear systems, which limits its reliability.
To address these limitations, several extensions of DMD have been proposed. These extensions of DMD follow at least three different approaches: Linearizing the dynamics by applying a coordinate transformation to the data, e.g., involving a data-driven approximation of transfer operators~\cite{Klus2018}; decomposing the dynamics into a linear part and a nonlinear forcing term; learning linear operators under specific constraints to ensure physical consistency. Examples of such approaches are Extended DMD (EDMD)~\cite{EDMD}, DMD with control (DMDc)\cite{DMDc}, and Pysics informed DMD (piDMD) \cite{Baddoo2023}, respectively. EDMD consists of applying DMD on vector-valued transformations of the measurements, which aim to linearize the dynamics evolution~\cite{time-delay,Klus2018,kooppde,Lusch2018, Mardt2018,williams2015}. DMDc extends DMD to systems with external inputs, allowing for the modeling of forced linear systems. 
piDMD learns linear operators under specific constraints to ensure physical consistency.

In this work we focus on DMD-based approaches building forced linear representations of the dynamics by adding an external forcing term to a learned linear model. The effectiveness of forced linear characterizations of nonlinear, even chaotic, dynamics was first demonstrated by a data-driven strategy called HAVOK~\cite{chao}, based on time-delay embedding. Other methods, inspired by HAVOK, can be found in~\cite{PCT} and~\cite{Khodkar2019}. In \cite{PCT} the forcing component is learned from the one-step prediction error of a linear model of the dynamics obtained through the DMD algorithm. In \cite{Khodkar2019} the forcing term is approximated as a linear combination of polynomial basis functions. Both approaches rely on time-delay embedding and Takens theorem~\cite{takens}. It is worth noting that the operator resulting from time-delayed measurements may not have a straightforward connection to the original coordinates, which can limit its interpretability. Moreover, it is unclear whether the effectiveness of these methods is due to the use of time delay embedding, the forced linear model they build or the combination of the two. That is, it is not clear if they would be equally effective using the original coordinates. Additionally, it's important to consider that the above mentioned methods are  motivated by either empirical results~\cite{chao}, or heuristic arguments~\cite{PCT}~\cite{Khodkar2019}, therefore lack of a sound theoretical motivation. 
Recently, a framework for forced nonlinear multiple-input multiple-output dynamical systems that integrates Hankel Dynamic Mode Decomposition with control and transfer operators was proposed in~\cite{math14040625}. 

While these strategies have shown to be effective in identifying nonlinear dynamical components and providing forced linear characterizations of nonlinear systems, their effectiveness and usefulness in building an interpretable linear operator characterizing the evolution of the linear component of the dynamics is still to be determined.
\section{Theoretical background}\label{theoretical_background}
In this section we recapture the key concepts behind Koopman operator theory and Dynamic Mode Decomposition (DMD) that are relevant to our work. We start with a simple introduction to discrete-time dynamical systems, then we present the Koopman operator and the concept of Koopman invariant subspaces. Next, we introduce DMD and its variant,  Dynamic Mode Decomposition with control (DMDc)~\cite{DMDc}.

\subsection{Discrete-time dynamical systems}
We consider discrete-time dynamical systems of the form 
\begin{equation}
    \label{dynamicalsystem}
    \bar{\mathbf{x}}=T(\mathbf{x}),
\end{equation}
where $\mathbf{x}, \bar{\mathbf{x}}\in \mathcal{M} \subset \mathbb{R}^n$, with $\mathcal{M}$ a set representing the
state space. We assume $\mathcal{M}$ to be a measure space with an associated measure $\mu$. 
Equation (\ref{dynamicalsystem}) defines the systems's evolution: for any state $\mathbf{x}$, the next state $\bar{\mathbf{x}}$ is given by the map $T:\mathcal{M} \rightarrow \mathcal{M}$.
We call any sequence of states $\{\mathbf{x}_m\} \subset \mathcal{M}$ such that $\mathbf{x}_{m+1}=T(\mathbf{x}_m)$ a \textit{trajectory} of the dynamical system. 

\subsection{The Koopman operator}
We can associate the dynamical system in ($\ref{dynamicalsystem}$) with an operator $\mathcal{K}$ acting on an infinite-dimensional Hilbert space,
 $\mathcal{H}:=L^{2}_{\mu}(\mathcal{M}; \mathbb{C})$ of functions from $\mathcal{M}$ into $\mathbb{C}$, square integrable with respect to the measure $\mu$. We can define $\mathcal{K}: \mathcal{H} \rightarrow \mathcal{H}$ as follows:
\begin{equation}
    \mathcal{K}g=g \circ T ,
\end{equation}
where $g \in \mathcal{H}$. The operator $\mathcal{K}$ is called the Koopman operator and the space $\mathcal{H}$ is called the observables space, where an
observable is a function $g \in \mathcal{H}$.  The Koopman operator advances the dynamics in time in the observables space, i.e.,  $g(\mathbf{x}_{m+1})= \mathcal{K}g(\mathbf{x}_m)$. 
The Koopman operator is linear on its domain. In general, defining
$\mathcal K g=g\circ T$ as an operator on all of $L^2_\mu(\mathcal M;\mathbb C)$
requires additional assumptions on $T$ and $\mu$. In this work we only use the
action of $\mathcal K$ on the finite-dimensional observable spaces introduced
below. We therefore assume that the observables considered, and their
compositions with $T$, belong to $L^2_\mu(\mathcal M;\mathbb C)$. Under this
standing assumption, all Koopman expressions appearing in the finite-dimensional
representations below are well defined. The use of an $L^2$ observable space is
standard in work relating DMD-type algorithms to Koopman theory~\cite{onconv}.

\subsection{Koopman invariant subspaces}
\label{Koopman invariant subspace}
Now, we explain how the Koopman theoretical framework can be used to find finite linear representations of nonlinear systems by introducing the concept of
Koopman invariant subspaces already investigated in~\cite{invariant}.
Consider a finite-dimensional subset $\mathcal{I} \subset \mathcal{H}$ of the observable space. We say that $\mathcal{I}$ is invariant under the action of the Koopman operator if for any $g \in \mathcal{I}$ we have that $\mathcal{K}g \in \mathcal{I}$. This property is very useful in applications, and now we see why. Assume we have linearly independent observables $\{g_i\}_{i=1}^N, \; g_i \in \mathcal{H}$, such that $\mathcal{H}_N:= \operatorname{span}\{g_i\}_{i=1}^N$ is a Koopman invariant subspace. Thus, we have that $\mathcal{K}g_i\in \mathcal{H}_N$ for $i=1,\dots,N$. Consequently, since the $\{g_i\}_{i=1}^N$  form a basis of $\mathcal{H}_N$, there exist coefficients $\{K_{ij}\}_{j=1}^N \subset \mathbb{C}$ such that
\begin{equation}
  \mathcal{K}g_i=\sum_{j=1}^N K_{ij}g_j, \quad i=1,..,N,
\end{equation}
which is written in a vector-valued fashion yields
\begin{equation}
  \label{finite_dimensional_approx}
  \mathcal{K} \mathbf{ \hat{\mathbf{g}}}= \mathbf{K} \mathbf{ \hat{\mathbf{g}}},
\end{equation}
where $ \mathbf{ \hat{\mathbf{g}}}=[g_1,\dots,g_N]^T$ and $(\mathbf{K}_{i,j})=K_{ij}$. Equation ($\ref{finite_dimensional_approx}$) tells us that if we, to study a dynamical system, choose a set of linearly independent observables that span a Koopman invariant subspace, we can represent the evolution of any trajectory $ \mathbf{x}_k\in \mathcal{M}$ in a linear fashion by considering the observation of such a trajectory through those observables, independently of how complicated the 
underlying dynamics is, i.e.,
\begin{equation}
  \label{linear representation}
   \mathbf{ \hat{\mathbf{g}}}( \mathbf{x}_{k+1})= \mathcal{K} \mathbf{ \hat{\mathbf{g}}}( \mathbf{x}_k)= \mathbf{K} \mathbf{ \hat{\mathbf{g}}}( \mathbf{x}_k).
\end{equation}

\subsection{Dynamic Mode Decomposition (DMD)}
\label{numerical_algo}
We now introduce the dynamic mode decomposition~\cite{schmidt}, which is the principal algorithm used within the Koopman operator framework to compute data-driven linear models.
DMD was introduced in the fluid dynamic community~\cite{schmidt,dmd222}, to extract patterns and predict the future state of 
fluid flows. After its original formulation, many variants of DMD have been developed to address new tasks \cite{DMDc, cDMD}, or to make it more efficient \cite{optDMD, mrDMD,HODMD, debiasingDMD, leastDMD}. Moreover, DMD is strongly connected with Koopman operator theory. In particular, it has been shown that it can compute a finite-dimensional approximation of the Koopman operator~\cite{ergoticDMD,onconv} by building a best linear model that advances the measurements in time.

Given measurements (or snapshots) $\{ \mathbf{x}_k\}_{k=0}^M\subset \mathbb{C}^n$ representing the temporal evolution of a dynamical system, DMD computes the
eigendecomposition of the matrix $\mathbf{A}$ that gives the best linear fit of the measurements $\{ \mathbf{x}_k\}$:
\begin{equation}
   \mathbf{x}_{k+1}\approx \mathbf{A}  \mathbf{x}_k.
\end{equation}
Turning to the details, DMD computes eigenvalues and eigenfunctions of the matrix $\mathbf{A}$ that solves the following minimization problem in the Frobenius norm~\cite{Klus2018} 
\begin{equation}
\label{minprob}
\min_{\mathbf{A} }\|\bar{\mathbf{X}}-\mathbf{AX}\|_F,
\end{equation}
where
\begin{equation}
\label{snapshots}
\mathbf{X}=[\mathbf{x}_0,\mathbf{x}_1,\dots,\mathbf{x}_{M-1}] \quad \text{and} \quad \bar{\mathbf{X}}=[\mathbf{x}_1,\mathbf{x}_2,\dots,\mathbf{x}_M]  
\end{equation} 
are the so-called snapshots matrices where each column represents a state of the dynamical system at different times. The columns are arranged chronologically, and the columns of $\bar{\mathbf{X}}$ are shifted one time step forward with respect to the columns of $\mathbf{X}$. The solution to $(\ref{minprob})$ is given by the matrix
\begin{equation}
\mathbf{A}=\bar{\mathbf{X}}\mathbf{X}^{\dagger},
\end{equation}
where $\mathbf{X}^{\dagger}$ is the Moore-Penrose inverse of the matrix $\mathbf{X}$. In principle, the matrix $\mathbf{A}$ can be directly computed from data. However, it may be advantageous to compute a low-rank approximation $\tilde{\mathbf{A}}$ of $\mathbf{A}$. The eigenvalues and eigenvectors of $\mathbf{A}$ called DMD eigenvalues and DMD modes, respectively, can be then recovered from the eigenvalues and eigenvectors of $\tilde{\mathbf{A}}$. DMD modes and eigenvalues can be later used to reconstruct the original measurements and the operator $\mathbf{A}$ solution of ($\ref{minprob}$). Algorithm~\ref{alg:DMD}, in Appendix E,
shows how DMD is implemented. For a general overview of the DMD approach see~\cite{dmd222,dmd}.

The DMD algorithm builds a best-fit linear model for the observed dynamics. However, for nonlinear systems, this linear assumption is often invalid. Extended Dynamic Mode Decomposition (EDMD)~\cite{EDMD} addresses this by applying DMD to nonlinearly transformed measurements $\{\hat{\mathbf{g}}(\mathbf{x}_i)\}_{i=1}^N$, where $\hat{\mathbf{g}} = [g_1, g_2, \dots, g_m]$ are chosen to ideally span a Koopman invariant subspace, thus enabling a linear representation as in (\ref{linear representation}). In practice, finding such observables is difficult, and for systems with multiple fixed points or general attractors, finite-dimensional Koopman invariant subspaces may not exist~\cite{invariant}. As a result, recent research also explores data-driven approaches for constructing meaningful forced linear representations of nonlinear dynamics~\cite{Brunton2019}.

\subsection{Dynamic Mode Decomposition with control (DMDc)}
\label{dmdcsection}

Dynamic Mode Decomposition with control (DMDc) is an extension of the classical
DMD algorithm that enables forced linear representations of nonlinear systems. It was developed in~\cite{DMDc} to address issues arising in the
context of system identification and control. DMDc is based on the assumption
that the snapshots $ \mathbf{x}_k \in \mathbb{C}^n, k=0,...,M$ arise from a
linear controlled system, which means that the linear dynamical system
connecting the future state $ \mathbf{x}_{k+1}$ relies on information from both the
current state $ \mathbf{x}_k \in \mathbb{C}^n$ and the current control $\mathbf{u}_k \in
\mathbb{C}^l$, i.e., 
\begin{equation}
\label{forcingeq}
 \mathbf{x}_{k+1}=\mathbf{A} \mathbf{x}_k+\mathbf{B}\mathbf{u}_k,
\end{equation}
where  $\mathbf{A}\in \mathbb{C}^{n \times n}$, $\mathbf{B}\in \mathbb{C}^{n,l}$ and $k=0,...,M-1$.
Equation $(\ref{forcingeq})$ can be also written in a matrix formulation as follows
\begin{equation}
\bar{\mathbf{X}}=\mathbf{AX} + \mathbf{B\Upsilon}\\,
\end{equation}
where $\bar{\mathbf{X}}, \mathbf{X}$ as in $(\ref{snapshots})$ and
$\mathbf{\Upsilon}=[\mathbf{u}_0,\mathbf{u}_1,...,\mathbf{u}_{M-1}]$. Thus, given $\bar{\mathbf{X}},\mathbf{X}$ and $\mathbf{\Upsilon}$ DMDc
algorithm computes the matrices $\mathbf{A} \in \mathbb{C}^{n \times n}$ and $\mathbf{B} \in
\mathbb{C}^{n \times l}$ that solve the following minimization problem
\begin{equation}
\label{miniproblemfor}
\min_{\mathbf{A},\mathbf{B}}\| \bar{\mathbf{X}}- \mathbf{AX}- \mathbf{B\Upsilon}\|_F.\\
\end{equation}
The analytical solution to $(\ref{miniproblemfor})$ is 
\begin{equation}
[\mathbf{A},\mathbf{B}]= \bar{\mathbf{X}}[\mathbf{X},\mathbf{\Upsilon}]^{\dagger}.
\end{equation}
To have a detailed explanation of the different steps performed by the algorithm
and a deeper understanding about the technicalities in it, see~\cite{DMDc}. To explicitly see the
connection between the minimization problem $(\ref{miniproblemfor})$  and the
solution of the DMDc algorithm look at~\cite{meko}. What is relevant in this work is that
the solution of the minimization problem $(\ref{miniproblemfor})$ provides us with matrices
$\mathbf{A}$ and $\mathbf{B}$ which, given the knowledge of the control terms $\mathbf{u}_k$, allow us to reconstruct
the evolution of the measurements $\{ \mathbf{x}_k\}$ as a forced system, i.e.,
\begin{equation}
  \label{characterisation_DMDc}
  \hat{\mathbf{x}}_{k+1}= \mathbf{A}\hat{\mathbf{x}}_k + \mathbf{B}\mathbf{u}_k,
\end{equation}
where $\hat{\mathbf{x}}_k$ represents the reconstructed time series.  
Notice that the characterization of the measurements' evolution as a forced linear system
performed in ($\ref{characterisation_DMDc}$) is only possible thanks to the knowledge
of the adequate control terms $\mathbf{u}_k$. In Section $\ref{proposed_method}$ we develop a numerical method
that, as DMDc, can find forced linear characterizations of dynamical systems, but
does not require any knowledge of the control terms $\mathbf{u}_k$, which can be fully recovered from 
data.
 
\section{Koopman-based forced linear representation of nonlinear dynamics}
~\label{Representing_chaos}
 Let us now introduce our theoretical framework, which is based on Koopman theory, to describe the evolution of nonlinear dynamical systems in a finite forced linear fashion once a finite observable space and a finite complement for its Koopman image have been fixed.
We use the concept of a Koopman almost-invariant subspace, an extension of the Koopman invariant subspace from Section~\ref{Koopman invariant subspace}.

\subsection{Koopman almost-invariant subspaces}
We say that a finite-dimensional linear subspace $\mathcal{E}\subset\mathcal{H}$
is almost-invariant under the action of the Koopman operator if there exists a
finite-dimensional linear subspace $\mathcal{V}\subset \mathcal{H}$ such that
\begin{equation}
  \label{intersection}
  \mathcal{E} \bigcap  \mathcal{V}=\{ 0 \},
\end{equation}
and such that for all $g \in\mathcal{E}$,
\begin{equation}
  \label{sum_spaces}
  \mathcal{K}g \in \mathcal{E} +  \mathcal{V},
\end{equation}
where $
\mathcal{E} +  \mathcal{V}:=\{g=h+\psi \;| \; h\in \mathcal{E}, \;\psi \in
\mathcal{V} \}$. We call $\mathcal{V}$ a \textit{defect space} of
$\mathcal{E}$. Notice that, since $\mathcal{K}g \in \mathcal{E} +  \mathcal{V}$
and $\mathcal{E} \bigcap  \mathcal{V}=\{ 0 \}$ there exists a unique
representation of $\mathcal{K}g$ in $\mathcal{E} +
\mathcal{V}$~\cite{Berberian1961}, i.e., 
\begin{equation}
  \mathcal{K}g=h+f,
\end{equation}
where $h\in \mathcal{E}$ and $f \in \mathcal{V}$ are uniquely defined for the chosen space $\mathcal V$. This uniqueness is therefore relative to a fixed defect space. Different complements can lead to different decompositions. Thus, the word ``defect'' should not be interpreted as defining a canonical object associated only with the dynamical system. It denotes a chosen complement used to separate the part of $\mathcal K\mathcal E$ represented in $\mathcal E$ from the part represented as forcing. The usefulness of such a decomposition depends on the additional modelling criteria used to choose $\mathcal V$, such as separation from $\mathcal E$, computability from data, numerical stability, and the interpretability of the resulting linear component. We now see that considering
Koopman almost-invariant subspaces enable us to build forced linear
representations of even chaotic systems.

\subsection{Forced linear representation}

Consider a trajectory $\{ \mathbf{x}_k\}_{k=0}^M\subset \mathcal{M}\subset \mathbb{R}^n$
of a dynamical system in the form of (\ref{dynamicalsystem}), and observables $\{g_i\}_{i=1}^N
\subset \mathcal{H}$ through which we observe the system. We define 
$\mathcal{H}_N:= \operatorname{span}\{g_i\}_{i=1}^N$, which by definition is a
finite-dimensional linear subspace of the observable space $\mathcal{H}$. Since every linear
finite-dimensional subspace of $\mathcal{H}$ is almost-invariant under the
action of the Koopman operator (Appendix A), there exists a defect space
$\mathcal{H_D}\subset\mathcal{H}$ of $\mathcal{H}_N$, which is a finite-dimensional
linear subspace with $\mathcal{H}_N \bigcap  \mathcal{H_D}=\{ 0 \}$, such that for
all $g \in \mathcal{H}_N$, $\mathcal{K}g \in \mathcal{H}_N + \mathcal{H_D}$. As a
consequence, once a defect space has been fixed, for each $i=1,\dots,N$ there exist uniquely defined $h_i\in \mathcal{H}_N$ and $ f_i \in \mathcal{H_D}$ such that
\begin{equation}
  \mathcal{K}g_i= h_i+ f_i. 
\end{equation}
 Furthermore, since $h_i \in \mathcal{H}_N$, by definition of $\mathcal{H}_N$ there exist 
 coefficients $K_{ij}\in \mathbb{C}$ such that
$h_i= \sum_{j=1}^N K_{ij} g_j$. Thus, we can write 
\begin{equation}
  \mathcal{K}g_i= \sum_{j=1}^N K_{ij} g_j + f_i, \quad i=1,\dots, N,
\end{equation}
which, in a vector-valued form, yields
\begin{equation}
  \label{forced_representation_in_H}
  \mathcal{K} \hat{\mathbf{g}}= \mathbf{K} \hat{\mathbf{g}} + \hat{\mathbf{f}},
\end{equation}
where $\hat{\mathbf{f}}=[f_1,\dots,f_N]^T$, $ \hat{\mathbf{g}}=[g_1,\dots,g_N]^T$ and $(\mathbf{K})_{i,j}=K_{ij}$.
Computing the vector-valued observable $ \hat{\mathbf{g}}$ along the trajectory $ \mathbf{x}_k$ we can re-write 
(\ref{forced_representation_in_H}) as
\begin{equation}
  \label{forced_representation}
   \hat{\mathbf{g}}( \mathbf{x}_{k+1})= \mathbf{K} \hat{\mathbf{g}}(\mathbf{x}_{k}) + \hat{\mathbf{f}}(\mathbf{x}_{k}).
\end{equation}
Equation~(\ref{forced_representation}) shows that, once a defect space has been fixed, the dynamics observed through the vector-valued observable $\hat{\mathbf g}$ admits a finite forced linear representation. Such a representation can be constructed for any finite set of observables $\{g_i\}_{i=1}^N$, since one can always choose a finite-dimensional complement to the Koopman image of $\mathcal H_N$. Thus, almost-invariance is not, by itself, a restrictive property and does not characterize a special class of nonlinear systems. Its role here is instead to make explicit the modelling choice involved in separating the Koopman action into a component represented in $\mathcal H_N$ and a defect contribution represented as forcing.

In many applications, the observables of primary interest are the state variables themselves. In this case, $N=n$ and $g_i(\mathbf{x}_k)=x_{k_i}$, where $x_{k_i}$ denotes the $i$-th component of $\mathbf{x}_k$. More generally, the same construction applies to any finite collection of observables. If the functions $\{g_i\}_{i=1}^N$ are linearly independent, they form a basis for $\mathcal H_N$; hence, for a fixed decomposition of $\mathcal K\mathcal H_N$ into $\mathcal H_N$ and the chosen defect space, the matrix $\mathbf K\in\mathbb C^{N \times N}$ in~(\ref{forced_representation}) is uniquely determined. Without fixing this defect space, neither the defect contribution nor the matrix $\mathbf K$ is canonical.
In this work, the central modelling task is therefore to select a defect space that is appropriate for the application, for example one that is separated from the observable space, computable from data, numerically stable, and yields an interpretable operator $\mathbf K$.

\section{Computing forced linear models}
\label{proposed_method}
We now use the concept of defect space within the Koopman-theretic formalism to motivate the computation of forced linear models of dynamical systems using measurement data and DMDc~\cite{DMDc}. Moreover, we explain how strategies suggested in previous works can be formulated within our theoretical framework.  Additionally, in this section, we propose a new approach for a forced linear characterization of nonlinear systems based on the introduced theoretical framework.

\subsection{Theoretical motivation}
\label{theoretical_motivation}

Let us consider measurements $\{ \mathbf{x}_k\}_{k=0}^M$ representing the evolution of a
dynamical system in the form of (\ref{dynamicalsystem}), and a collection of
observables $\{g_i\}_{i=1}^N$ through which we observe the system. The aim is to
find a forced linear representation for the temporal evolution of the vector-valued
observations $\{ \hat{\mathbf{g}}(\mathbf{x}_k)\}_{k=0}^M$ of the measurements, with
$ \hat{\mathbf{g}}=[g_1,\dots,g_N]^T$. From
(\ref{forced_representation}), we know that there exists a forced linear
representation of the dynamics of the observed measurements.
Furthermore, for each $f_i$ from (\ref{forced_representation}), $i=1,\dots,N$, we have $f_i \in \mathcal{H_D}$ where
$\mathcal{H_D} \subset \mathcal{H}$ is a defect space of
$\mathcal{H}_N:=\operatorname{span}\{g_i\}_{i=1}^N$ with dimension $D$. 

Since $\mathcal{H_D}$ is a linear
finite-dimensional subspace of $\mathcal{H}$, there exists an orthogonal basis
$\{\psi_l\}_{l=1}^D$ such that $\mathcal{H_D}:= \operatorname{span}\{\psi_l\}_{l=1}^D$. Thus,
for each $i=1,\dots,N$ there exist $\{B_{i,l}\}_{l=1}^D \subset \mathbb{C}$ such
that
\begin{equation}
  \label{forcing_expansion}
  f_i = \sum_{l=1}^{D} B_{i,l}\psi_l.
\end{equation}
Substituting (\ref{forcing_expansion}) into (\ref{forced_representation}) yields
\begin{equation}
    \label{forced_equation}
     \hat{\mathbf{g}}( \mathbf{x}_{k+1})= \mathbf{K} \hat{\mathbf{g}}( \mathbf{x}_k) + \mathbf{B}\hat{\boldsymbol{\psi}}( \mathbf{x}_k),
\end{equation}
where, additionally to (\ref{forced_representation}), 
$\hat{\boldsymbol{\psi}}:=[\psi_1,\dots,\psi_D]^T$ and $\mathbf{B} \in \mathbb{C}^{N \times D}$,
with $(\mathbf{B})_{i,j}=B_{i,j}$. 
From (\ref{forced_equation}) we observe that, given the values of a set of basis functions $\{\psi_l\}_{l=1}^D$ of the defect space
$\mathcal{H}_D$, computed along the trajectory $\{ \mathbf{x}_k\}$, we can obtain the matrices $\mathbf{K}\in
\mathbb{C}^{N\times N}$ and $\mathbf{B}\in \mathbb{C}^{N\times D}$ that best fit ($\ref{forced_equation}$) by solving the following minimization problem
\begin{align}
    \label{minimization_problem_defect}
    \begin{split}
    \mathbf{K}, \mathbf{B}&= \argmin_{\substack{\tilde{\mathbf{K}}\in \mathbb{C}^{N\times N},\\ \tilde{\mathbf{B}}\in \mathbb{C}^{N\times D}}}
    \sum_{k=0}^{M-1}\| \hat{\mathbf{g}}( \mathbf{x}_{k+1})-\tilde{\mathbf{K}} \hat{\mathbf{g}}( \mathbf{x}_k)- \tilde{\mathbf{B}}\hat{\boldsymbol{\psi}}( \mathbf{x}_k)\|_2^2\\
     &= \argmin_{\substack{\tilde{\mathbf{K}}\in \mathbb{C}^{N\times N},\\ \tilde{\mathbf{B}}\in \mathbb{C}^{N\times D}}}
     \| \hat{\mathbf{g}}(\bar{\mathbf{X}})-\tilde{\mathbf{K}} \hat{\mathbf{g}}(\mathbf{X})- \tilde{\mathbf{B}}\hat{\boldsymbol{\psi}}(\mathbf{X})\|_F^2,
    \end{split}
    \end{align}
where $\|\cdot\|_F$ is the Frobenius norm, and
\begin{equation*} 
\hat{\mathbf{g}}(\bar{\mathbf{X}})=[ \hat{\mathbf{g}}(\mathbf{x}_1),..., \hat{\mathbf{g}}(\mathbf{x}_M)], \quad
     \hat{\mathbf{g}}(\mathbf{X})=[ \hat{\mathbf{g}}(\mathbf{x}_0),..., \hat{\mathbf{g}}(\mathbf{x}_{M-1})] \quad 
\end{equation*}
and
\begin{equation*}
 \quad
    \hat{\boldsymbol{\psi}}(\mathbf{X})=[\hat{\boldsymbol{\psi}}(\mathbf{x}_0),...,\hat{\boldsymbol{\psi}}(\mathbf{x}_{M-1})].
\end{equation*}
 The minimization problem in
($\ref{minimization_problem_defect}$) can be solved via the DMDc algorithm
discussed in Subsection~\ref{numerical_algo}.
The operator $\mathbf{K}$ governs the linear part of the dynamics while
the forcing $\hat{\mathbf{f}}( \mathbf{x}_k)=\mathbf{B}\hat{\boldsymbol{\psi}}( \mathbf{x}_k)$ incorporates the
nonlinearities of the system. The operator $\mathbf K$ obtained from~(\ref{minimization_problem_defect}) is also relative to the chosen forcing coordinates. 
If the rows of $\hat{\mathbf g}(\mathbf X)$ and $\hat{\boldsymbol\psi}(\mathbf X)$ are linearly dependent, then the matrices $\mathbf K$ and $\mathbf B$ are not uniquely identifiable without an additional convention, such as the Moore--Penrose minimum-norm solution used by DMDc. 
Even when the combined regression matrix has full row rank, the fitted operator $\mathbf K$ depends on the forcing coordinates included in the regression. 
For this reason, the role of the defect space is not merely to improve the residual fit, but to define the particular splitting of the observed dynamics into a linear autonomous component and a forcing component.

Motivated by these observations, we seek a data-driven choice of forcing
functions \sloppy ${\{\psi_l\}_{l=1}^D \subset \mathcal H}$ whose span is separated from the
observable space $\mathcal H_N$ and can be used to model the part of the
evolution not represented by the linear term $\mathbf K\hat{\mathbf g}$. In
particular, we aim to compute from data the values of such functions along the trajectory
points $\{\mathbf x_k\}_{k=0}^{M-1}$ and use them as forcing coordinates in the
least-squares problem \eqref{minimization_problem_defect}. This yields a forced
linear characterization of the observed measurements, with the quality of the
representation depending on the chosen forcing space and the available samples. 
\subsection{Modelling the defect space}
\label{modelling}
We now use the concept of defect space, within the Koopman-theoretic framework, to motivate the DMDc-based approaches developed in \cite{Khodkar2019}  and \cite{PCT}. We introduce the approaches and relate them to different choices of the defect space.
The strategy developed in \cite{Khodkar2019} can be associated with the choice
of a polynomial basis for the defect space, while the method proposed in
\cite{PCT} can be associated with a defect space spanned by forcing functions
discovered from the one-step prediction error of a linear model of the dynamics
obtained through the DMD algorithm. 



The following is based on an arbitrary choice of observables $\{g_i\}_{i=1}^N$ through which we study the dynamics. However, we note that to work effectively, the methods from previous works require time-delay observables, i.e.,
$ \hat{\mathbf{g}}( \mathbf{x}_k):=[\mathbf{x}_k, \mathbf{x}_{k+1},\dots, \mathbf{x}_{k+d}]^T$ where $d$ is the length of the delay. 
\subsubsection{Defect space induced by polynomial basis functions}
  The choice of defect space associated with the method developed in
  \cite{Khodkar2019} is the one spanned by polynomial combinations of the
  observables considered to study the dynamics. Assume that we have chosen observables $\{g_i\}_{i=1}^N$, then we can consider all those functions
  $\psi_{I_1,...,I_N}$ of the form 
  \begin{equation}
    \label{equation_degree}
    \psi_{I_1,...,I_N}= g_1^{I_1}g_2^{I_2}...g_N^{I_N}, \quad 1<I_1+I_2+...+I_N \leq m,
  \end{equation}
  where $m$ is a parameter defining the maximum degree of the polynomial functions
  $\psi_{I_1,...,I_N}$ and $I_1,I_2,...,I_N \in \mathbb{N}$, where $\mathbb{N}$ includes 0. Note that, even if
$\{g_i\}_{i=1}^N$ are linearly independent, the functions
$\psi_{I_1,\dots,I_N}$ need not be linearly independent in general, and
$\operatorname{span}\{\psi_{I_1,\dots,I_N}\} \cap \operatorname{span}\{g_i\}_{i=1}^N$
need not be equal to $\{0\}$. For example, if $g_1(x)=x$ and $g_2(x)=x^2$, then $g_1$ and $g_2$ are linearly independent, but the polynomial combination $g_1^2$ coincides with $g_2$. Hence, the functions $\psi_{I_1,\dots,I_N}$ need not be linearly independent, and the intersection $\operatorname{span}\{\psi_{I_1,\dots,I_N}\}\cap \operatorname{span}\{g_i\}_{i=1}^N$ may be nontrivial.  Nevertheless, enlarging the set of linearly independent functions $\psi_{I_1,\dots,I_N}$ increases the size of the approximation space and can improve the accuracy of the forced linear representation of the dynamics.
\paragraph{Computing polynomial basis functions from data}
Given the vector-valued observed measurements $\{ \hat{\mathbf{g}}( \mathbf{x}_k)\}_{k=0}^{M-1}$, for each
  $k$ compute the vectors $\hat{\boldsymbol{\psi}}( \mathbf{x}_k)=[\psi_1( \mathbf{x}_k), \psi_2( \mathbf{x}_k), ... ,\psi_D( \mathbf{x}_k)]^T$,
  where $\psi_l$, $l=1,...,D$ are polynomial basis functions obtained from 
  the observables, i.e., 
  \begin{equation}
    \psi_l( \mathbf{x}_k)= g_1^{I_{l_1}}( \mathbf{x}_k)g_2^{I_{l_2}}( \mathbf{x}_k)...g_N^{I_{l_N}}( \mathbf{x}_k).
  \end{equation}
  The exponents $I_{l_1}, ..., I_{l_N}$ of the observables in each basis function are hyperparameters which need to be chosen a priori. Notice that the effectiveness of the method relies on the choice of these hyperparameters. 
  
  In case the underlying differential equation is not known, the authors suggest making this selection based on the Pearson correlation coefficients
  between the time series of the discrete temporal derivative of each observable $g_i$
  and the time series of all other observables. Observables with correlations higher than a certain threshold should be present in a basis function. 
  The degrees of the basis functions have to be chosen without a data-driven strategy. 
  
\subsubsection{Defect space induced by discovered forcing based on DMD's error} 
The method used in \cite{PCT} first builds a linear approximation of
the action of the Koopman operator onto the observables, i.e., given a vector-valued
observation $ \hat{\mathbf{g}}=[g_1,g_2,\dots,g_N]^T$ it finds the matrix $\mathbf{K}$ such that 
  \begin{equation}
    \mathcal{K} \hat{\mathbf{g}} \approx \mathbf{K} \hat{\mathbf{g}}.
  \end{equation}
  The functions forcing the dynamics are then computed as the step-wise error made by the linear model
  in approximating the action of the Koopman operator, i.e.
  \begin{equation}
    \hat{\boldsymbol{\psi}}=\mathcal{K} \hat{\mathbf{g}} -\mathbf{K} \hat{\mathbf{g}},
  \end{equation}
where $\hat{\boldsymbol{\psi}}:=[\psi_1,\psi_2,\dots,\psi_N]^T$. The functions $\{\psi_i\}_{i=1}^N$ do not represent a basis of a defect 
space of $\mathcal{H}_N:=\operatorname{span}\{g_i\}_{i=1}^N$, they are rather functions that belong to a defect space that can be computed
only after the linear prediction has been performed. Thus, in this method a linear model is forced with whatever functions
are needed to effectively predict the action of the Koopman operator on the observables chosen to study the dynamics. 

This method strongly relies on the DMD's ability to build a reasonable linear model.
Unfortunately, DMD presents several closure issues related to modelling nonlinear systems \cite{Wu2021}, which lessen the effectiveness and reliability of this strategy.

\paragraph{Computing forcing based on DMD's error}
Given the vector-valued observed measurements $\{ \hat{\mathbf{g}}( \mathbf{x}_k)\}_{k=0}^M$, use DMD to solve the following minimization problem
\begin{equation}
  \mathbf{K} = \argmin_{\tilde{\mathbf{K}}\in \mathbb{C}^{N \times N}}
                       \| \hat{\mathbf{g}}(\bar{\mathbf{X}})-\tilde{\mathbf{K}} \hat{\mathbf{g}}(\mathbf{X})\|_F,
\end{equation}
where $ \hat{\mathbf{g}}(\bar{\mathbf{X}})$ and $ \hat{\mathbf{g}}(\mathbf{X})$ are as in (\ref{observed_snapshots}). After that, for each $k=0,\dots,M-1$, compute the stepwise linear prediction error as
\begin{equation}
    \hat{\boldsymbol{\psi}}(\mathbf{x}_k):=\hat{\mathbf{g}}(\mathbf{x}_{k+1})-\mathbf{K}\hat{\mathbf{g}}(\mathbf{x}_{k}).
\end{equation}
\subsection{Defect space induced by a correlation operator} 
\label{correlation_basis_functions}
Our choice of basis functions for the defect space is inspired by Proper Orthogonal Decomposition \cite{Hol1996}, where a correlation matrix is used to extract orthogonal directions that optimally represent the dominant structure of the data. In our setting, we exploit a similar idea for a different purpose: rather than selecting modes for approximation alone, we use the spectral structure of a correlation operator associated with the observables to identify functions that are naturally separated from the observable subspace and can therefore serve as basis functions for the defect space. The construction below should be understood as a data-driven model for the defect contribution rather than as a proof that an exact Koopman defect space has been recovered. 

 To compute an effective set
 of basis functions for the defect space  from the observed measurements $\{ \hat{\mathbf{g}}( \mathbf{x}_k)\}_{k=0}^M$, $ \hat{\mathbf{g}}=[g_1,\dots,g_N]^T$, we start by considering the
kernel function $c:\mathcal{M}\times\mathcal{M}\rightarrow \mathbb{C}$ defined
as
\begin{equation}
  \label{Kernel_function}
    c(\mathbf{x},\mathbf{y}):= \sum_{i=1}^N g_i(\mathbf{x})\overline{g_i(\mathbf{y})},
\end{equation}
where $\mathbf{x},\mathbf{y} \in \mathcal{M}$ and $\overline{g_i(\mathbf{y})}$ is the conjugate of $g_i(\mathbf{y})$. Notice that since $\{g_i\}_{i=1}^N\subset \mathcal{H}:= L^2_\mu(\mathcal M;\mathbb C)$, we have $c \in L^2_{\mu}(\mathcal{M} \times \mathcal{M}; \mathbb{C})$.
As a matter of fact
\begin{equation}
\begin{aligned}
\int_{\mathcal M}\int_{\mathcal M}|c(\mathbf x,\mathbf y)|^2\,d\mu(\mathbf x)\,d\mu(\mathbf y)
&\le \int_{\mathcal M}\int_{\mathcal M}\sum_{i=1}^N |g_i(\mathbf x)|^2|g_i(\mathbf y)|^2\,d\mu(\mathbf x)\,d\mu(\mathbf y) \\
&\quad + \int_{\mathcal M}\int_{\mathcal M}\sum_{i\ne j}|g_i(\mathbf x)\overline{g_i(\mathbf y)}|\,|g_j(\mathbf x)\overline{g_j(\mathbf y)}|\,d\mu(\mathbf x)\,d\mu(\mathbf y) \\
&\le \sum_{i=1}^N \|g_i\|^4 + \sum_{i\ne j}\|g_i\|^2\|g_j\|^2 < \infty,
\end{aligned}
\end{equation}
where the second inequality follows from Hölder's inequality. Moreover, the kernel $c(\mathbf{x},\mathbf{y})$  satisfies the Hermitian symmetry relation:
\begin{equation}
\label{kernel}
c(\mathbf{x},\mathbf{y})
= \sum_{i=1}^N g_i(\mathbf{x})\overline{g_i(\mathbf{y})}
= \overline{\sum_{i=1}^N g_i(\mathbf{y})\overline{g_i(\mathbf{x})}}
= \overline{c(\mathbf{y},\mathbf{x})},
\qquad \mathbf{x},\mathbf{y}\in\mathcal M.
\end{equation}
We define the correlation operator $\mathcal{C}:\mathcal{H}\rightarrow
\mathcal{H}$ as:
\begin{equation}
\label{operator}
(\mathcal C h)(\mathbf{x})= \int_{\mathcal{M}} c(\mathbf{x},\mathbf{y})\,h(\mathbf{y})\,d\mu(\mathbf{y}),
\qquad \forall h\in\mathcal H.
\end{equation}
Since for the kernel $c \in L^2_{\mu}(\mathcal{M} \times \mathcal{M}; \mathbb{C})$ holds,
the operator $\mathcal{C}$ defined in ($\ref{operator}$) is a Hilbert-Schmidt
operator \cite[pg 1009]{dunford}. From
operator theory we know that every Hilbert-Schmidt operator is compact \cite[Theorem 6 pg 1012]{dunford},
and from the symmetry of the kernel we have that $\mathcal{C}$
is also self-adjoint (Appendix B). Since the correlation operator $\mathcal{C}$ is compact and self-adjoint, 
the number of eigenvalues of $\mathcal{C}$ is countable, and they are all
real. Furthermore, there exists an orthonormal basis of the separable Hilbert space $\mathcal{H}$
consisting of eigenfunctions under $\mathcal{C}$. More precisely, let us define the eigenspaces
\begin{equation}
  \mathcal{H}_{\lambda}:=\{\psi \in \mathcal{H}\;|\; \mathcal{C}\psi=\lambda \psi\},
\end{equation}
so that $\mathcal{H}_{\lambda} \setminus \{0\}$ consists of the eigenfunctions of
$\mathcal{C}$ with eigenvalue $\lambda$. Note that $0$ is the zero element of the Hilbert space $\mathcal H$, that is, the equivalence class of functions equal to zero almost everywhere with respect to $\mu$.  As a consequence of the properties of the operator $\mathcal{C}$, the spaces $\mathcal{H}_{\lambda}$ are
mutually orthogonal, each $\mathcal{H}_{\lambda}$ for $\lambda \neq 0 $ is
finite-dimensional, and there exists an orthonormal basis of $\mathcal{H}$ given by
the union of orthonormal bases of the $\mathcal{H}_{\lambda}$'s \cite[Theorem 2.3 pg 39]{Knapp2017}. Notice
that $\mathcal{H}_0$ coincides with the null space, $\nullspace(\mathcal{C})$, of the
correlation operator, which we here consider as the space of eigenfunctions
associated with eigenvalue $\lambda=0$.

Ideally, since the eigenfunctions of the correlation operator form an
orthonormal basis of $\mathcal{H}$, our aim is to consider a finite subset of
linearly independent eigenfunctions $\{\psi_l\}_{l=1}^D$ such that, for each
$g\in\mathcal{H}_N:=\operatorname{span}\{g_i\}_{i=1}^N$,
\[
\mathcal{K}g \in \operatorname{span}\{\psi_l\}_{l=1}^D+\mathcal{H}_N,
\qquad
\operatorname{span}\{\psi_l\}_{l=1}^D\cap\mathcal{H}_N=\{0\}.
\]
After that, we can use them to find a solution to the minimization problem in
$\eqref{minimization_problem_defect}$ using DMDc. Note that the space $\mathcal{H}_N$
coincides with the span of the eigenfunctions of $\mathcal{C}$ associated with
nonzero eigenvalues, that is,
\[
\mathcal H_N
=
\operatorname{span}\{\psi\in\mathcal H:\mathcal C\psi=\lambda\psi,\ \lambda\neq 0\},
\]
and hence $\nullspace(\mathcal C)=\mathcal H_N^\perp$; see Appendix B. Thus, to
ensure that
$\operatorname{span}\{\psi_l\}_{l=1}^D\cap\mathcal{H}_N=\{0\}$, we consider
eigenfunctions of $\mathcal{C}$ associated with the zero eigenvalue.

The full null space $\nullspace(\mathcal C)$ therefore provides a natural orthogonal
complement to the observable space $\mathcal H_N$, and it contains all
components of $\mathcal K\mathcal H_N$ that are orthogonal to $\mathcal H_N$.
This motivates using zero-eigenvalue eigenfunctions of $\mathcal C$ as forcing
coordinates. In computations, however, only finitely many such eigenfunctions can
be approximated from data. Therefore, we consider the computable
zero-eigenvalue eigenfunctions of $\mathcal C$ and call the sampled functions
obtained in this way the \textit{correlation basis functions}. This finite
empirical basis should be understood as a computable modelling choice for the
forcing coordinates, motivated by the exact orthogonal decomposition, rather
than as a guaranteed recovery of the exact complement needed to represent
$\mathcal K\mathcal H_N$.

\paragraph{Estimating the correlation basis functions from data}
Let $\{\hat{\mathbf g}(\mathbf x_k)\}_{k=0}^{M}\subset\mathbb C^{N}$ be the observed measurements, where $\hat{\mathbf g}=[g_1,\dots,g_N]^T$. We seek linearly independent functions in the null space of the correlation operator associated with $\{g_i\}_{i=1}^N$, that is, functions $\psi\in\mathcal H$ satisfying
\begin{equation}
\label{eigenvalue_problem}
(\mathcal C\psi)(\mathbf x)=\int_{\mathcal M} c(\mathbf x,\mathbf y)\psi(\mathbf y)\,d\mu(\mathbf y)=0.
\end{equation}
We want to estimate the values of such functions along the trajectory given 
by the measurements $\{ \mathbf{x}_k\}_{k=0}^M$. To achieve that, we replace the measure $\mu$ by the empirical measure
\begin{equation}
\label{eq:empirical_measure}
\mu_M := \frac{1}{M}\sum_{j=0}^{M-1}\delta_{\mathbf x_j},
\end{equation}
supported on the sampled points $\{\mathbf x_j\}_{j=0}^{M-1}$. We construct the correlation matrix using only the first $M$ snapshots
rather than all $M+1$ measurements. This choice is consistent with the DMDc optimization problem, which is formulated on the $M$ one-step transitions $\mathbf{x}_k \rightarrow \mathbf{x}_{k+1}$, $k=0,\dots,M-1$. Accordingly, the discrete correlation operator is constructed from the points $\{\mathbf{x}_j\}_{j=0}^{M-1}$, so that the resulting basis functions are evaluated on the same $M$ samples entering the forcing term. The final snapshot $\mathbf{x}_M$ is used only in the shifted output data and does not define an additional transition.
Then, for $k=0,\dots,M-1$, we approximate the action of $\mathcal{C}$ along the sampled trajectory by the equal-weight empirical average
\[
(\mathcal{C}\psi)(\mathbf{x}_k)
\approx \frac{1}{M}\sum_{j=0}^{M-1} c(\mathbf{x}_k,\mathbf{x}_j)\psi(\mathbf{x}_j).
\]
This corresponds to approximating the integral with respect to $\mu$ by integration against the empirical measure $\mu_M$ in (\ref{eq:empirical_measure}).
As a consequence, we have that the finite-dimensional approximation of the 
operator $\mathcal{C}$ is given by the matrix $\mathbf{C}\in \mathbb{C}^{M\times M}$ such that 
$(\mathbf{C})_{i,j}=c(\mathbf{x}_i, \mathbf{x}_j)$ for all $i,j=0,\dots,M-1$, where $c(\mathbf{x}_i, \mathbf{x}_j)$ is the kernel
function defined in (\ref{Kernel_function}), computed at points $\mathbf{x}_i,\; \mathbf{x}_j$ of the trajectory.
Thus, a discrete approximation  of (\ref{eigenvalue_problem}) is then given by
\begin{equation}
    \frac{1}{M}\sum_{j=0}^{M-1}c( \mathbf{x}_k, \mathbf{x}_j)\psi_i( \mathbf{x}_j)= 0,
\end{equation} 
where $\psi_i( \mathbf{x}_j)$ is the $j$-th entry of the $i$-th eigenvector of the 
matrix $\mathbf{C}$ associated
with zero eigenvalue.  The $l$-th zero-eigenvalue eigenvector, $\bar{\boldsymbol{\psi}}_l=[\psi_l(\mathbf x_0),\psi_l(\mathbf x_1),\dots,\psi_l(\mathbf x_{M-1})]^T$ of the matrix $\mathbf{C}$ can be regarded as a sampled zero-eigenvalue mode of the correlation operator $\mathcal C$. Accordingly, its entries are used as sampled values of a candidate forcing coordinate along the trajectory $\{\mathbf x_k\}_{k=0}^{M-1}$.

In practice, we want to estimate the null space of the matrix $\mathbf{C} \in \mathbb{C}^{M\times M}$. The estimation 
of the null space can be performed by resorting to a singular value decomposition (SVD) of $\mathbf{C}$.
Indeed, for any matrix  $\mathbf{A} \in \mathbb{C}^{m \times n}$ with $r$ non-zero singular
values, the null space of $\mathbf{A}$ is spanned by the right singular vectors corresponding
to the $n$-$r$ zero singular values of $\mathbf{A}$ \cite{Horn2012}. 

$\mathbf{C}\in\mathbb{C}^{M\times M}$ is Hermitian positive semidefinite, therefore, its singular value decomposition is 
\[
\mathbf{C}=\mathbf{\Psi}\mathbf{\Sigma}\mathbf{\Psi}^*,
\]
where $\mathbf{\Psi}$ is a unitary matrix whose columns are the orthonormal eigenvectors of $\mathbf{C}$, $\mathbf{\Psi}^*$ its conjugate transpose, and $\mathbf{\Sigma}=\operatorname{diag}(\sigma_1,\dots,\sigma_M)$, with $\sigma_i\ge 0$ are both the eigenvalues and singular values of $\mathbf{C}$. 
We are only interested in those eigenvectors associated with
zero singular values, that is, those eigenvectors in the null space of $\mathbf{C}$. 

To estimate the null space in finite precision, one should not rely on exact
zero singular values. Because of numerical roundoff, sampling effects, and
ill-conditioning, directions that analytically belong to the null space may
correspond to computed small but nonzero singular values. We therefore define a
\emph{numerical null space} by introducing a relative threshold $\alpha$ for the
singular values $\sigma_j$ of $\mathbf C$. A default choice for $\alpha$ is
machine precision, i.e., $\alpha\approx 2.22\times 10^{-16}$ in double
precision. With this choice, our criterion agrees with the thresholding used in
SciPy's \texttt{null\_space} function \cite{scipy}, which treats singular values
satisfying
\begin{equation}
\label{eq:alpha}
   \sigma_j \leq \alpha M\max_i{\sigma_i}
\end{equation}
as numerically zero.

The choice of $\alpha$ is therefore part of the numerical method. 
In the experiments below we report the value of $\alpha$ used in each case. A principled, problem-independent choice of this threshold is  out of scope of this paper and considered to be future work.

\section{The proposed numerical procedure: Correlation-Basis enhanced DMD (CB-DMD)}
This section introduces our numerical procedure to compute forced linear models of nonlinear dynamics based on our choice of the basis functions of the defect space and DMDc. The steps of the procedure are summarized in Algorithm \ref{alg:Method}.  We recall that our final goal is to compute linear operators that are less prone to spurious modes by capturing in the forcing term the nonlinear components of the dynamics.
\subsection{Numerical procedure of CB-DMD}
  \begin{algorithm}[t]
    \caption{Correlation-Basis DMD (CB-DMD)}\label{alg:Method}
    \textbf{Input} Observed measurements $\{ \hat{\mathbf{g}}( \mathbf{x}_k)\}_{k=0}^M \subset \mathbb{C}^N$\\
    \textbf{Output} Linear operator $\mathbf{K} \in \mathbb{C}^{N \times N}$, forcing matrix $\mathbf{B}\Psi(\mathbf{X}) \in \mathbb{C}^{N \times M}$, and reconstructed 
    observed measurements $\{\tilde{\mathbf{g}}( \mathbf{x}_k)\}_{k=1}^M$.
    \begin{algorithmic}[1]
  
      \STATE Arrange the data into matrices $ \hat{\mathbf{g}}(\bar{\mathbf{X}})$ and $ \hat{\mathbf{g}}(\mathbf{X})$ as in
      (\ref{observed_snapshots}).
      \ Compute the
      correlation matrix $\mathbf{C} \in \mathbb{C}^{M \times M}$ as follows:
      \begin{equation}
        \mathbf{C}:= \frac{1}{M} \hat{\mathbf{g}}(\mathbf{X})^T\overline{ \hat{\mathbf{g}}(\mathbf{X})},
      \end{equation}
      \STATE Compute the singular value decomposition (SVD) of $\mathbf{C}$,
      \begin{equation}
        \mathbf{C}= \mathbf{\Psi\Sigma\Psi}^*,
      \end{equation} 
     \STATE Let $r$ denote the number of singular values regarded as nonzero according to the thresholding criterion in (\ref{eq:alpha}), and define $D:=M-r.$
    For each $i=r+1,\dots,M$, let $\bar{\boldsymbol{\psi}}_i$ be the $i$-th column of $\mathbf{\Psi}$, and define
    \[
    \mathbf{\Psi}(\mathbf{X})
    :=[\bar{\boldsymbol{\psi}}_{r+1},\bar{\boldsymbol{\psi}}_{r+2},\dots,\bar{\boldsymbol{\psi}}_M]^T
    \in\mathbb{C}^{D\times M}.
    \]    
      \STATE Use DMDc to solve the following minimization
      problem in the Frobenius norm:
      \begin{align}
        \mathbf{K},\mathbf{B} &= \argmin_{\substack{\tilde{\mathbf{K}}\in \mathbb{C}^{N \times N},\\ \tilde{\mathbf{B}}\in \mathbb{C}^{N \times D}}}
                             \| \hat{\mathbf{g}}(\bar{\mathbf{X}})-\tilde{\mathbf{K}} \hat{\mathbf{g}}(\mathbf{X})- \tilde{\mathbf{B}}\mathbf{\Psi}(\mathbf{X})\|_F.
      \end{align}
      \STATE Given an initial condition, $ \hat{\mathbf{g}}(\mathbf{x}_0)$, use the learned linear operator to predict the dynamics trajectory as follows:
      \begin{equation}
        \tilde{\mathbf{g}}( \mathbf{x}_{k+1})=\mathbf{K}\tilde{\mathbf{g}}(\mathbf{x}_{k})
      \end{equation}
     
  \end{algorithmic}
  \end{algorithm}

\paragraph{Step 1} 
Gather the measurements representing a discrete time
realization of the dynamical system we analyze, choosing a set of observables of
interest through which we want to study the dynamics, and arrange the data
into matrices.
Specifically, given measurements $\{ \mathbf{x}_k\}_{k=0}^M$ representing
the evolution of a dynamical system, we can choose any finite set of
observables $\{g_i\}_{i=1}^N \subset \mathcal{H}$ and study the dynamics of the
observed measurements $\{ \hat{\mathbf{g}}( \mathbf{x}_k)\}_{k=0}^M$, where
$ \hat{\mathbf{g}}=[g_1,\dots,g_N]^T$.  In this work, we mainly focus on analysing the
dynamics in the original state variables, i.e., $ \hat{\mathbf{g}}( \mathbf{x}_k)= \mathbf{x}_k$.
However, the procedure applies for any given choice of observables
$\{g_i\}_{i=1}^N$. Once the set of observables
through which we want to study the dynamics has been chosen, we build the matrices
\begin{equation}
  \label{observed_snapshots}
   \hat{\mathbf{g}}(\mathbf{X})=[ \hat{\mathbf{g}}(\mathbf{x}_0), \hat{\mathbf{g}}(\mathbf{x}_1),\dots, \hat{\mathbf{g}}(\mathbf{x}_{M-1})] \quad \text{and} \quad  \hat{\mathbf{g}}(\bar{\mathbf{X}})=[ \hat{\mathbf{g}}(\mathbf{x}_1), \hat{\mathbf{g}}(\mathbf{x}_2),\dots, \hat{\mathbf{g}}(\mathbf{x}_M)],  
\end{equation} 
 where each column represents an observed state of the dynamical system at
  different times. The columns are arranged chronologically, and the columns of
  $ \hat{\mathbf{g}}(\bar{\mathbf{X}})$ are shifted one time step forward in the future with
  respect to the columns of $ \hat{\mathbf{g}}(\mathbf{X})$.
\paragraph{Step 2}
Compute the numerical approximation of the correlation 
operator $\mathcal{C}$.
Given the matrix $ \hat{\mathbf{g}}(\mathbf{X})$ defined in ($\ref{observed_snapshots}$), compute the
correlation matrix $\mathbf{C} \in \mathbb{C}^{M \times M}$ as follows
\begin{equation}
  \mathbf{C}:= \frac{1}{M} \hat{\mathbf{g}}(\mathbf{X})^T\overline{ \hat{\mathbf{g}}(\mathbf{X})},
\end{equation}
where $\overline{ \hat{\mathbf{g}}(\mathbf{X})}$ is the conjugate of $ \hat{\mathbf{g}}(\mathbf{X})$.
\paragraph{Step 3}
Compute the values of the correlation basis functions along the system trajectory
given by the measurements $\{\mathbf{x}_k\}_{k=0}^M$. Given the matrix
$\mathbf{C}\in\mathbb{C}^{M\times M}$ from the previous step, compute its singular value decomposition
\begin{equation}
  \mathbf{C}=\mathbf{\Psi}\mathbf{\Sigma}\mathbf{\Psi}^*,
\end{equation}
where $\mathbf{\Psi}\in\mathbb{C}^{M\times M}$ is a unitary matrix whose columns are the eigenvectors of $\mathbf{C}$,
$\mathbf{\Psi}^*$ is its conjugate transpose, and
$ \mathbf{\Sigma}=\operatorname{diag}(\sigma_1,\dots,\sigma_M)$
contains the singular values of $\mathbf{C}$. Let $r$ denote the number of singular values regarded as nonzero according to the thresholding criterion in (\ref{eq:alpha}), and define $D:=M-r.$
For each $i=r+1,\dots,M$, let
\begin{equation}
  \bar{\boldsymbol{\psi}}_i=
  [\psi_i(\mathbf{x}_0),\psi_i(\mathbf{x}_1),\dots,\psi_i(\mathbf{x}_{M-1})]^T
  \in\mathbb{C}^M
\end{equation}
be the $i$-th column of $\mathbf{\Psi}$. These vectors are associated with the numerical null space of $\mathbf{C}$ and can be regarded as discrete approximations of zero-eigenvalue modes of the correlation operator evaluated along the sampled trajectory. We then define
\begin{equation}
  \mathbf{\Psi}(\mathbf{X})
  :=
  [\bar{\boldsymbol{\psi}}_{r+1},\bar{\boldsymbol{\psi}}_{r+2},\dots,\bar{\boldsymbol{\psi}}_M]^T
  \in\mathbb{C}^{D\times M}.
\end{equation}
Thus, the $k$-th column of $\mathbf{\Psi}(\mathbf{X})$ collects the values of the $D$ selected correlation basis functions at the sample point $\mathbf{x}_k$.
\paragraph{Step 4}
Use DMDc to identify the forced linear model on the sampled trajectory by solving
\begin{align}
  \label{minimization_K_B}
  \mathbf{K}, \mathbf{B}
  &=
  \argmin_{\substack{\tilde{\mathbf{K}}\in\mathbb{C}^{N\times N},\\
  \tilde{\mathbf{B}}\in\mathbb{C}^{N\times D}}}
  \left\|
  \hat{\mathbf{g}}(\bar{\mathbf{X}})
  -\tilde{\mathbf{K}}\hat{\mathbf{g}}(\mathbf{X})
  -\tilde{\mathbf{B}}\mathbf{\Psi}(\mathbf{X})
  \right\|_F.
\end{align}
This yields a matrix $\mathbf{K}\in\mathbb{C}^{N\times N}$ describing the linear component of the dynamics and a matrix $\mathbf{B}\in\mathbb{C}^{N\times D}$ mapping the correlation basis functions into the forcing term. Equivalently, for each $k=0,\dots,M-1$, the fitted model has the form
\begin{equation}
  \hat{\mathbf{g}}(\mathbf{x}_{k+1})
  \approx
  \mathbf{K}\hat{\mathbf{g}}(\mathbf{x}_k)
  +\mathbf{B}\hat{\boldsymbol{\psi}}(\mathbf{x}_k),
\end{equation}
where $\hat{\boldsymbol{\psi}}(\mathbf{x}_k)\in\mathbb{C}^D$ is the $k$-th column of $\mathbf{\Psi}(\mathbf{X})$.
\paragraph{Step 5}
Given an initial condition $\hat{\mathbf{g}}(\mathbf{x}_0)$, compute the predicted trajectory recursively using the learned autonomous linear operator $\mathbf{K}$:
\begin{equation}
\tilde{\mathbf{g}}(\mathbf{x}_{k+1})=\mathbf{K}\tilde{\mathbf{g}}(\mathbf{x}_k),
\qquad k=0,\dots,M-1,
\end{equation}
with
$\tilde{\mathbf{g}}(\mathbf{x}_0)=\hat{\mathbf{g}}(\mathbf{x}_0).$
This step provides an autonomous prediction of the dynamics based solely on the operator $\mathbf{K}$. The forcing term is used during training to separate residual components from the fitted linear part. It is not extrapolated unless an additional model for the forcing coordinates is supplied.

\subsection{Exact and numerical defect spaces in CB-DMD}
We now clarify how the correla\-tion-based forcing coordinates used in CB-DMD relate to the defect-space framework and how they affect the linear operator $\mathbf K$ learned via DMDc. At the continuous level, the correlation operator $\mathcal C$ associated with the observables $\{g_i\}_{i=1}^N$ satisfies $\nullspace(\mathcal C)=\mathcal H_N^\perp$ under the linear independence assumptions of Section~\ref{correlation_basis_functions}. Thus, $\nullspace(\mathcal C)$ provides a natural defect-space model: it is orthogonal to the span of the chosen observables and represents directions not resolved by their correlation structure.

CB-DMD replaces this continuous decomposition by a finite-sample analogue. Let
\[
\mathbf G:=\hat{\mathbf g}(\mathbf X),\qquad
\mathbf Y:=\hat{\mathbf g}(\bar{\mathbf X}),\qquad
\mathbf C:=\frac{1}{M}\mathbf G^T\overline{\mathbf G}.
\]
If $\mathbf Q\in\mathbb C^{M\times D}$ contains the selected correlation-basis vectors and $\mathbf U:=\mathbf\Psi(\mathbf X)=\mathbf Q^T$, then the augmented regression problem is
\[
\min_{\mathbf K,\mathbf B}
\|\mathbf Y-\mathbf K\mathbf G-\mathbf B\mathbf U\|_F .
\]
The selected columns of $\mathbf Q$ are those singular directions of $\mathbf C$ that are treated as zero according to the relative threshold $\alpha$ introduced in~\cref{eq:alpha}. Hence the forcing space is an empirical numerical null space, not an exact analytic one. The parameter $\alpha$ controls which null or near-null directions are included: increasing $\alpha$ enlarges the forcing space, while decreasing it makes the construction closer to the exact empirical null space. These selected directions serve as empirical forcing coordinates motivated by the continuous defect space $\nullspace(\mathcal C)$.

Let us consider the exact empirical null space case. Suppose that the columns of $\mathbf Q$ are orthonormal and satisfy $\mathbf C\mathbf Q= \mathbf 0$. Since $\mathbf C=M^{-1}\overline{\mathbf G}^{\,*}\overline{\mathbf G}$, we obtain
\[
\mathbf 0=\mathbf Q^*\mathbf C\mathbf Q
  =\frac1M(\overline{\mathbf G}\mathbf Q)^*
  (\overline{\mathbf G}\mathbf Q),
\]
and hence $\overline{\mathbf G}\mathbf Q=\mathbf 0$. Taking complex conjugates gives $\mathbf G\overline{\mathbf Q}=\mathbf 0$. Because $\mathbf U=\mathbf Q^T$, we have $\mathbf U^*=\overline{\mathbf Q}$, and consequently $\mathbf G\mathbf U^*=\mathbf 0$. Let $\mathbf R:=\mathbf Y-\mathbf K\mathbf G-\mathbf B\mathbf U$ denote the residual of the augmented least-squares problem. The normal equation associated with variations in $\mathbf K$ is $\mathbf R\mathbf G^*=\mathbf 0$, or equivalently
\[
(\mathbf Y-\mathbf K\mathbf G-\mathbf B\mathbf U)\mathbf G^*=0 .
\]
Using $\mathbf U\mathbf G^*=\mathbf 0$, this reduces to
\[
\mathbf K\mathbf G\mathbf G^*=\mathbf Y\mathbf G^* .
\]
This is precisely the normal equation for the standard DMD problem $\min_{\mathbf K}\|\mathbf Y-\mathbf K\mathbf G\|_F$. Thus, among the least-squares minimizers, the minimum-Frobenius-norm solution is $\mathbf K=\mathbf Y\mathbf G^\dagger$, so the fitted linear operator coincides with the standard DMD operator computed from the same snapshot data.

Therefore, exact empirical null space forcing does not alter the autonomous DMD operator. Changes in $\mathbf K$ can arise when $\alpha$ selects numerical near-null directions, when rank truncations or other thresholding decisions are introduced, or when finite-precision effects produce a small but nonzero coupling $\mathbf G\mathbf U^*$. In such cases, the augmented normal equations remain coupled:
\[
\mathbf K\mathbf G\mathbf G^*
+
\mathbf B\mathbf U\mathbf G^*
=
\mathbf Y\mathbf G^*,
\qquad
\mathbf K\mathbf G\mathbf U^*
+
\mathbf B\mathbf U\mathbf U^*
=
\mathbf Y\mathbf U^* .
\]
Thus, CB-DMD should be interpreted as a finite-data, $\alpha$-dependent residual-splitting procedure. The value of $\alpha$ is an essential numerical parameter and is reported in the experiments. We note that this is a theoretical analysis performed for completeness. In our experiments, we never observe an empirical null space with eigenvalues whose moduli are exactly zero.

\section{Numerical results}
\label{results}
To demonstrate the effectiveness of CB-DMD, we consider three nonlinear dynamical systems. The Nonlinear Schroedinger (NLS) equation and the post-transient low-Reynolds-number cylinder wake are chosen because they fit our setting of sustained periodic dynamics. We also include the Cubic-quintic Ginzburg-Landau (CQGL) equation as a more challenging test case, since it can exhibit ``a wide range of nonlinear spatio-temporal dynamics, including spatio-temporal periodicity and chaos''~\cite{kooppde}.

Periodic dynamics provide a useful diagnostic for detecting clearly spurious eigenvalues. For observables restricted to an invariant periodic orbit, Koopman eigenvalues associated with the sustained oscillation lie on the unit circle. Thus eigenvalues with modulus substantially larger than one indicate artificial growth in autonomous reconstructions, while strongly decaying eigenvalues are not part of the sustained oscillatory component. In finite data, complex eigenvalues inside the unit circle may also reflect transients, noise, or numerical regularization. We therefore use distance from the unit circle as a practical diagnostic rather than as an absolute mathematical criterion. Similar diagnostics are used in \cite{kooppde} for the NLS equation and in \cite{resDMD} for cylinder wake flow.

For CB-DMD, the automatic selection of eigenvectors in the null space of the correlation operator, based on the relative threshold in \eqref{eq:alpha} with $\alpha$ set to machine precision, works well for the NLS equation and the cylinder wake flow, but is too restrictive for the more challenging scenario provided by the CQGL equation. For this reason, we report the CQGL results using $\alpha=10^{-8}$ instead. For completeness, the appendix also includes the NLS and cylinder wake results for this choice of $\alpha$, as well as the CQGL results obtained with the default threshold. 

\subsection{Baselines}
We compare our proposed approach with standard DMD (no forcing) and with the methods of \cite{Khodkar2019} (polynomial forcing) and \cite{PCT} (error forcing), each corresponding to a particular choice of defect space.  Although the approaches in \cite{Khodkar2019} and \cite{PCT} can effectively capture the dynamics when time-delay observables are used, here we focus on methods that extract meaningful dynamical patterns directly from the measured observables. We therefore assess their effectiveness without relying on time-delay embeddings. 

Since the full baseline operators often produce diverging predictions, we also consider improved versions obtained by removing all DMD modes associated with eigenvalues of modulus greater than one, which may represent growing components of the dynamics.  For linear error forcing, this modification is applied both to the initial operator computed by standard DMD and to the final operator computed by DMDc (Error forcing improved). That is, the first filtering removes unstable modes before constructing the error-forcing term, while the second filtering removes unstable modes from the final operator. We note that retaining modes with modulus only slightly above one may still be beneficial in practice, since the corresponding growth can be small enough to improve numerical reconstruction over the time horizon of interest. In  Appendix D, we therefore also determine, for each method, the best eigenvalue-magnitude threshold from the set
$ \{0.99,\; 0.999,\; 0.9999,\; 1.0,\; 1.0001,\; 1.001,\; 1.01\},$ and compare the corresponding results, which align with those reported in this section.

As in \cite{Khodkar2019}, the polynomial basis functions used in polynomial forcing are selected from Pearson correlation coefficients (PCCs) between the observables and their time derivatives. Here, $\dot g_i$ denotes the time derivative of the observable $g_i$, which we approximate from the data by finite differences in time. For each observable $g_i$, we first identify the measurements $g_j$, $j\neq i$, whose PCC with $\dot g_i$ exceeds a prescribed threshold. We then form all polynomial terms, up to a chosen degree, from products of $g_i$ with these selected measurements, including the corresponding monomials. To keep the number of forcing vectors computationally manageable, additional restrictions are imposed in each subsection, and the polynomial degree is limited to two.

In our experiments, we do not truncate the singular value decompositions used in the DMD variants, but instead compute the full linear operators. This avoids confounding the effect of the defect space with that of SVD truncation, which can itself reduce the number of spurious eigenvalues. Moreover, since we are not aware of a principled approach to choose a truncation threshold, including truncation would make a thorough comparison more difficult.

\subsection{Nonlinear Schroedinger equation}
\label{nls}
The NLS equation is an important PDE describing wave propagation. It's applications include describing the evolution of surface gravity water waves or the propagation of a heat pulse in a solid \cite{NLS_app}.
The evolution equation is given by
\begin{align} iq_t + \frac{1}{2}q_{xx} + |q|^2 q = 0,
\end{align}
where $q=q(x, t)$ denotes a function in space and time. 
We set the initial conditions to $$q(x, 0) = 2 \sech(x),$$
where $\sech(x):= \frac{2}{e^{x} + e^{-x}}$ denotes the hyperbolic secant, which leads to the periodic 2-solition solution \cite{kooppde}. 
The space domain $[-15, 15]$ is discretized with $n=512$ points. 
We employ the explicit Runge-Kutta (4,5)-formula after Fourier transforming the space variable to integrate the solution. 
For computing the linear operators, we use $M = 1000$ slices comprising the time interval $t \in[0,\pi]$. The prediction interval is given by the subsequent $500$ time steps up to time $\frac{3}{2} \pi$. 

 \begin{figure}
    \centering
      \includegraphics[width=0.30\textwidth]{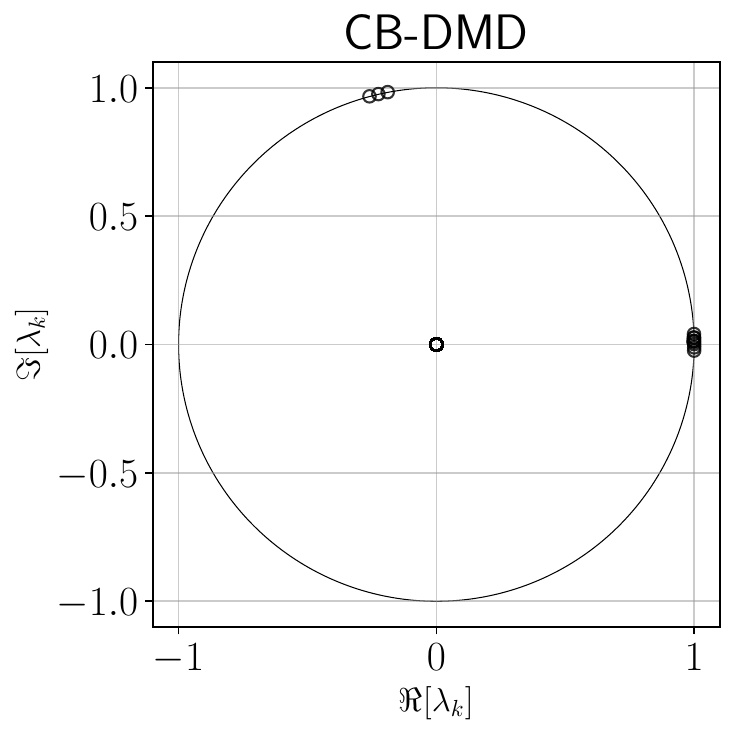}
      \includegraphics[width=0.30\textwidth]{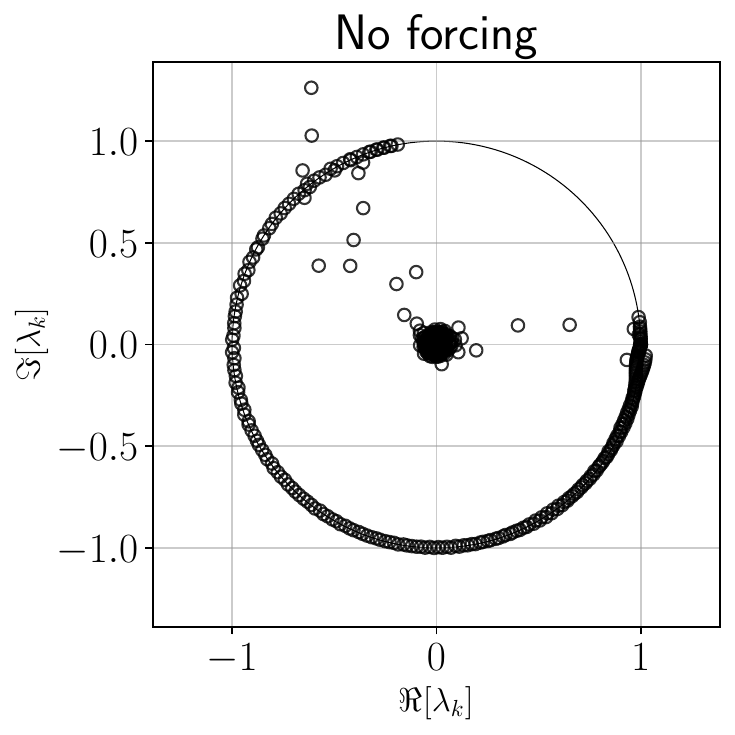}\\
     \includegraphics[width=0.30\textwidth]{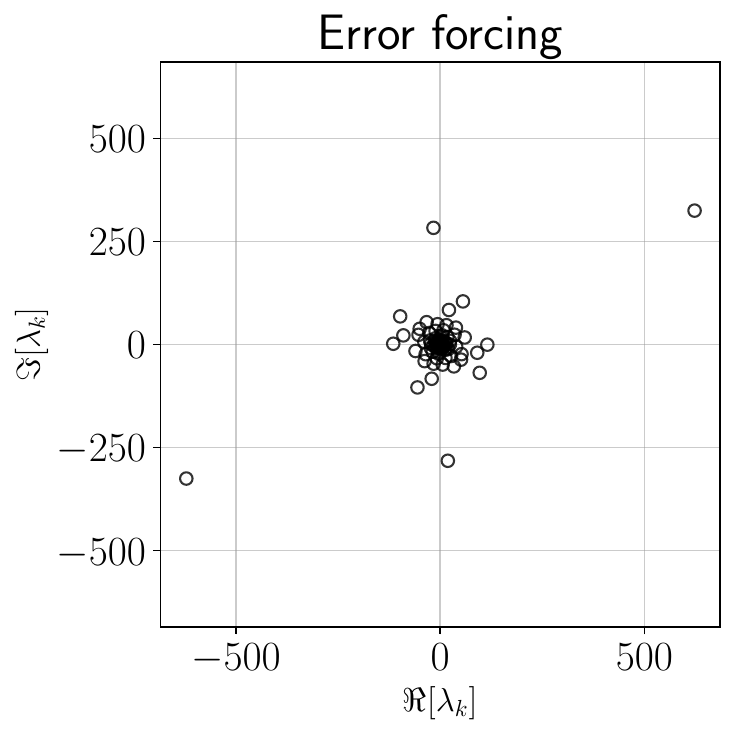}
     \includegraphics[width=0.30\textwidth]{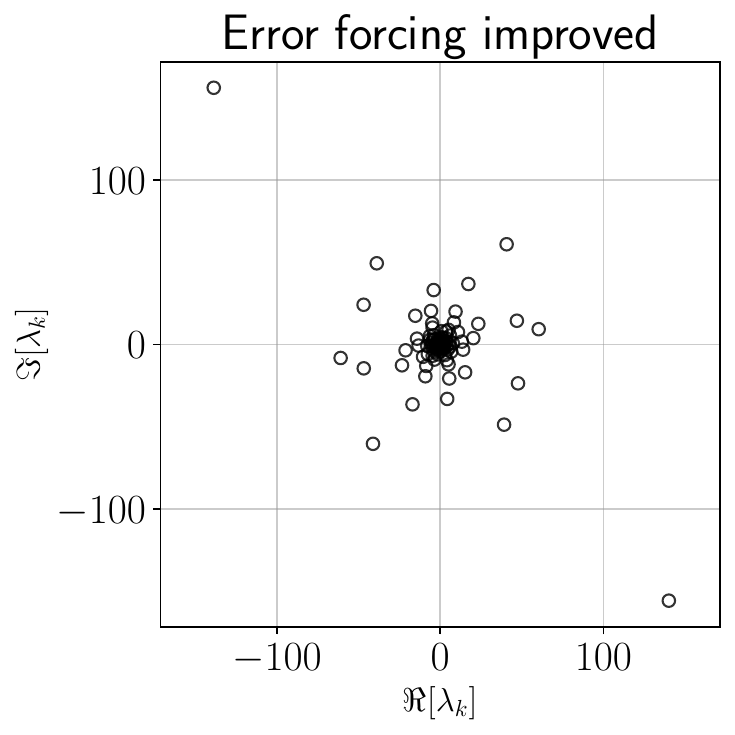}
      \includegraphics[width=0.30\textwidth]{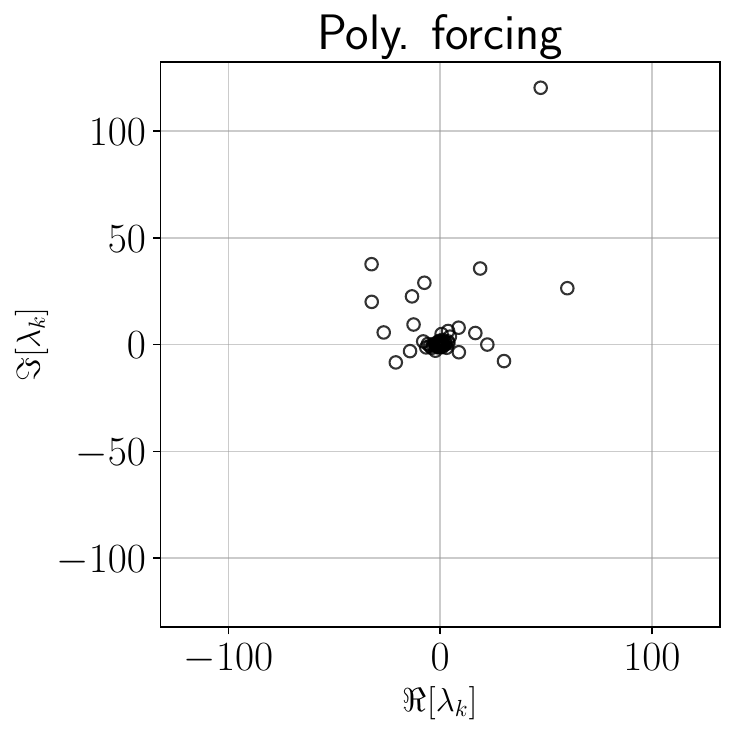}
  \caption{Spectra of the operators learned for the NLS dynamics using standard DMD, polynomial forcing, error forcing, and CB-DMD.  Improved variants are obtained by removing DMD modes associated with eigenvalues of modulus greater than one; for error forcing, this is done both for the initial standard-DMD operator and for the final DMDc operator. Since the dynamics are periodic, relevant eigenvalues are expected near the unit circle, while eigenvalues away from it are interpreted as spurious. The circular line indicates the unit circle. CB-DMD yields the most concentrated spectrum, with only $12$ numerically relevant eigenvalues ($|\lambda|>10^{-3}$), all close to the unit circle, that is, satisfying $0.99\leq|\lambda|\leq 1.01$.}
    \label{fig:NLS_spectrum}
  \end{figure}
For poly. forcing, the outer 106 grid points on both ends of the one dimensional spatial grid are excluded from the PCC computation since they are  close to zero during the entire time period of interest. A reasonable PCC threshold is visually determined to be 0.3.

\begin{table}
\centering
\caption{Comparison of eigenvalue properties across DMD methods for the NLS system. We consider an eigenvalue relevant if $|\lambda|>10^{-3}$. We also report eigenvalues close to the unit circle, defined by $0.99<|\lambda|\leq 1.01$, and eigenvalues outside the unit circle, defined by $|\lambda|>1$.}
\label{tab:eigenvalue_comparison_NLS}
\begin{tabular}{lcccc}
\hline
\textbf{Method} 
& \makecell{\textbf{Relevant}\\\textbf{eigenvalues}} 
& \makecell{\textbf{Near}\\\textbf{unit circle}} 
& \makecell{\textbf{Eigenvalues with}\\\boldmath{$|\lambda|>1$}} 
& \makecell{\textbf{Largest}\\\boldmath{$|\lambda|$}} \\
\hline
CB-DMD                        & 12  & 12 & 4   & 1.0010   \\
Standard DMD                  & 512 & 235 & 77  & 1.4027   \\
Error Forcing                 & 512 & 102 & 218 & 702.5536 \\
Error Forcing (Improved)      & 512 & 86 & 163 & 209.3115 \\
Polynomial Forcing            & 241 & 16 & 72  & 129.1700 \\
\hline
\end{tabular}
\end{table}
\begin{figure}
  \centering
    \includegraphics[width=0.7\textwidth]{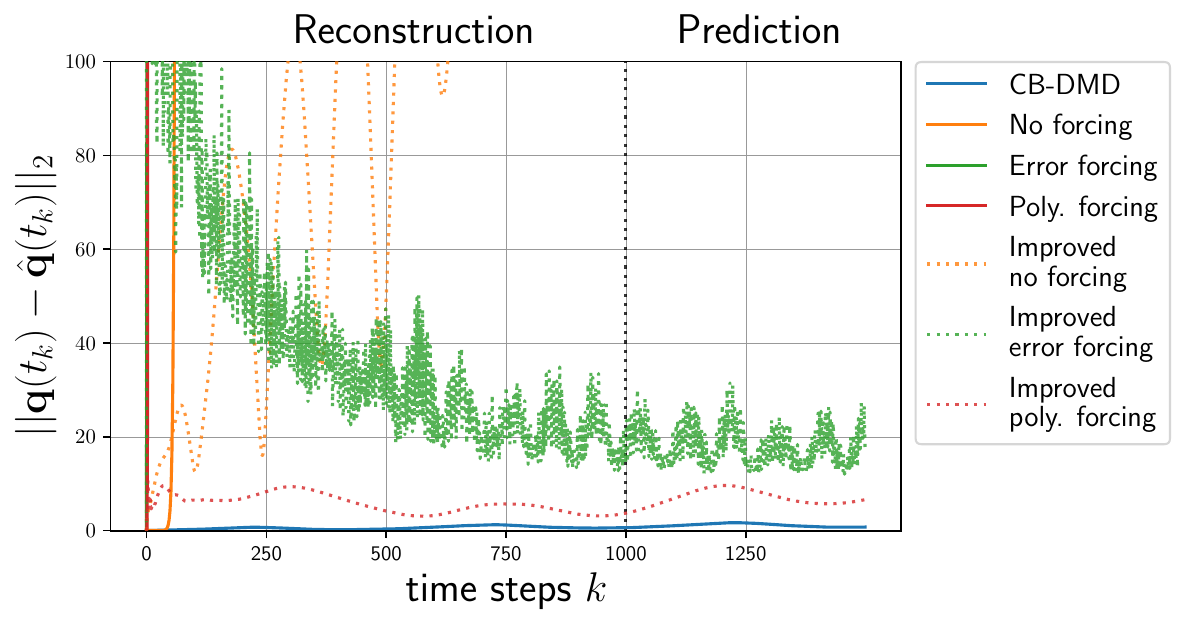}
\caption{Absolute error of the reconstructions and forward predictions for the NLS dynamics obtained with the linear operators learned by CB-DMD, the baseline methods, and their improved variants. The first $1000$ time steps correspond to reconstruction, and the following $500$ time steps correspond to prediction. The vertical dashed black line marks the start of the prediction period.}
  \label{fig:NLS_error}
  \end{figure}
The spectral results in Fig.~\ref{fig:NLS_spectrum} show a clear advantage for CB-DMD. Since the dynamics are periodic, the eigenvalues of the learned linear operator, which matter for the dynamics, are expected to lie on or very near the unit circle. Eigenvalues away from the unit circle are therefore interpreted as spurious. In particular, eigenvalues with modulus greater than one can indicate artificial growth in the reconstructed dynamics. The main feature of the CB-DMD spectrum is that it contains only a small number of numerically relevant eigenvalues, and that these eigenvalues are all located close to the unit circle.

More precisely, let us call an eigenvalue relevant if its modulus satisfies $|\lambda|>10^{-3}$. This threshold serves as a qualitative cutoff between spectral components whose magnitude indicates a meaningful contribution to the dynamics and those whose effect is negligible. In Table~\ref{tab:eigenvalue_comparison_NLS} we summarize the behaviour of the eigenvalues.
CB-DMD has only $12$ relevant eigenvalues, and all $12$ satisfy $0.99\leq|\lambda|\leq 1.01$, that is, they are close to the unit circle. The largest eigenvalue modulus for CB-DMD is $1.0010$, indicating only a negligible deviation from the unit circle and therefore suggesting numerical error rather than genuine growth in the dynamics.
Overall, CB-DMD yields a spectrum with few relevant modes, all concentrated near the unit circle, which indicates a much cleaner separation between the true dynamics and spurious effects.
All other approaches have significantly more relevant eigenvalues, and much larger eigenvalue moduli.
\begin{figure}
    \centering
    \includegraphics[width =0.32\textwidth]{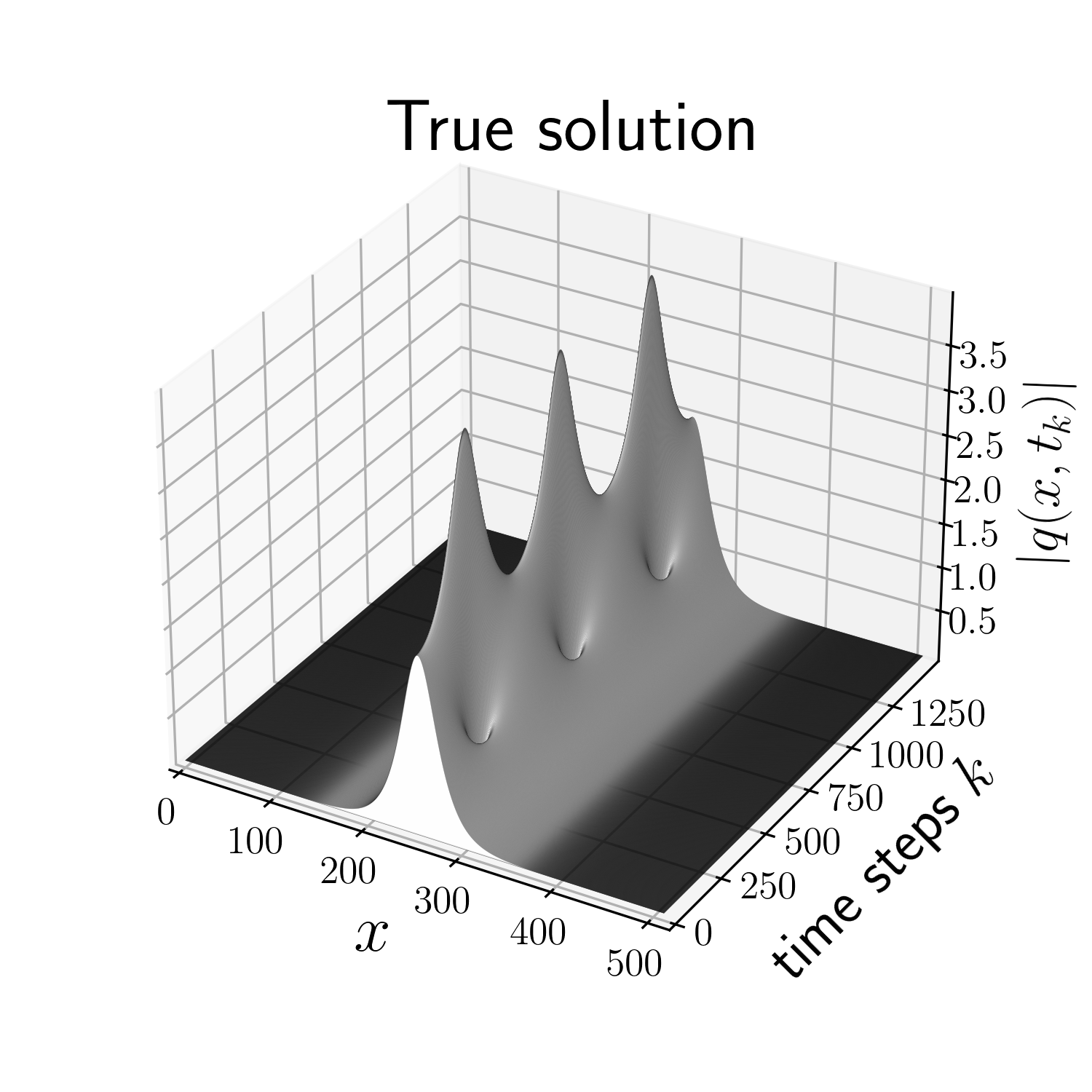}
    \includegraphics[width =0.32\textwidth]{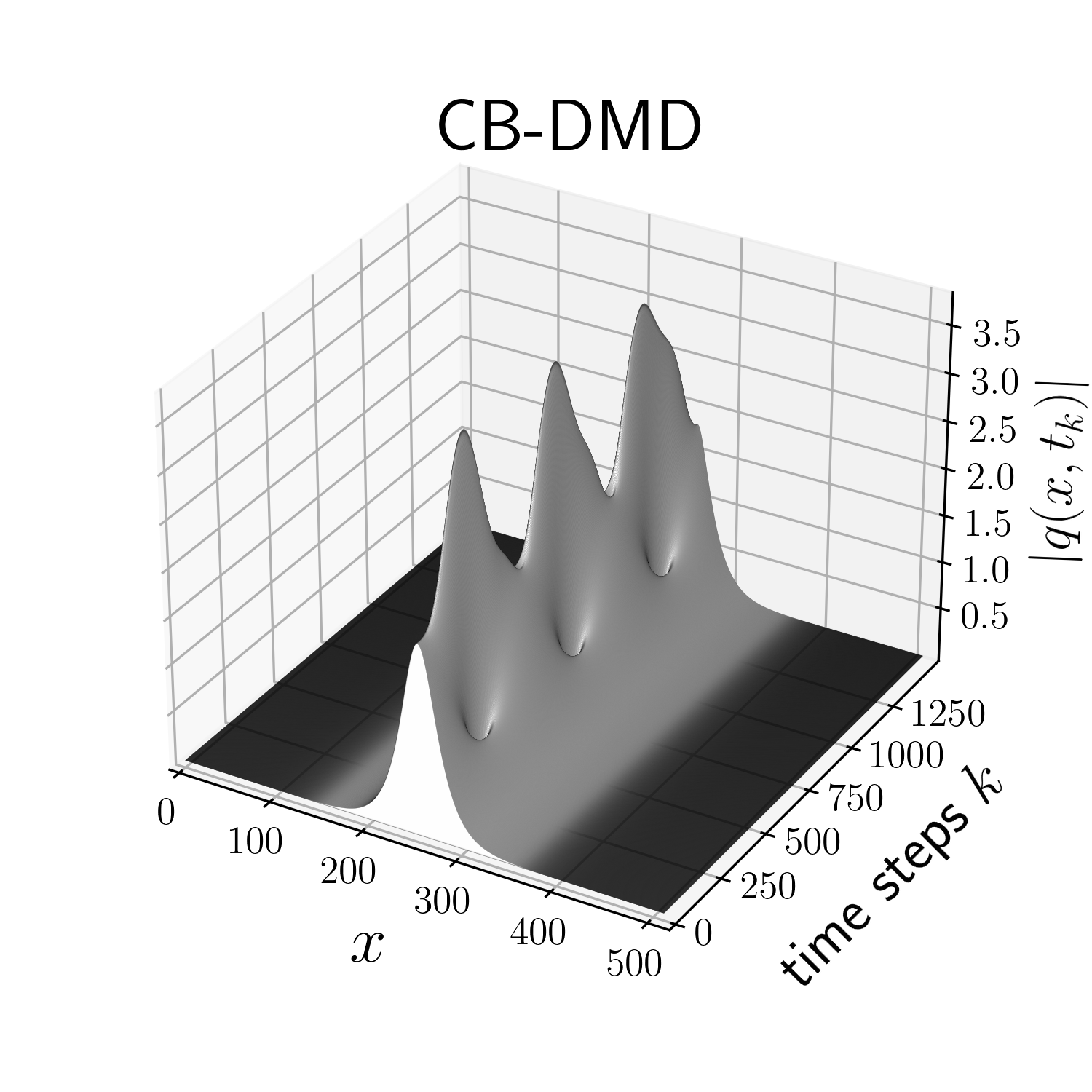}
    \includegraphics[width =0.32\textwidth]{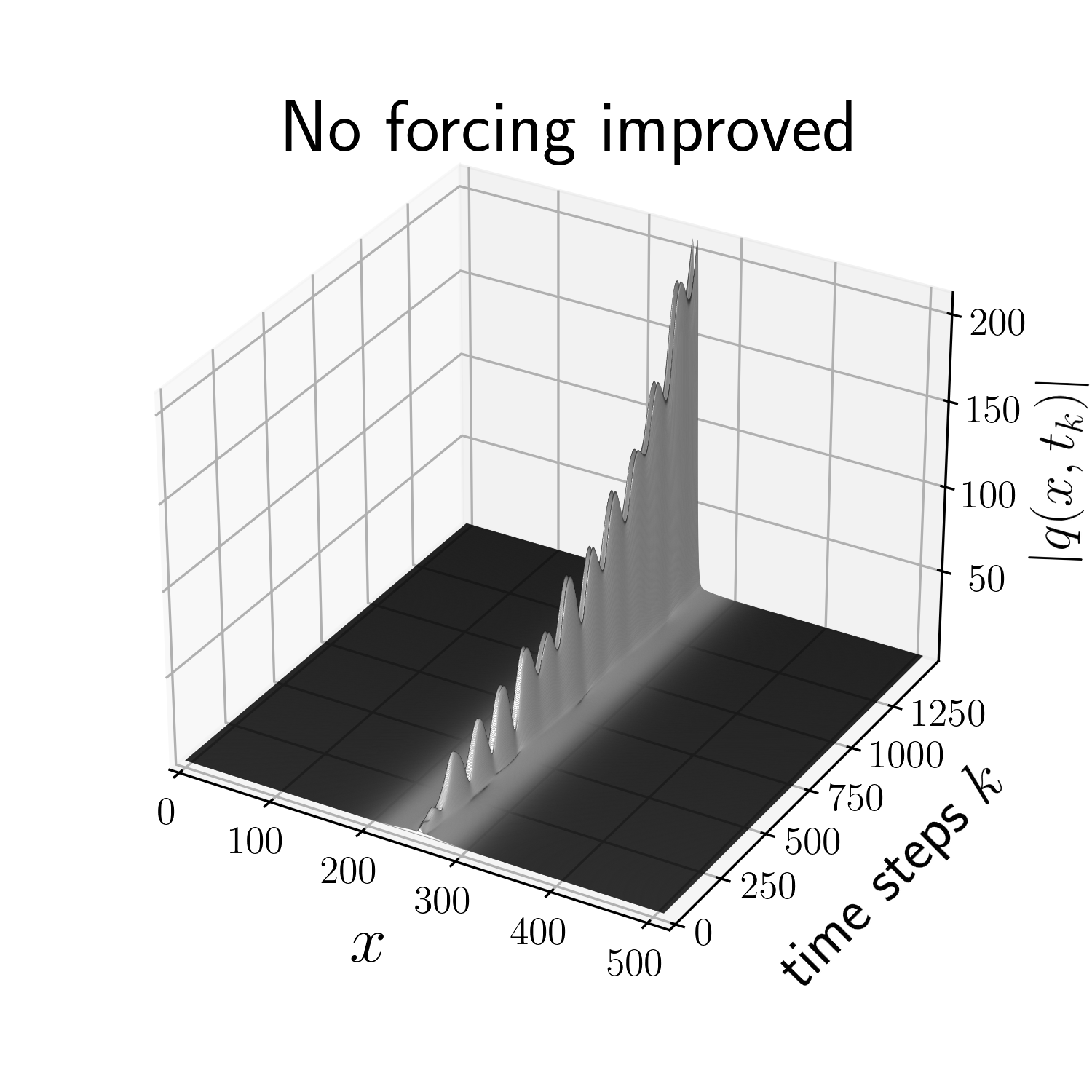}\\
    \includegraphics[height=0.32\textwidth]{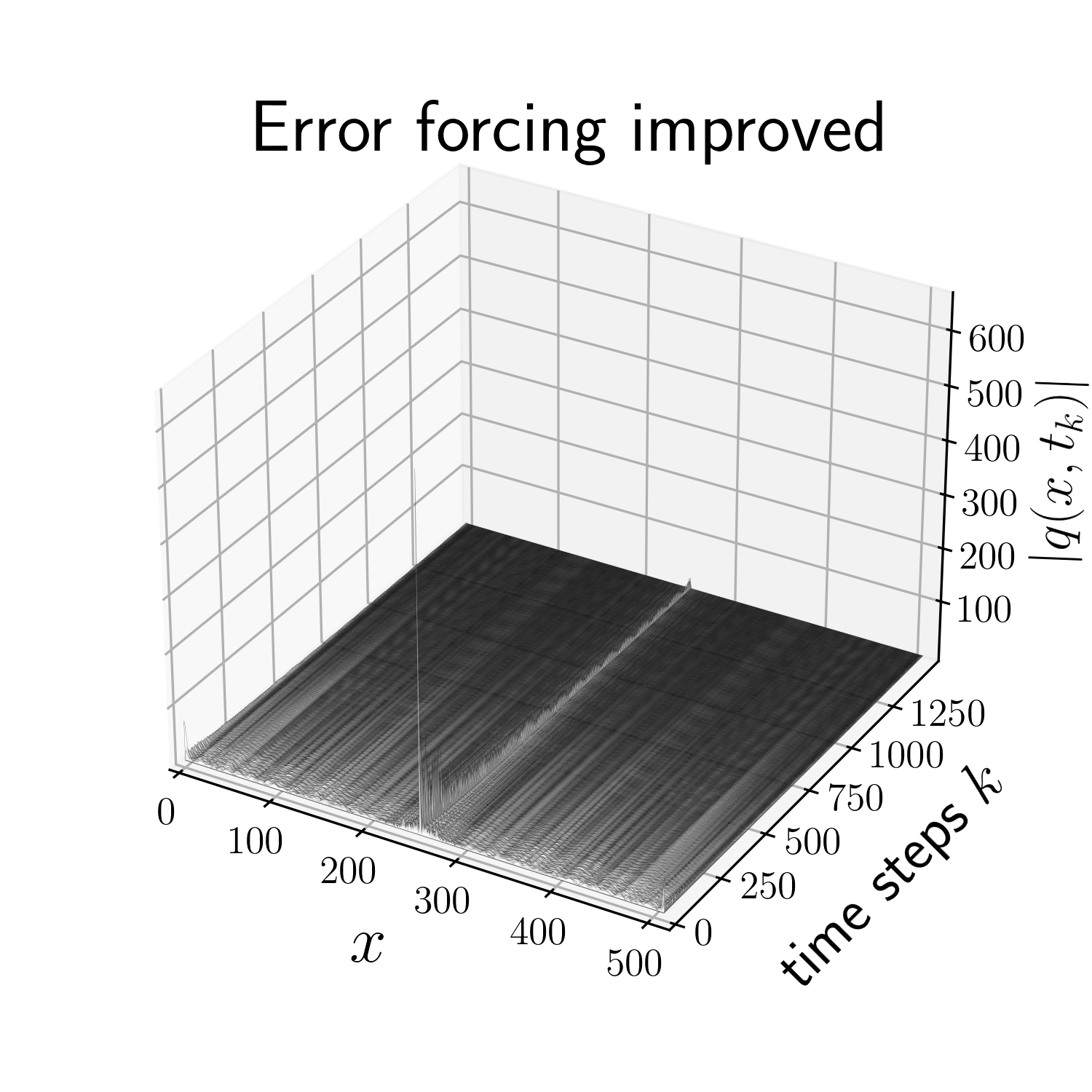}
    \includegraphics[width =0.32\textwidth]{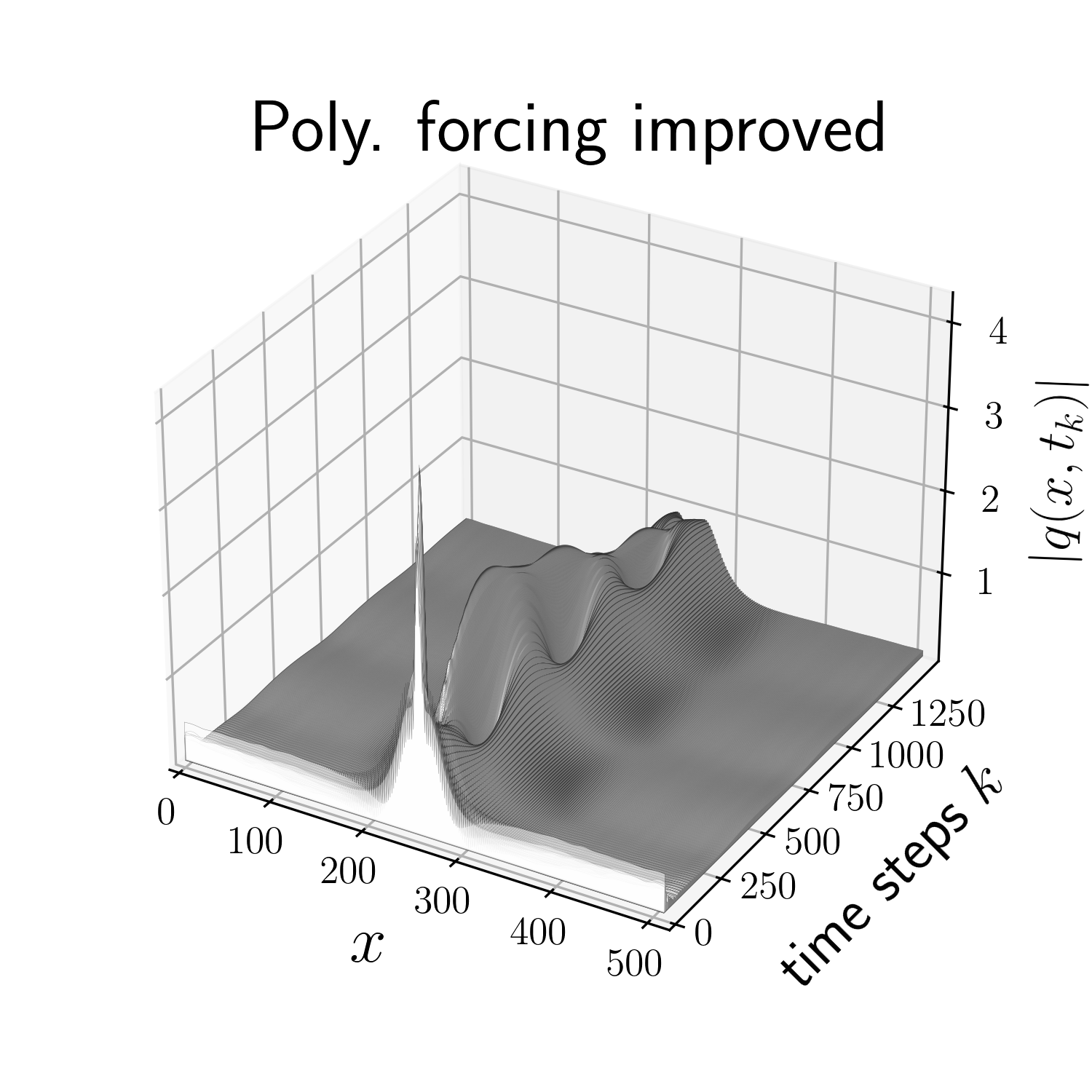}
\caption{Reconstructions and forward predictions of the NLS dynamics obtained with CB-DMD, the baseline methods, and their improved variants. The operators are learned from $1000$ snapshots, whereas the following $500$ time steps correspond to prediction. Only non-diverging results are shown. CB-DMD gives the most accurate reconstruction and the most reliable prediction of the periodic wave pattern.}
  \label{fig:NLS_recon}
  \end{figure}

This spectral behavior is reflected in the reconstruction and prediction errors shown in Fig.~\ref{fig:NLS_error}. The operators from standard DMD, its improved variant, polynomial forcing, and error forcing generally lead to diverging behavior. Although the improved polynomial-forcing model does not diverge, it is still less accurate than CB-DMD. The improved error-forcing model first reaches an error of about $600$ and then decreases, but this is caused by the reconstruction decaying toward zero rather than by correctly reproducing the solution. In contrast, CB-DMD gives the most accurate non-diverging reconstruction and the most reliable forward prediction. This is also clear in Fig.~\ref{fig:NLS_recon}, showing only the non-diverging results. Here, CB-DMD most closely matches the true periodic wave pattern. Overall, for the NLS example, CB-DMD gives the most stable and accurate linear representation among the tested methods.

\subsection{Cylinder wake flow}
The flow behind a circular cylinder is a canonical problem in engineering. It's properties are of interest since self-sustained oscillations can cause vibrations, noise and resonance in practice \cite{Bagheri2013}.

\begin{figure}[t]
    \centering
      \includegraphics[width = \textwidth]{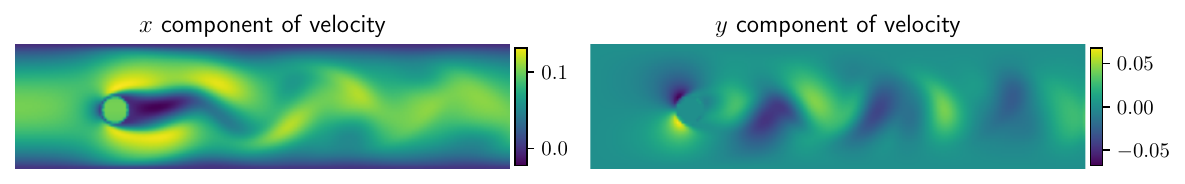}
    \caption{Example snapshot of the cylinder flow data. }
    \label{fig:Cylinder_example}
  \end{figure}

For our experiments, we use the simulations from \cite{resDMD}. The Reynolds number is set to $Re = 100$ and originally the domain is discretized using a mesh of size 400x100.
Due to computational constraints we use a sub-sampled version of their data, reducing the mesh to 200x50. The velocity is given in 2D, yielding 20K measurements per snapshot. An exemplary snapshot is shown in Fig.~\ref{fig:Cylinder_example}.
   \begin{figure}
    \centering
      \includegraphics[width=0.3\textwidth]{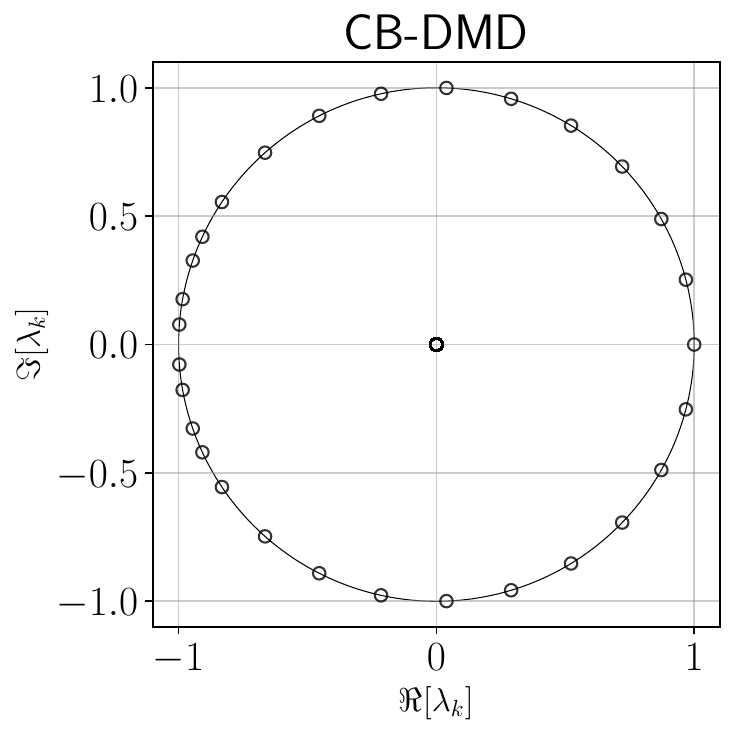}
      \includegraphics[width=0.3\textwidth]{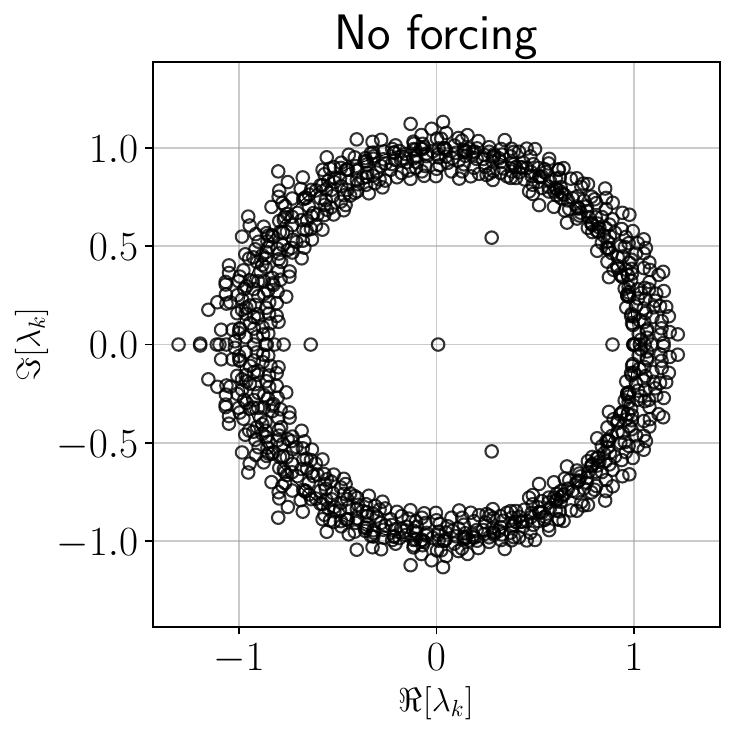}\\
     \includegraphics[width=0.3\textwidth]{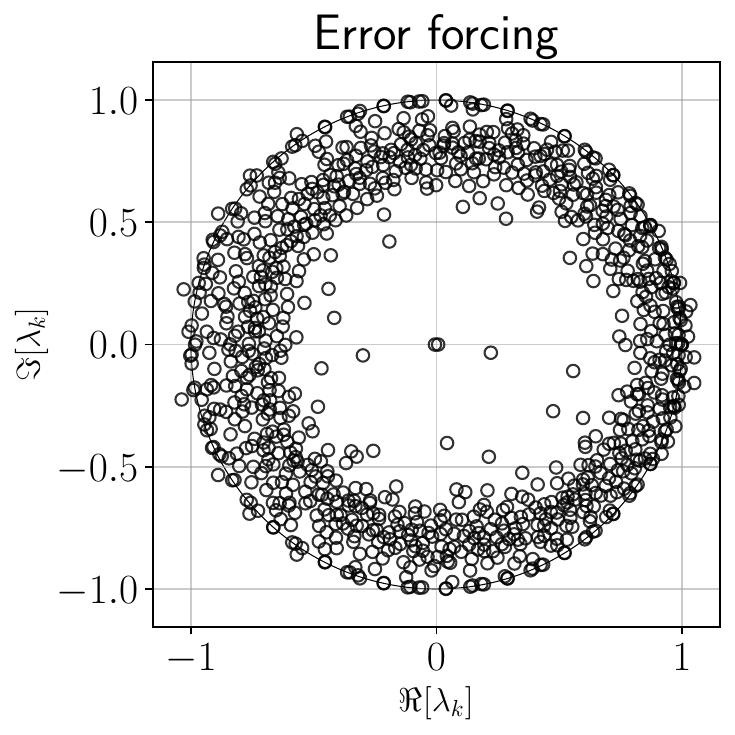}
      \includegraphics[width=0.3\textwidth]{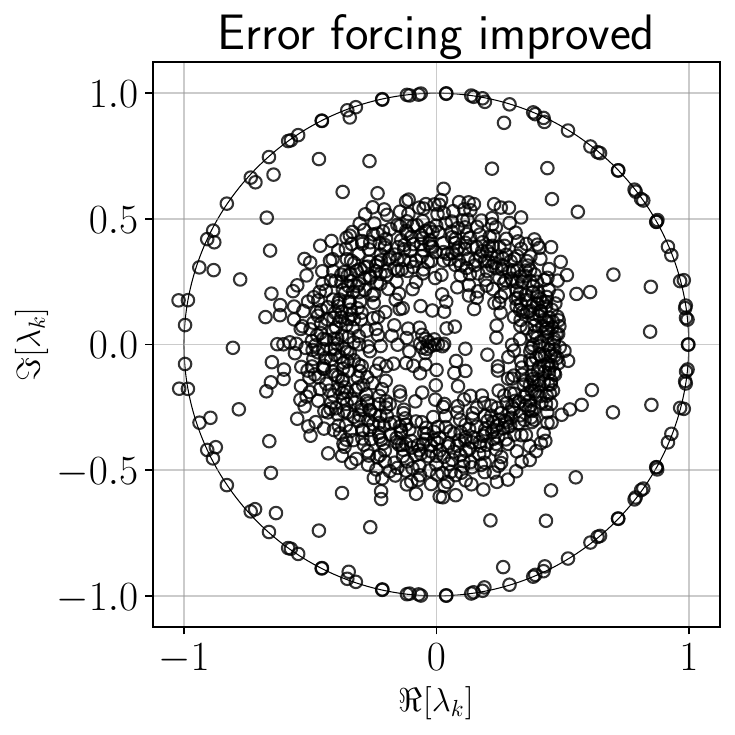}
      \includegraphics[width=0.3\textwidth]{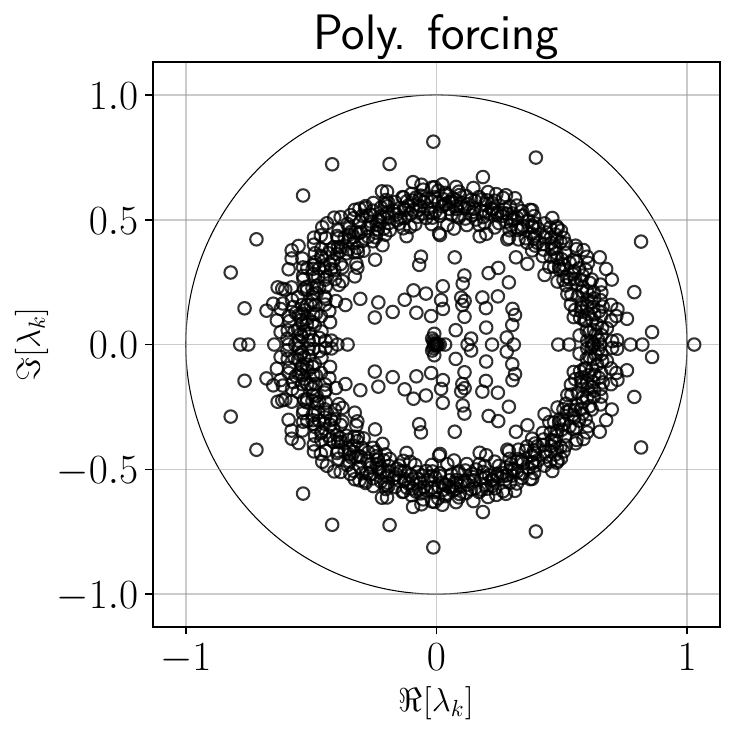}
    \caption{Spectra of the operators learned for the cylinder wake dynamics using standard DMD, polynomial forcing, error forcing, and CB-DMD. Improved variants are obtained by removing DMD modes associated with eigenvalues of modulus greater than one; for error forcing, this is done both for the initial standard-DMD operator and for the final DMDc operator. Relevant eigenvalues are expected near the unit circle, while eigenvalues away from it are interpreted as spurious. The circular line indicates the unit circle. CB-DMD again yields the most concentrated spectrum, with $29$ numerically relevant eigenvalues ($|\lambda|>10^{-3}$), all satisfying $0.99\leq|\lambda|\leq 1.01$.}
    \label{fig:Cylinder_spectrum}
  \end{figure}
Like the authors of \cite{resDMD} we consider the post-transient regime, starting 3000 time steps from the beginning of the simulations. To compute the linear operators we use the subsequent $M = 1000$ time steps and the prediction error is computed over the following 500 time steps. 

In the setting described in  \cite{Bagheri2013}, the spectrum of the Koopman operator of the cylinder wake flow forms a lattice structure on the unit circle in the complex plane.

Even with our initial down-sampling it is not feasible to compute for polynomial forcing the PCCs of all combinations of measurements. Hence, we compute the forcing using only every tenth original grid point. The PCC threshold is set to 0.95.

Note that although the cylinder snapshots belong to a state space of dimension $20$K, the DMD spectra contain only $999$ eigenvalues because the training set consists of $1000$ snapshots, hence $999$ snapshot pairs. DMD constructs $\mathbf{\hat{g}(X)}$ and $\mathbf{\hat{g}(\overline{X})}$ as in (\ref{observed_snapshots}), with $\mathbf{\hat{g}(X)}\in\mathbb{R}^{20000\times 999}$ and $\operatorname{rank}(\mathbf{\hat{g}(X)})\leq 999$. Therefore, the DMD based operators are identified only on the snapshot subspace spanned by the columns of $\mathbf{\hat{g}(X)}$, and the reduced operator has dimension at most $999$, independently of the ambient state dimension.

Fig.~\ref{fig:Cylinder_spectrum} shows the spectra of the operators obtained with CB-DMD and the baseline methods for the cylinder wake flow dynamics. According to the theoretical analysis  proposed in \cite{Bagheri2013}, the physically meaningful eigenvalues are expected to be concentrated near the unit circle, while eigenvalues far from it are interpreted as spurious. 
The spectral results show that CB-DMD produces the most concentrated spectrum around the unit-circle structure. 


\begin{table}
\centering
\caption{Comparison of eigenvalue properties across DMD methods for the cylinder wake system. We consider an eigenvalue relevant if $|\lambda|>10^{-3}$. We also report eigenvalues close to the unit circle, defined by $0.99<|\lambda|\leq 1.01$, and eigenvalues outside the unit circle, defined by $|\lambda|>1$.}
\label{tab:eigenvalue_comparison_cylinder}
\begin{tabular}{lcccc}
\hline
\textbf{Method} 
& \makecell{\textbf{Relevant}\\\textbf{eigenvalues}} 
& \makecell{\textbf{Near}\\\textbf{unit circle}} 
& \makecell{\textbf{Eigenvalues with}\\\boldmath{$|\lambda|>1$}} 
& \makecell{\textbf{Largest}\\\boldmath{$|\lambda|$}} \\
\hline
CB-DMD                        & 29  & 29  & 4   & 1.0000 \\
Standard DMD                  & 999 & 196 & 463 & 1.3054             \\
Error Forcing                 & 999 & 145 & 72  & 1.0619 \\
Error Forcing (Improved)      & 999 & 100 & 24  & 1.0364             \\
Polynomial Forcing            & 993 & 0   & 1   & 1.0287             \\
\hline
\end{tabular}
\end{table}

Table~\ref{tab:eigenvalue_comparison_cylinder} summarizes the eigenvalue properties across the different methods. As before, we call an eigenvalue relevant if $|\lambda|>10^{-3}$, which provides a qualitative cutoff for identifying spectral components with a non-negligible contribution to the dynamics. CB-DMD retains only $29$ relevant eigenvalues, and all of them lie close to the unit circle, satisfying $0.99<|\lambda|\leq 1.01$. In contrast, the baseline methods retain substantially larger relevant spectra. Standard DMD, error forcing, and its improved version all retain $999$ relevant eigenvalues, but only a fraction of these lie close to the unit circle. Polynomial forcing retains $993$ relevant eigenvalues, none of which satisfy the near-unit-circle criterion, indicating that its relevant spectrum is almost entirely separated from the expected unit-circle structure. The comparison of eigenvalues outside the unit circle leads to the same conclusion. CB-DMD has only $4$ eigenvalues with $|\lambda|>1$, and its largest eigenvalue modulus is $1.0000000000089522$, indicating only a negligible numerical deviation from the unit circle. By contrast, the baseline methods exhibit larger deviations from the unit circle, either through many eigenvalues with $|\lambda|>1$ or through noticeably larger maximum eigenvalue moduli.
\begin{figure}
  \centering
    \includegraphics[width=0.7\textwidth]{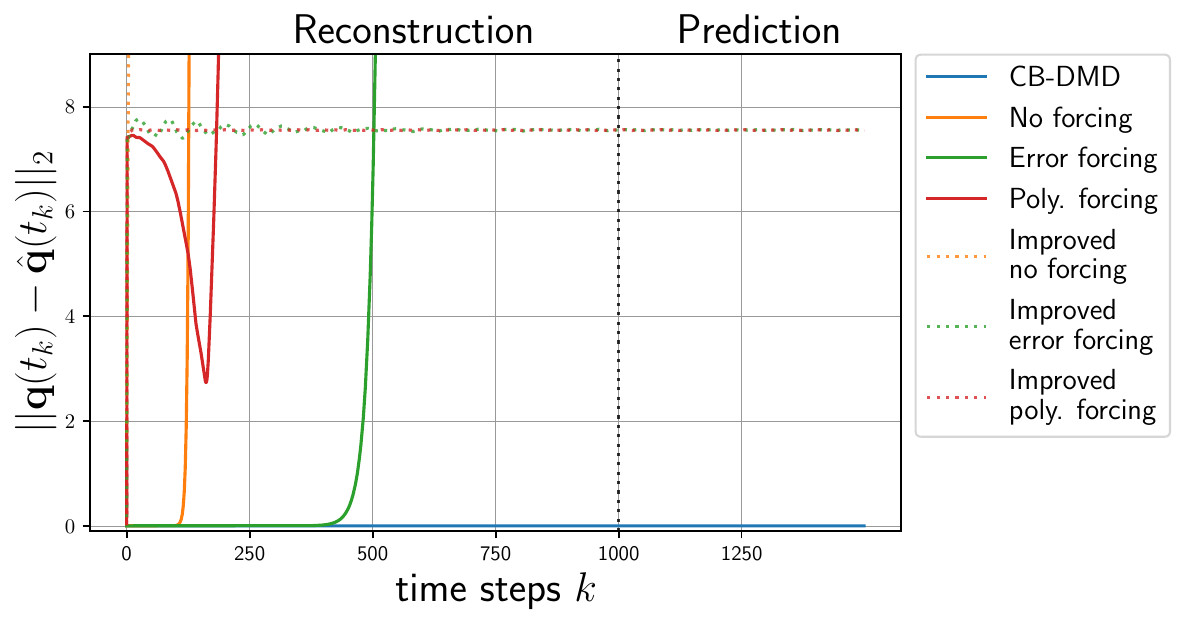}
  \caption{Absolute error of the reconstructions and forward predictions for the cylinder wake dynamics obtained with the linear operators learned by CB-DMD, the baseline methods, and their improved variants. The first $1000$ time steps correspond to reconstruction, and the following $500$ time steps correspond to prediction. The vertical dashed black line marks the start of the prediction period.}
  \label{fig:Cylinder_error}
  \end{figure}
  \begin{figure}
  \centering
    \includegraphics[width=0.32\textwidth]{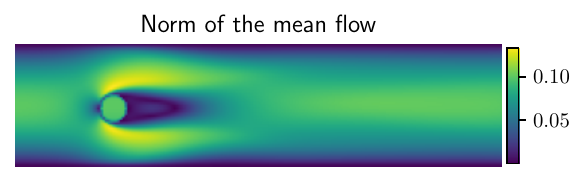}
    \includegraphics[width=0.32\textwidth]{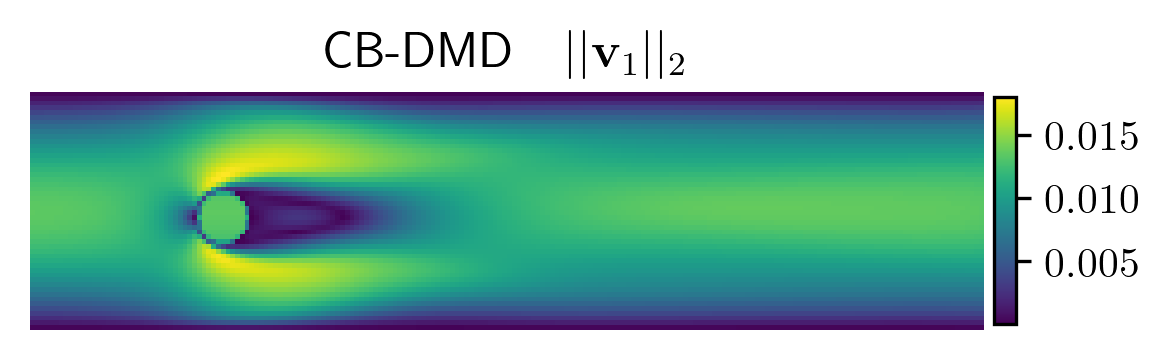}
    \includegraphics[width=0.32\textwidth]{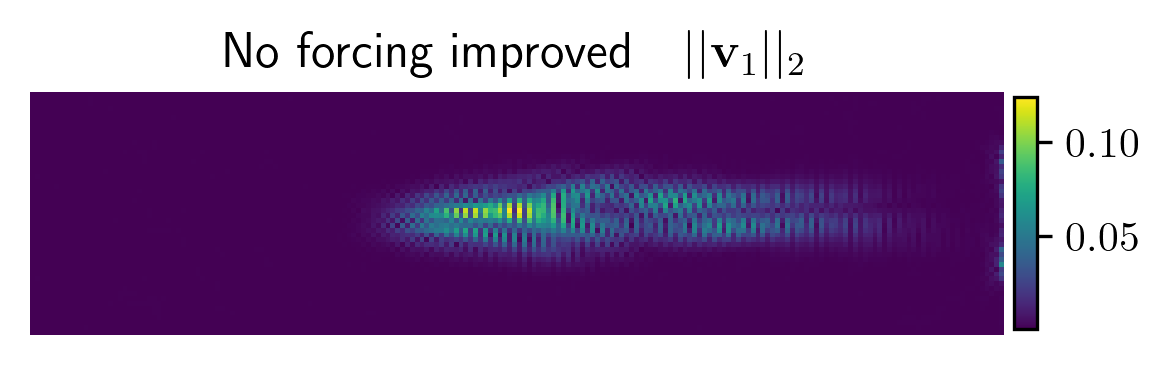}\\
    \includegraphics[width=0.32\textwidth]{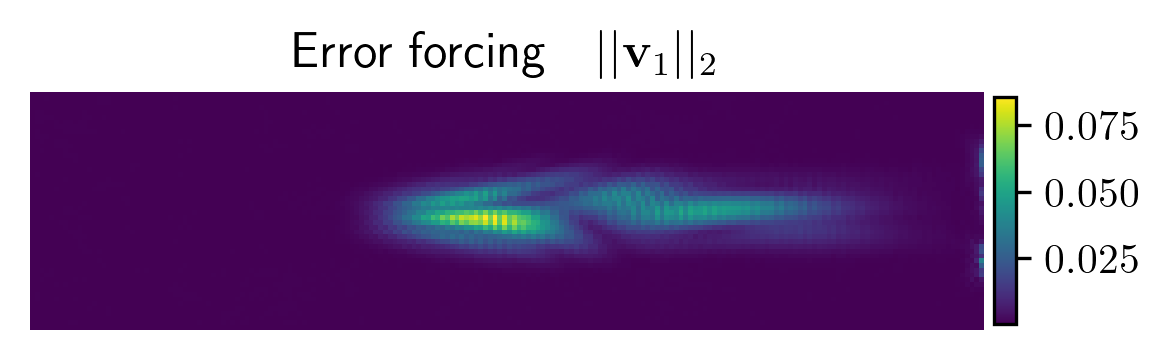}
    \includegraphics[width=0.32\textwidth]{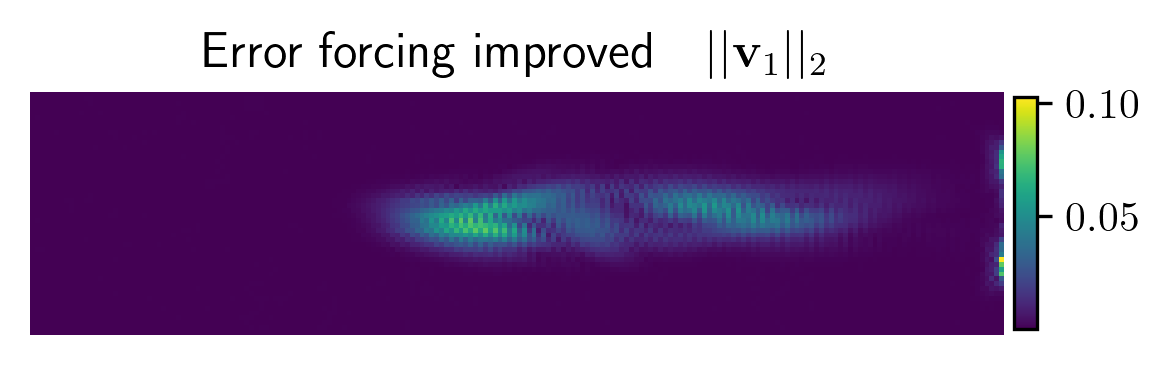}
    \includegraphics[width=0.32\textwidth]{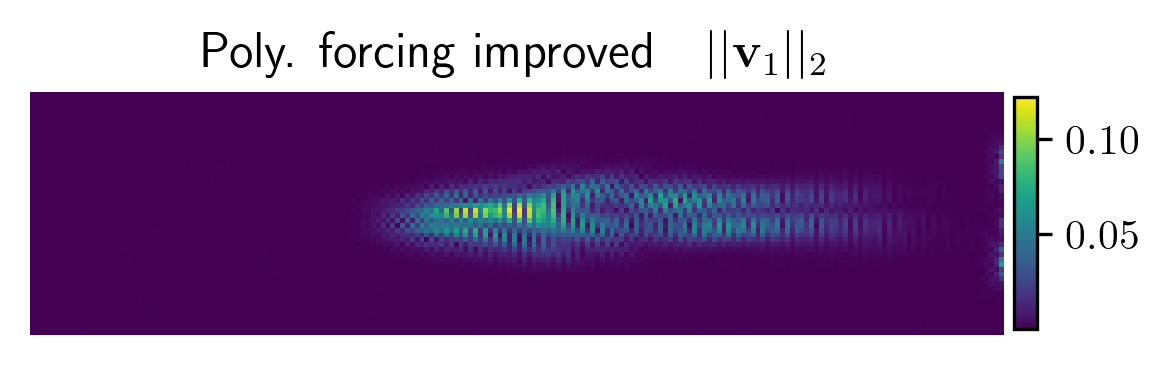}
  \caption{Euclidean norm of the mean flow and dominant modes for the cylinder wake dynamics. The first panel shows the Euclidean norm of the temporal mean of the training snapshots. The remaining panels show the pointwise Euclidean norm of the mode $\mathbf{v}_1$ associated with the eigenvalue of largest modulus for each linear operator whose reconstruction remains bounded. The CB-DMD mode norm resembles the mean-flow norm, indicating that the dominant CB-DMD mode captures a physically meaningful flow structure.}
  \label{fig:Cylinder_modes}
  \end{figure}
  
Overall, CB-DMD provides the most physically consistent spectral representation of the cylinder wake dynamics. Its relevant eigenvalues are few in number, all concentrated near the unit circle, and its largest deviation from the unit circle is several orders of magnitude smaller than those observed for the baseline methods.



Fig.~\ref{fig:Cylinder_error} reports the linear reconstruction and prediction errors for the cylinder wake dynamics obtained with CB-DMD, the baseline methods, and their improved variants. CB-DMD achieves the smallest errors in both reconstruction and prediction. In particular, its error remains close to zero throughout the time interval, with a maximum value of $1.02\times 10^{-5}$.

The error behavior of the baseline methods reflects either unstable growth or bounded but inaccurate dynamics. Reconstructing the dynamics not considering components associated with eigenvalues outside the unit circle prevents divergence for the improved baselines, but does not yield an accurate representation of the cylinder wake dynamics. In particular, improved DMD has an average error of approximately $33$ over the considered time interval, due to the collapse of the reconstructed dynamics toward zero. Similarly, improved error forcing and improved poly forcing do not recover the correct flow evolution. Instead, they stabilize on periodic motions with magnitudes that differ substantially from those of the true solution. We remark that the improved baselines reconstruct and predict the dynamics after removing modes associated with eigenvalues satisfying $|\lambda|>1$, since these modes are treated as spurious and may induce unstable growth. No method-specific eigenvalue-threshold tuning is used in the main comparison. As shown in Appendix~D, choosing the best threshold a posteriori can make standard DMD and error forcing achieve reconstruction and prediction errors comparable to those of CB-DMD; however, this requires additional fine-tuning and is therefore excluded from the default setting.

Fig.~\ref{fig:Cylinder_modes} shows the modes associated with the eigenvalue with the largest moduls for CB-DMD and for the non-diverging alternatives. The CB-DMD mode captures the mean flow structure of the dynamics, whereas the modes produced by the other non-diverging operators are nearly vanishing. The modes associated with the largest-modulus eigenvalues of the diverging operators are not reported, since they are almost constant and also close to zero.

\subsection{Cubic-quintic Ginzburg-Landau equation}
The complex cubic–quintic Ginzburg--Landau (CQGL) equation originates in the theory of superconductivity \cite{Ginzburg2009} and has since found applications in areas such as biological pattern formation \cite{Morales2015}.
The evolution equation is given by
\begin{align} iq_t + (\frac{1}{2} - i\tau)q_{xx} + i\kappa q_{xxxx} +(1-i\beta)|q|^2q + (\nu - i \sigma)|q|^4 q -i\gamma q = 0,
\end{align} 
where we set parameters to $$\tau = 0.08, \kappa = 0, \beta = 0.66, \nu = -0.1, \sigma = -0.1, \gamma = -0.1,$$ and use the initial condition 
$$q(x, 0) = 2 \sech(x),$$
yielding a periodic solution. 
The space domain $[-10, 10]$ is discretized with $n=512$ points. 
We obtain the solution to the CGQL equation by applying the explicit Runge-Kutta (4,5)-formula, after Fourier transforming the space variable to integrate the solution. 
To compute the linear operators we use the first $M = 1000$ snapshots comprising the time interval $t \in [0, 40]$. The prediction interval is given by the subsequent $500$ time steps up to time $60$. 
 \begin{figure}[t]
    \centering
      \includegraphics[width=0.30\textwidth]{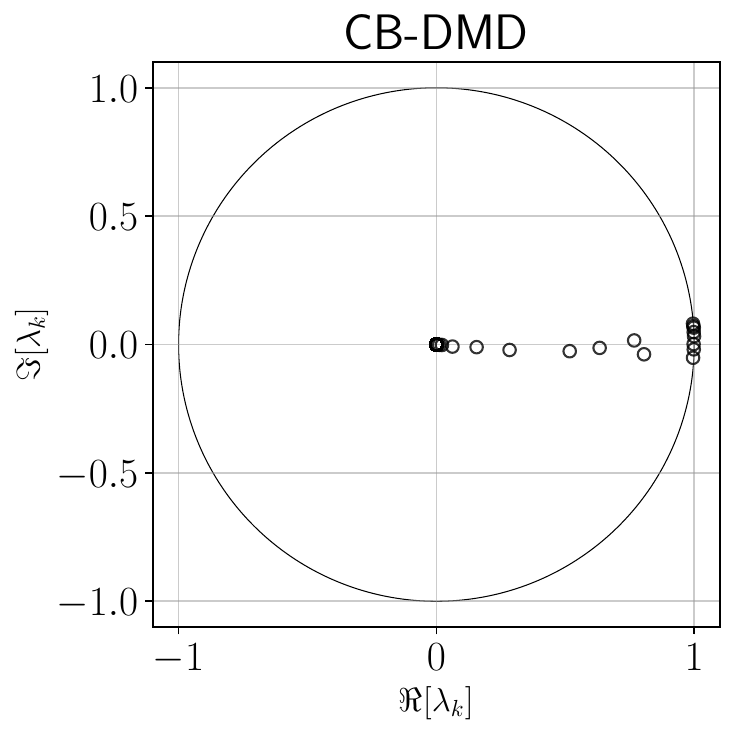}
      \includegraphics[width=0.3\textwidth]{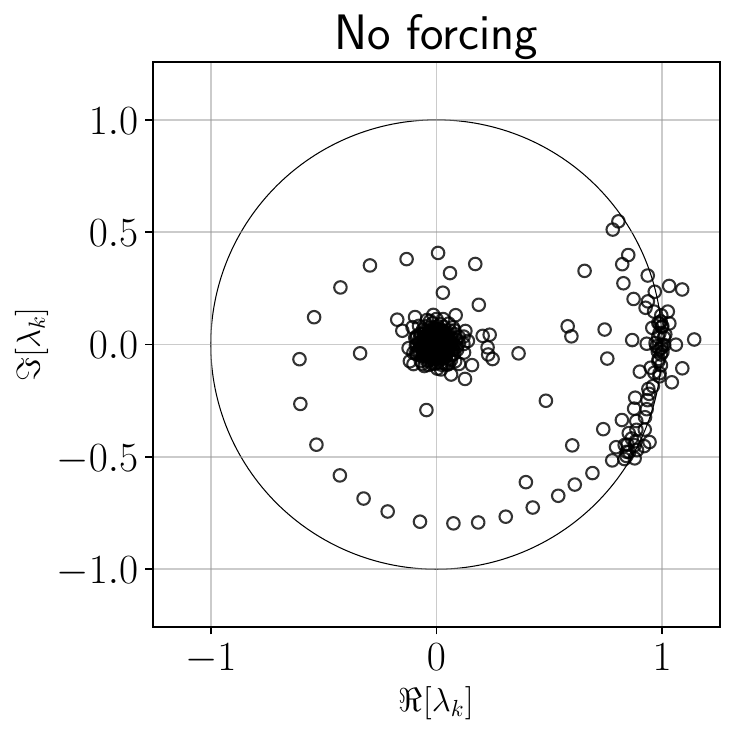}\\
     \includegraphics[width=0.3\textwidth]{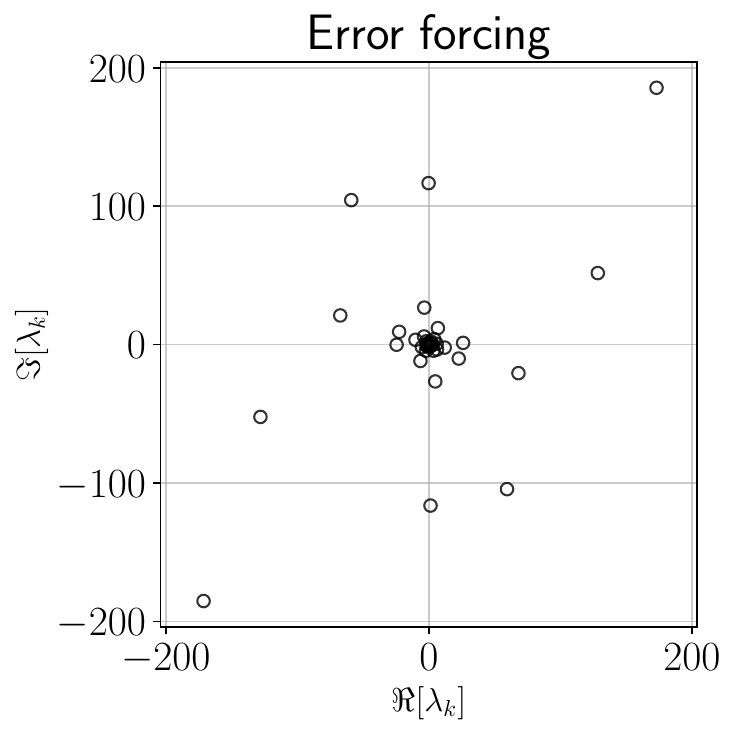}
     \includegraphics[width=0.3\textwidth]{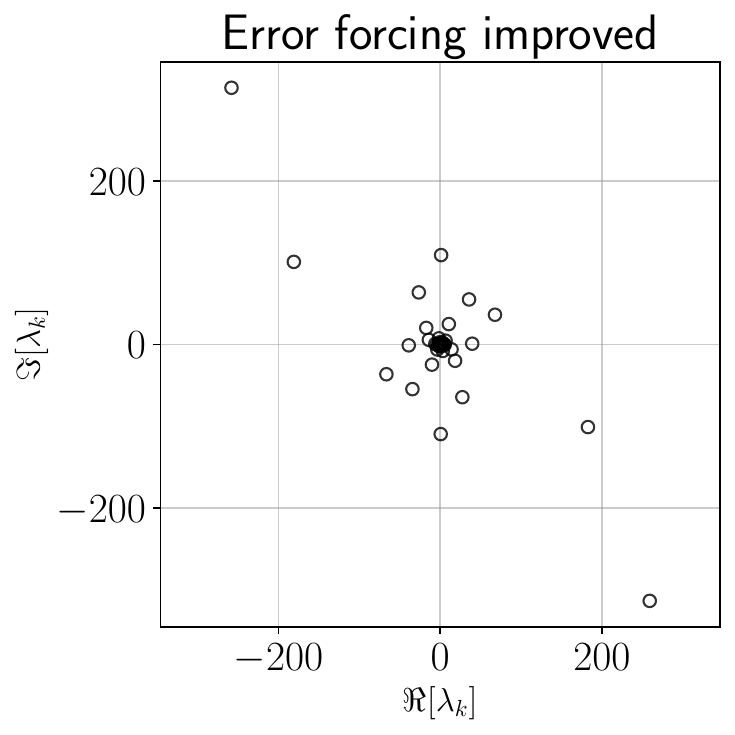}
      \includegraphics[width=0.3\textwidth]{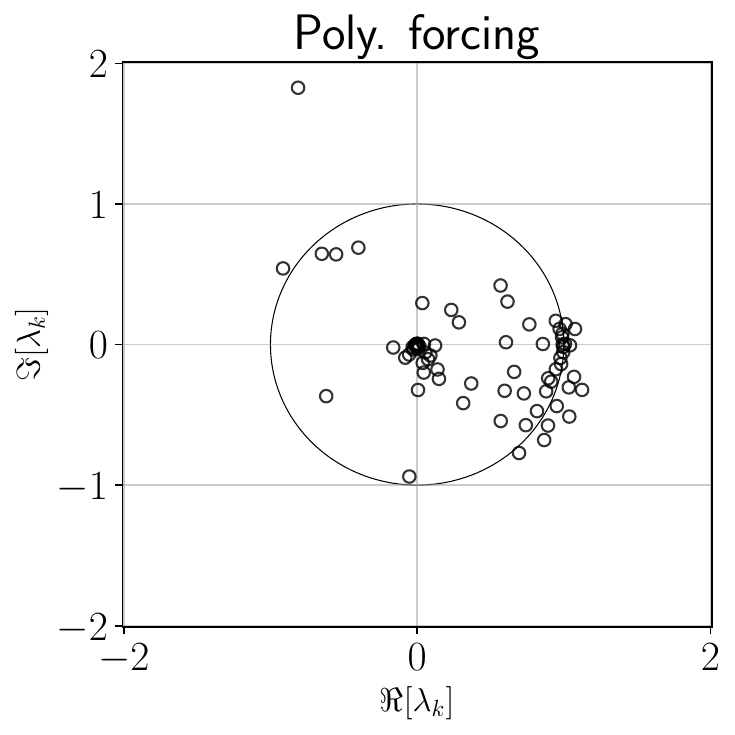}
  \caption{Spectra of the operators learned for the CQGL dynamics using standard DMD, polynomial forcing, error forcing, and CB-DMD. Improved variants are obtained by removing DMD modes associated with eigenvalues of modulus greater than one; for error forcing, this is done both for the initial standard-DMD operator and for the final DMDc operator. Relevant eigenvalues are expected near the unit circle, while eigenvalues away from it are interpreted as spurious. The circular line indicates the unit circle. CB-DMD again yields the most concentrated spectrum, with $19$ numerically relevant eigenvalues ($|\lambda|>10^{-3}$), of which $8$ satisfy $0.99\leq|\lambda|\leq 1.01$.}
    \label{fig:CQGL_spectrum}
  \end{figure}
Similarly to the NLS equation, for poly. forcing we can simplify the computation of the polynomial basis functions by excluding the outer 106 grid points on both ends of the one dimensional spatial grid, since they are close to zero during the entire time period of interest. Again, a reasonable PCC threshold is determined to be 0.2.
\begin{table}
\centering
\caption{Comparison of eigenvalue properties across DMD methods for the CQGL system. We consider an eigenvalue relevant if $|\lambda|>10^{-3}$. We also report eigenvalues close to the unit circle, defined by $0.99\leq|\lambda|\leq 1.01$, and eigenvalues outside the unit circle, defined by $|\lambda|>1$.}
\label{tab:eigenvalue_comparison_CQGL}
\begin{tabular}{lcccc}
\hline
\textbf{Method} 
& \makecell{\textbf{Relevant}\\\textbf{eigenvalues}} 
& \makecell{\textbf{Near}\\\textbf{Unit circle}} 
& \makecell{\textbf{Eigenvalues with}\\\boldmath{$|\lambda|>1$}} 
& \makecell{\textbf{Largest}\\\boldmath{$|\lambda|$}} \\
\hline
CB-DMD                        & 19  & 8   & 1  & 1.0010   \\
Standard DMD                  & 512 & 18   & 20 & 1.1428   \\
Error Forcing                 & 130 & 0   & 52 & 253.7228 \\
Error Forcing (Improved)      & 131 & 3   & 51 & 407.0758 \\
Polynomial Forcing            & 107 &  6   & 14 & 1.9986   \\
\hline
\end{tabular}
\end{table}

The spectral results shown in Fig.~\ref{fig:CQGL_spectrum} again favor CB-DMD. Table~\ref{tab:eigenvalue_comparison_CQGL} summarizes the eigenvalue properties across the different methods. As before, we call an eigenvalue relevant if $|\lambda|>10^{-3}$.
CB-DMD retains only $19$ relevant eigenvalues, and $8$ of them lie close to the unit circle, satisfying $0.99\leq|\lambda|\leq 1.01$. 

In contrast, the baseline methods retain substantially larger relevant spectra. Standard DMD retains all $512$ eigenvalues as relevant, but only a small fraction of them lie close to the unit circle. Error forcing and its improved version retain fewer relevant eigenvalues than standard DMD, but their spectra are still much less concentrated near the unit circle. In particular, the linear error forcing operator has no relevant eigenvalues satisfying the near-unit-circle criterion, while its improved version and polynomial forcing contain only a few such eigenvalues.

The comparison of eigenvalues outside the unit circle leads to the same conclusion. CB-DMD has only one eigenvalue with $|\lambda|>1$, and its largest eigenvalue modulus is $1.0010$, indicating only a negligible deviation from the unit circle. By contrast, the baseline methods exhibit either more eigenvalues outside the unit circle or much larger maximum eigenvalue moduli. This is especially pronounced for the forcing-based approaches, whose largest eigenvalue moduli are several orders of magnitude larger than that of CB-DMD.

Overall, CB-DMD provides the cleanest spectral representation for the CQGL system. It produces a compact set of relevant modes, places the dominant spectral content closest to the expected unit-circle structure, and avoids the large artificial growth effects observed in the baseline methods.

\begin{figure}[t]
  \centering
    \includegraphics[width=0.7\textwidth]{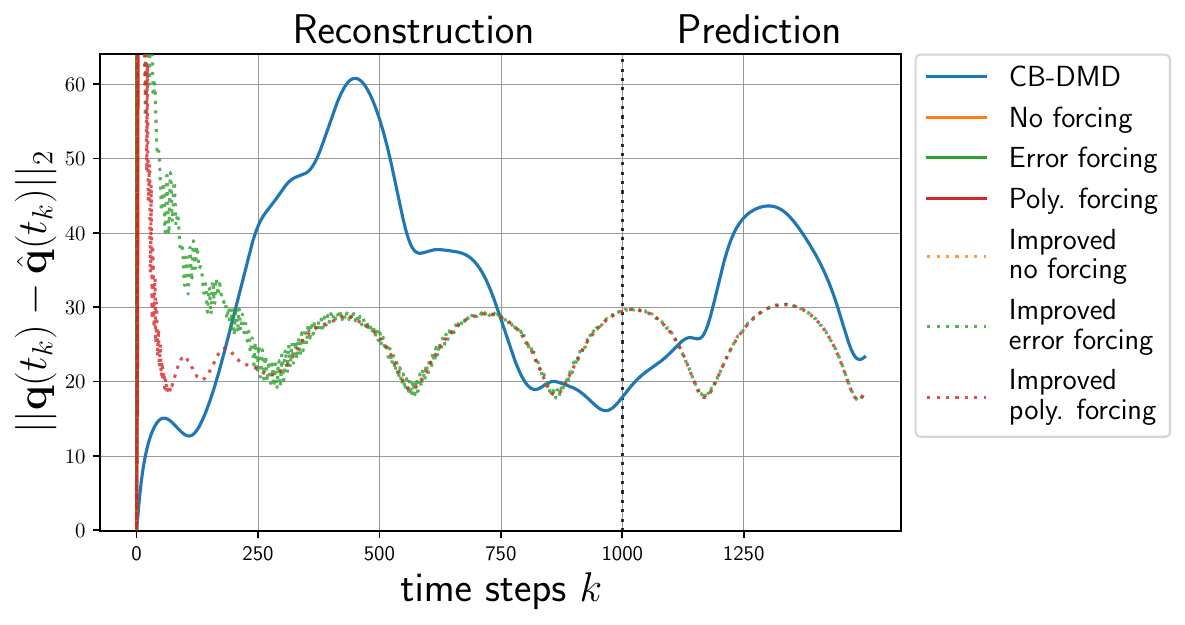}
\caption{Absolute error of the reconstructions and forward predictions for the CQGL dynamics obtained with the linear operators learned by CB-DMD, the baseline methods, and their improved variants. The first $1000$ time steps correspond to reconstruction, and the following $500$ time steps correspond to prediction. The vertical dashed black line marks the start of the prediction period.}
  \label{fig:CQGL_error}
  \end{figure}

\begin{figure}[t]
\centering
    \includegraphics[width =0.32\textwidth]{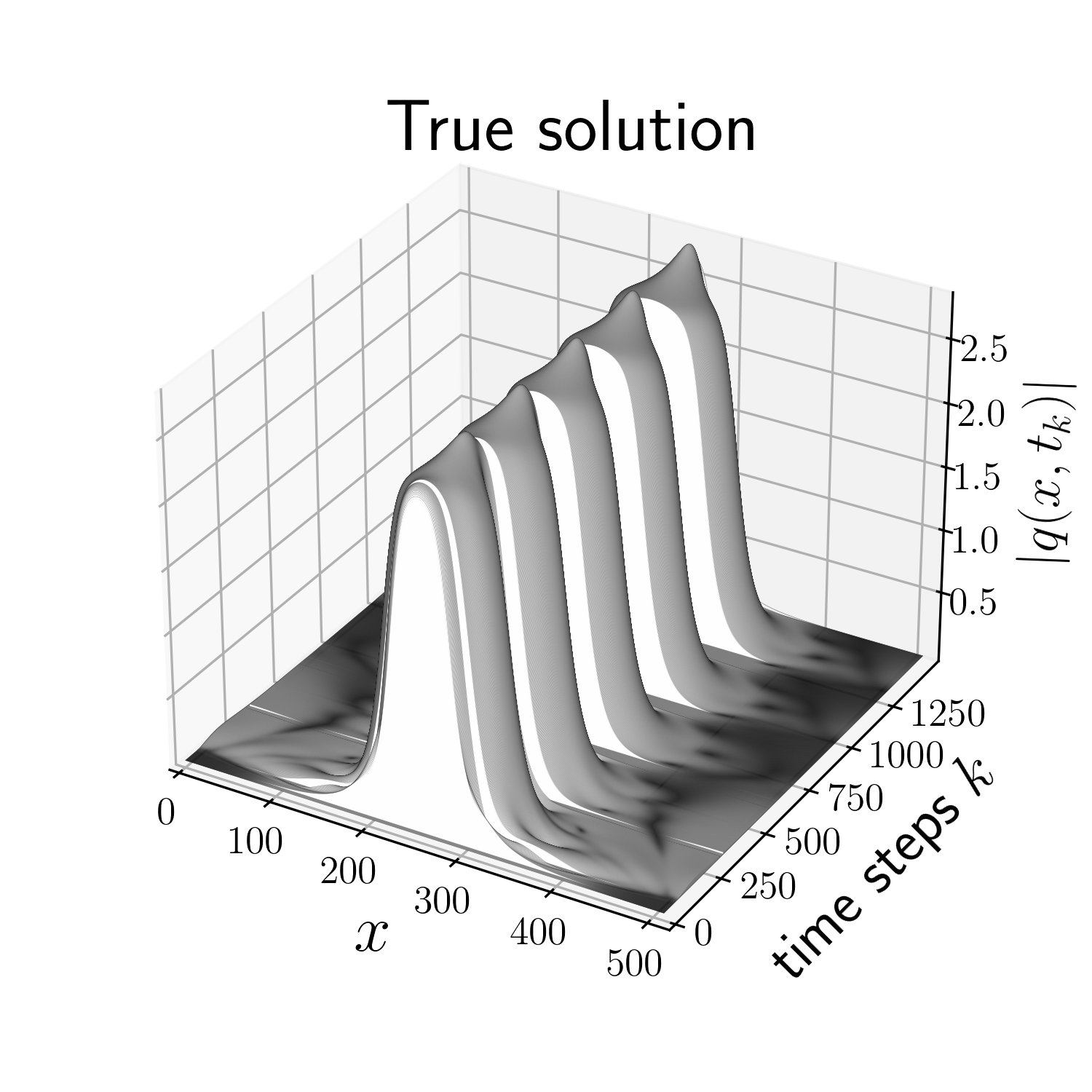}
    \includegraphics[width =0.32\textwidth]{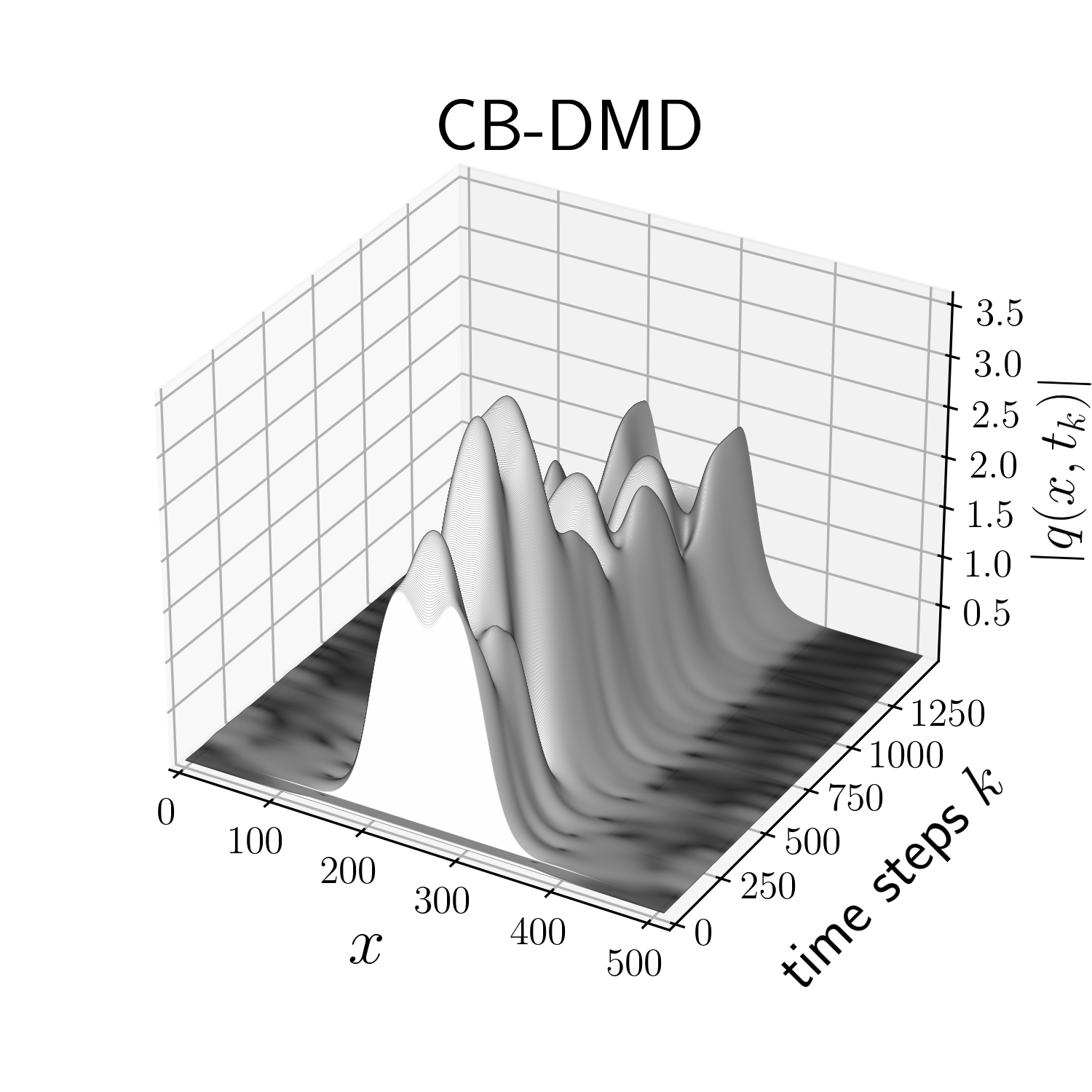}
    \includegraphics[width =0.32\textwidth]{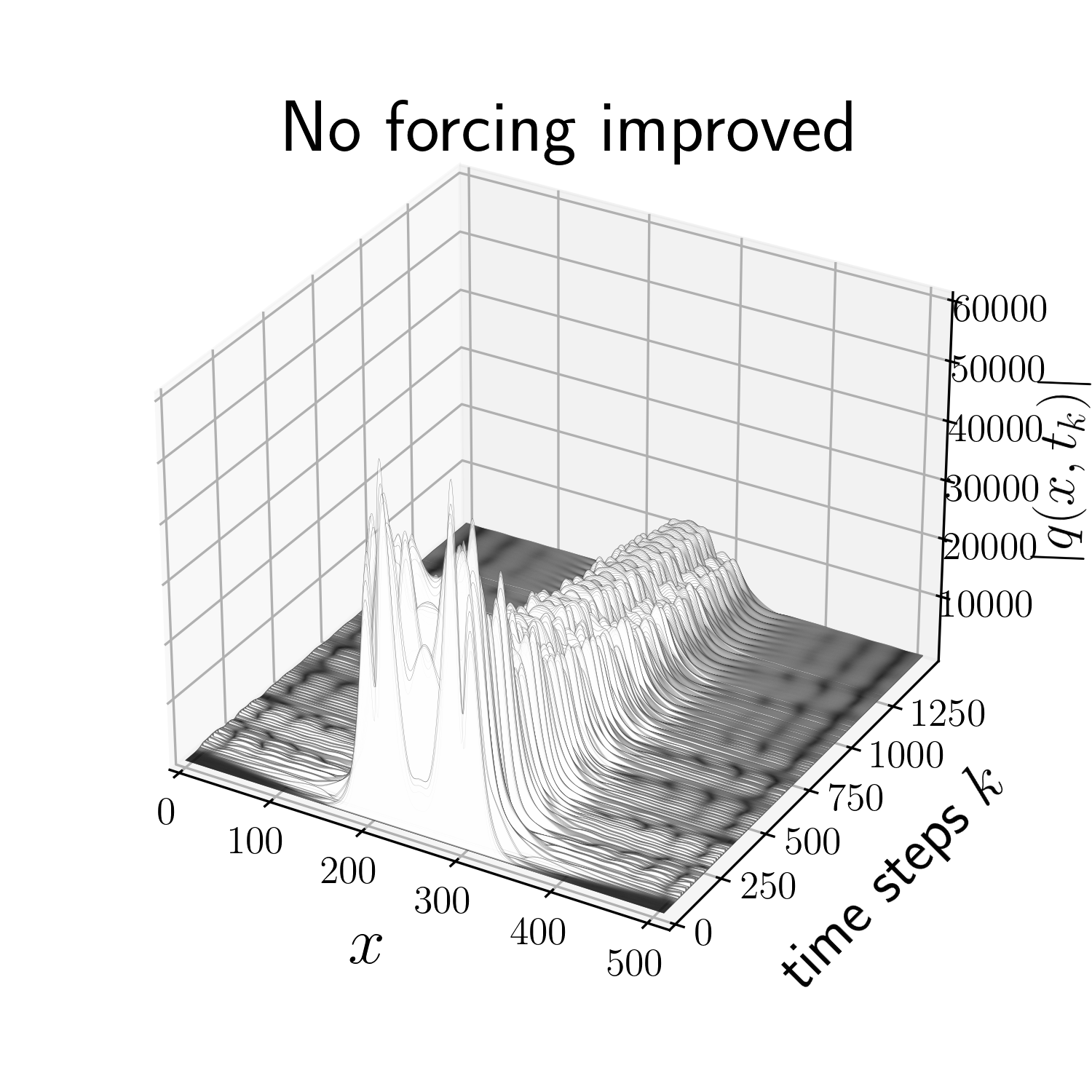}\\
    \includegraphics[height=0.32\textwidth]{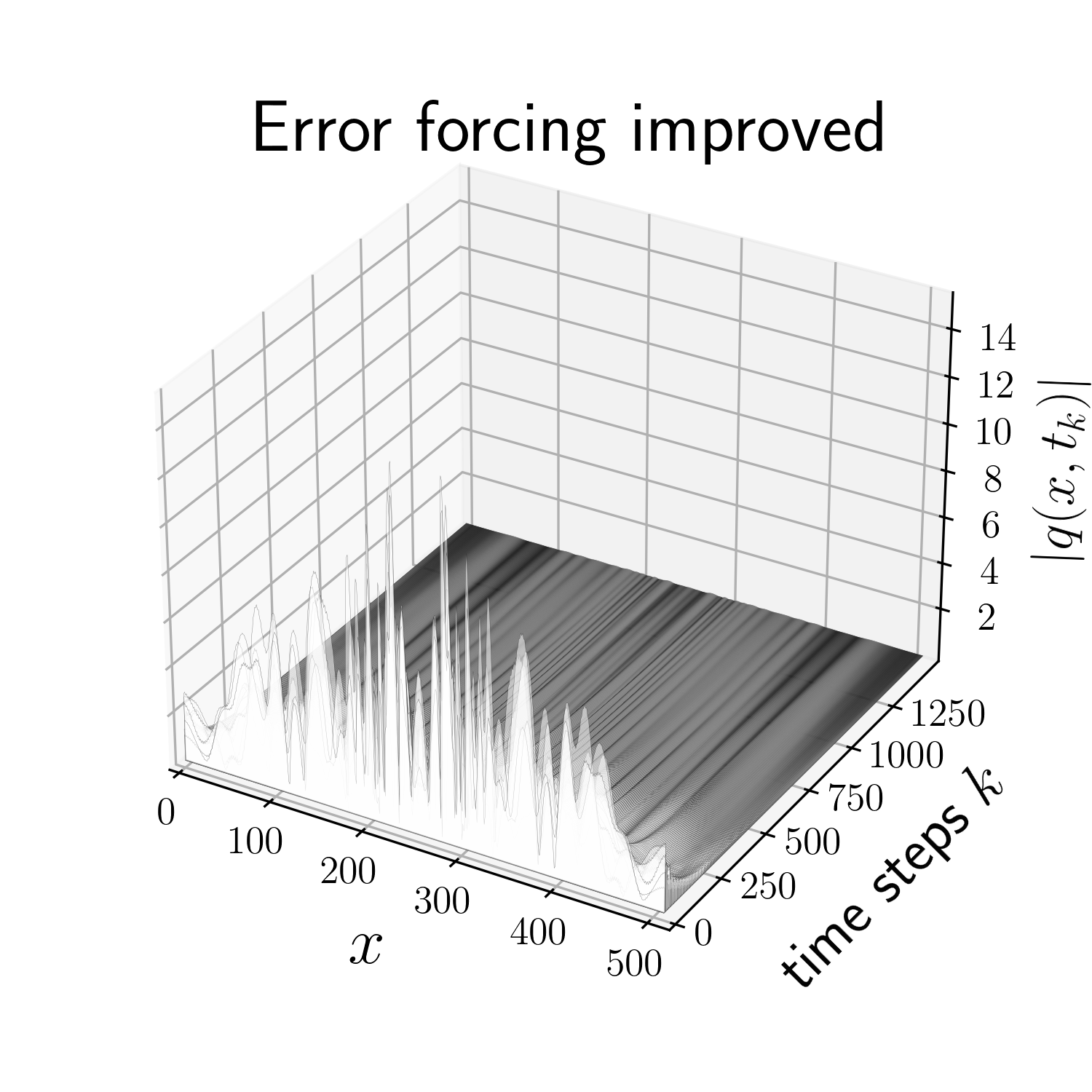}
    \includegraphics[height=0.32\textwidth]{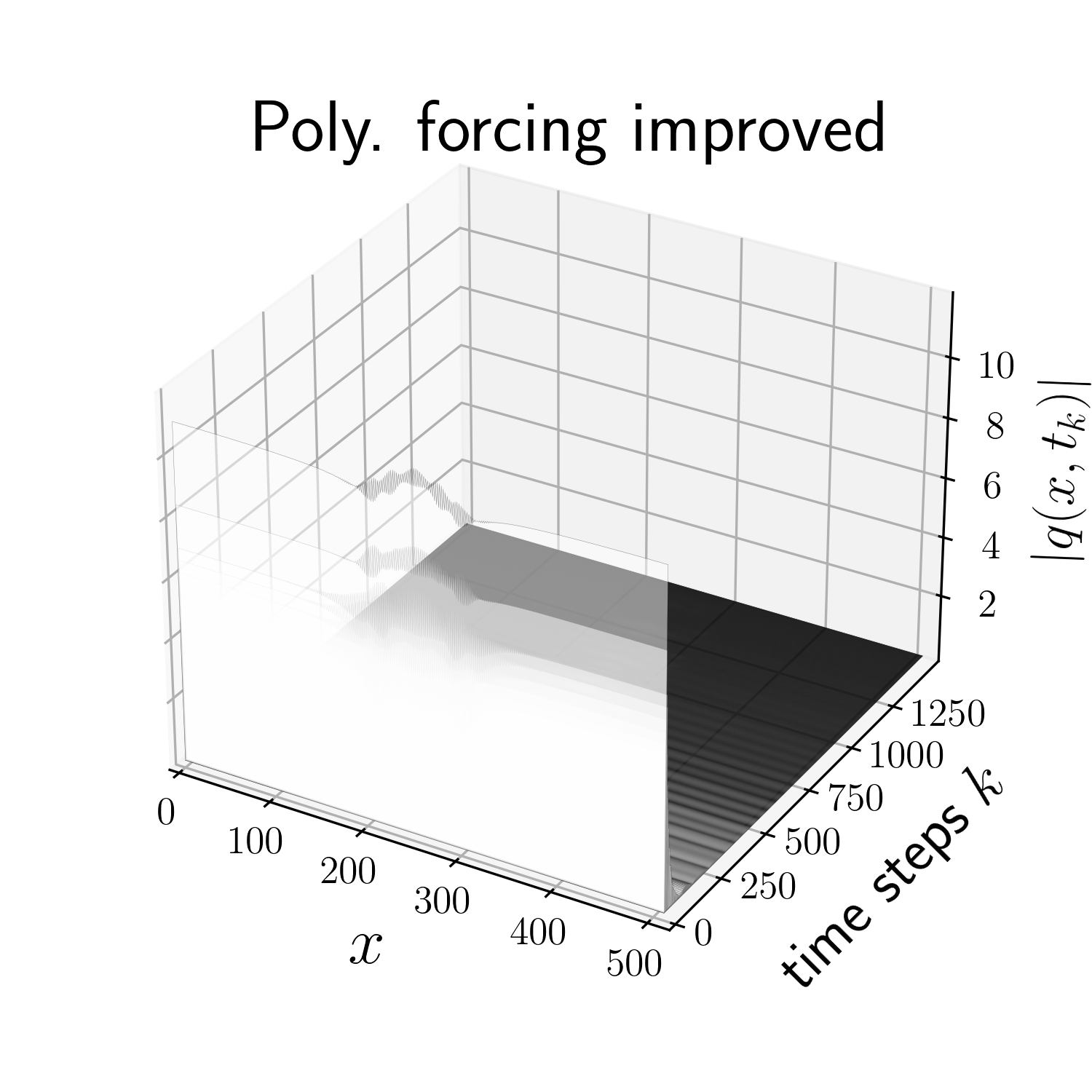}
\caption{Reconstructions and forward predictions of the CQGL dynamics obtained with CB-DMD, the baseline methods, and their improved variants. The operators are learned from $1000$ snapshots, whereas the following $500$ time steps correspond to prediction. Only non-diverging results are shown. CB-DMD gives the most accurate reconstruction and the most reliable prediction of the periodic wave pattern.}
  \label{fig:CQGL_recon}
  \end{figure}

This spectral behavior is reflected in the reconstruction and prediction errors shown in Fig.~\ref{fig:CQGL_error}. The methods leading to relevant results are CB-DMD, improved linear error forcing, and the improved polynomial-forcing model. However, for improved linear error forcing and the improved polynomial-forcing method, the reported error curves are misleading: after a short transient, the corresponding reconstructions and predictions decay rapidly to zero, so the error essentially reduces to the norm of the true solution itself. This is visible in Fig.~\ref{fig:CQGL_recon}, illustrating the non-diverging reconstructions and predictions of the dynamics obained with the various approaches. 
A related effect occurs for improved standard DMD: The magnitude of the reconstructed dynamics is about three orders of magnitude larger than that of the true solution. CB-DMD is the only approach that avoids divergence, preserves the periodic structures, and matches the magnitude of the CQGL solution in reconstruction and prediction.
  
\section{Future work}
A useful feature of the present framework is that it can be formulated for any prescribed finite set of observables, including the original state coordinates. This avoids requiring an exactly Koopman-invariant observable subspace, but it also makes the choice of defect or forcing coordinates a central modelling decision.

This flexibility suggests several directions for future work. One direction is to identify observables or coordinate systems that reduce the size or energy of the fitted defect contribution. Another is to develop a principled rule for choosing the numerical threshold $\alpha$ in the correlation-basis construction and to analyze the stability of the resulting near-null-space directions under sampling noise and finite data.

Such extensions would make the forced linear representation more predictive and interpretable, and would clarify when the fitted linear component can be regarded as a reliable approximation of Koopman dynamics rather than only as a useful residual-splitting model.

\section{Conclusion}
We considered a theoretical framework for the forced linear representation of dynamical systems based on Koopman operator theory and the notion of defect space.  We showed that the approaches developed in \cite{PCT} and \cite{Khodkar2019} can be motivated by different modelling choices for the defect space. We proposed a specific defect space computable from the basis function of the null space of a correlation operator. We have shown that such basis functions can be approximated from data without any knowledge of the underlying differential equations and for high-dimensional measurements. Based on our choice of the defect space, we propose CB-DMD, an algorithm for learning effective linear operators for periodic nonlinear dynamical systems. The method aims to reduce the appearance of spurious modes while preserving a linear representation that captures the dominant periodic behavior. 

Across the nonlinear systems considered, CB-DMD consistently outperforms standard DMD, polynomial-forcing DMD, and error-forcing DMD. The experiments suggest that CB-DMD produces cleaner spectra, with dynamically relevant eigenvalues concentrated near the unit circle, and yields more stable, non-diverging reconstructions and predictions. These results suggest that correlation-basis forcing effectively separates unresolved nonlinear components from the autonomous linear operator, leading to linear models with fewer spurious modes for periodic nonlinear dynamics.

\nocite{}

\bibliographystyle{siamplain}
\bibliography{Koopmanlib.bib}

\section*{Appendix A}
We show that any finite-dimensional linear subspace $\mathcal E\subset \mathcal H$ is almost-invariant under the Koopman operator. That is, we show that there exists a finite-dimensional linear subspace $\mathcal V\subset \mathcal H$ such that $\mathcal E\cap \mathcal V=\{0\}$ and $\mathcal K g\in \mathcal E+\mathcal V$ for all $g\in \mathcal E$, where, as in the main text, $\mathcal E+\mathcal V:=\{h+\psi\mid h\in\mathcal E,\ \psi\in\mathcal V\}$.

Since $\mathcal E$ is finite dimensional and $\mathcal K$ is linear, its image $\mathcal K\mathcal E:=\{\mathcal K g:\ g\in \mathcal E\}$ is also a finite-dimensional linear subspace of $\mathcal H$. Now consider the subspace $\mathcal W:=\mathcal K\mathcal E\cap \mathcal E$. Since $\mathcal W$ is a subspace of the finite-dimensional space $\mathcal K\mathcal E$, we can choose a basis $\{w_1,\dots,w_r\}$ of $\mathcal W$ and extend it to a basis
\[
\{w_1,\dots,w_r,v_1,\dots,v_s\}
\]
of $\mathcal K\mathcal E$. Defining $\mathcal V:=\operatorname{span}\{v_1,\dots,v_s\}$, we obtain $\mathcal V\subset \mathcal K\mathcal E$ and $\mathcal K\mathcal E=\mathcal W+\mathcal V$. Moreover, if $f\in \mathcal W\cap \mathcal V$, then $f$ can be written both as a linear combination of the vectors $w_1,\dots,w_r$ and as a linear combination of the vectors $v_1,\dots,v_s$. Since the combined family is a basis of $\mathcal K\mathcal E$, it is linearly independent, so both representations must be trivial. Hence $f=0$, and therefore $\mathcal W\cap \mathcal V=\{0\}$.

We claim that $\mathcal E\cap \mathcal V=\{0\}$. Indeed, if $f\in \mathcal E\cap \mathcal V$, then, since $\mathcal V\subset \mathcal K\mathcal E$, we also have $f\in \mathcal E\cap \mathcal K\mathcal E=\mathcal W$. Hence $f\in \mathcal W\cap \mathcal V$, which implies $f=0$.

Finally, let $g\in \mathcal E$. Then $\mathcal K g\in \mathcal K\mathcal E$, and therefore $\mathcal K g=w+v$ for some $w\in \mathcal W$ and $v\in \mathcal V$. Since $\mathcal W\subset \mathcal E$, it follows that $\mathcal K g\in \mathcal E+\mathcal V$. Thus, $\mathcal E$ is almost-invariant under the Koopman operator. Note that, the same procedure shows that any finite-dimensional subspace of $\mathcal H$ is almost-invariant under any linear operator, and hence in particular under the Koopman operator.

\section*{Appendix B} 
We reformulate a standard finite-rank result for the correlation operator in a form adapted to our setting. 
We aim to show that the span of the observables \(\{g_i\}_{i=1}^N\) coincides with the span of the eigenfunctions of the correlation operator associated with nonzero eigenvalues.

Let
$
\mathcal{H}:=L^2_\mu(\mathcal{M};\mathbb{C})
$
be the Hilbert space of square-integrable complex-valued functions on \(\mathcal{M}\), endowed with inner product
$
\langle f,h\rangle_{\mathcal H}:=\int_{\mathcal M} f(\mathbf{x})\overline{h(\mathbf{x})}\,d\mu(\mathbf{x}).
$
Let \(\{g_i\}_{i=1}^N\subset \mathcal H\) be the observables chosen for the dynamical system, and define the correlation kernel
\[
c(\mathbf{x},\mathbf{y}):=\sum_{i=1}^N g_i(\mathbf{x})\overline{g_i(\mathbf{y})}.
\]
The associated correlation operator \(\mathcal C:\mathcal H\to\mathcal H\) is defined by
\[
(\mathcal C h)(\mathbf{x})
:=\int_{\mathcal M} c(\mathbf{x},\mathbf{y})\,h(\mathbf{y})\,d\mu(\mathbf{y}),
\qquad \forall h\in\mathcal H.
\]
Using the definition of \(c\), we obtain
\[
(\mathcal C h)(\mathbf{x})
=\sum_{i=1}^N g_i(\mathbf{x})
\int_{\mathcal M}\overline{g_i(\mathbf{y})}\,h(\mathbf{y})\,d\mu(\mathbf{y})
=\sum_{i=1}^N g_i(\mathbf{x})\langle h,g_i\rangle_{\mathcal H}.
\]
Since \(\operatorname{span}\{g_1,\dots,g_N\}\) is finite-dimensional, it follows that the range of \(\mathcal{C}\) is finite-dimensional as well. Hence \(\mathcal{C}\) is a finite-rank operator. Indeed, every finite-rank operator on a
Hilbert space is Hilbert-Schmidt. Moreover, for every \(h,k\in\mathcal H\),
\begin{align*}
\langle \mathcal C h,k\rangle_{\mathcal H}
&=
\left\langle
\sum_{i=1}^N g_i\,\langle h,g_i\rangle_{\mathcal H},
\,k
\right\rangle_{\mathcal H} 
=
\sum_{i=1}^N \langle h,g_i\rangle_{\mathcal H}\,\langle g_i,k\rangle_{\mathcal H} \\
&=
\sum_{i=1}^N \langle h,g_i\rangle_{\mathcal H}\,
\overline{\langle k,g_i\rangle_{\mathcal H}} 
=
\left\langle
h,\sum_{i=1}^N g_i\,\langle k,g_i\rangle_{\mathcal H}
\right\rangle_{\mathcal H} 
=
\langle h,\mathcal C k\rangle_{\mathcal H}.
\end{align*}
so \(\mathcal C\) is self-adjoint. Also,
\[
\langle \mathcal C h,h\rangle_{\mathcal H}
=\sum_{i=1}^N |\langle h,g_i\rangle_{\mathcal H}|^2 \ge 0,
\]
hence \(\mathcal C\) is positive.
Let us define
$
\mathcal H_N:=\operatorname{span}\{g_i\}_{i=1}^N
$
and
\[
S:=\operatorname{span}\{\varphi\in\mathcal H \;|\; \mathcal C\varphi=\lambda\varphi
\text{ for some } \lambda\neq 0\}.
\]
We now show that
$
\mathcal H_N=S.
$
First, from
\[
(\mathcal C h)(\mathbf{x})=\sum_{i=1}^N g_i(\mathbf{x})\langle h,g_i\rangle_{\mathcal H},
\]
it follows immediately that
\[
\operatorname{Ran}(\mathcal C)\subseteq \mathcal H_N.
\]
Here \(\operatorname{Ran}(\mathcal C)\) denotes the range of \(\mathcal C\), that is,
$
\operatorname{Ran}(\mathcal C):=\{\mathcal C h \;:\; h\in\mathcal H\},
$
namely the set of all functions that can be obtained as \(\mathcal C h\) for some \(h\in\mathcal H\).
Next, consider the restriction
\[
\mathcal C|_{\mathcal H_N}:\mathcal H_N\to\mathcal H_N.
\]
We claim that this restriction is injective. Indeed, let \(h\in\mathcal H_N\) and assume that \(\mathcal C h=0\). Then
\[
0=\langle \mathcal C h,h\rangle_{\mathcal H}
=\sum_{i=1}^N |\langle h,g_i\rangle_{\mathcal H}|^2.
\]
Hence \(\langle h,g_i\rangle_{\mathcal H}=0\) for all \(i=1,\dots,N\), which implies that \(h\in\mathcal H_N^\perp\). Since \(h\in\mathcal H_N\) as well, we conclude that \(h=0\). Thus \(\mathcal C|_{\mathcal H_N}\) is injective.

Because \(\mathcal H_N\) is finite-dimensional, injectivity implies surjectivity. Therefore
\[
\operatorname{Ran}(\mathcal C|_{\mathcal H_N})=\mathcal H_N.
\]
Since \(\operatorname{Ran}(\mathcal C|_{\mathcal H_N})\subseteq \operatorname{Ran}(\mathcal C)\subseteq \mathcal H_N\), we conclude that
\[
\operatorname{Ran}(\mathcal C)=\mathcal H_N.
\]

We now prove that \(S\subseteq \mathcal H_N\). Let \(\varphi\in\mathcal H\) be an eigenfunction of \(\mathcal C\) associated with an eigenvalue \(\lambda\neq 0\). Then
\[
\varphi=\frac{1}{\lambda}\mathcal C\varphi,
\]
and therefore \(\varphi\in \operatorname{Ran}(\mathcal C)=\mathcal H_N\). This shows that
\begin{equation}
\label{S_to_Hn}
S\subseteq \mathcal H_N.
\end{equation}
To prove the converse inclusion, recall that \(\mathcal C\) is compact and self-adjoint on \(\mathcal H\). By the spectral theorem for compact self-adjoint operators, \(\mathcal H\) decomposes orthogonally as
\[
\mathcal H= \nullspace(\mathcal C) + S,
\]
where $\nullspace(\mathcal C):= \{f \in \mathcal H \; | \; \mathcal C f =0\}$ is the null space of $\mathcal C$ and  $S:=\operatorname{span}\{\varphi\in\mathcal H \;|\; \mathcal C\varphi=\lambda\varphi
\text{ for some } \lambda\neq 0\}$ is the span of the eigenfunctions of \(\mathcal C\) associated with nonzero eigenvalues. We have already shown that
\[
\operatorname{Ran}(\mathcal C)=\mathcal H_N.
\]
Let \(h\in\mathcal H_N\). Then there exists \(f\in\mathcal H\) such that
$
h=\mathcal C f.
$
Using the above orthogonal decomposition, we may write
\[
f=f_0+f_S,
\qquad f_0\in \nullspace(\mathcal C),\quad f_S\in S.
\]
Therefore
\[
h=\mathcal C f=\mathcal C f_0+\mathcal C f_S=\mathcal C f_S.
\]
Since \(S\) is spanned by eigenfunctions associated with nonzero eigenvalues, it is invariant under \(\mathcal C\), and hence \(\mathcal C f_S\in S\). Thus \(h\in S\). Since \(h\in\mathcal H_N\) was arbitrary, we conclude that
\begin{equation}
\label{Hn_to_S}
    \mathcal H_N\subseteq S.
\end{equation}
Combining the two inclusions in (\ref{S_to_Hn}) and (\ref{Hn_to_S}), we obtain
\[
\mathcal H_N=S.
\]

In conclusion, the span of the observables \(\{g_i\}_{i=1}^N\) coincides with the span of the eigenfunctions of the correlation operator associated with nonzero eigenvalues. Note that this also implies that $\nullspace(\mathcal C) = \mathcal{H}_N^{\perp}$. If, in addition, \(c(\mathbf{x},\mathbf{y})\) is continuous on \(\mathcal M\times\mathcal M\) and \(\mathcal M\) is compact, then Mercer’s theorem yields the representation
\[
c(\mathbf{x},\mathbf{y})
=\sum_j \lambda_j \varphi_j(\mathbf{x})\overline{\varphi_j(\mathbf{y})},
\]
where \(\{\varphi_j\}\) is an orthonormal family of eigenfunctions associated with the nonzero eigenvalues \(\lambda_j>0\). 
\section*{Appendix C}
In Section~\ref{results}, we compare CB-DMD with several related methods. CB-DMD uses a theoretically motivated choice of forcing vectors derived in Section~\ref{correlation_basis_functions}, namely the numerical null space of a correlation matrix. In all cases, this null space is determined using the bound in \eqref{eq:alpha}; the only difference is the value of the relative threshold $\alpha$. We consider two choices for this parameter: a default value given by machine precision, as in SciPy's \texttt{null\_space} implementation \cite{scipy}, and a data-dependent value selected by hyperparameter tuning. The default choice works well for NLS and Cylinder wake wave systems, but for the CQGL equation, it is too restrictive and leads to diverging dynamics, so we use a larger value of $\alpha$ instead. In this appendix, we report the results for the other two systems using this larger threshold, as well as the CQGL results obtained with the default one.
\subsection*{Nonlinear Schroedinger equation}
\begin{figure}
  \centering
    \includegraphics[height=5cm]{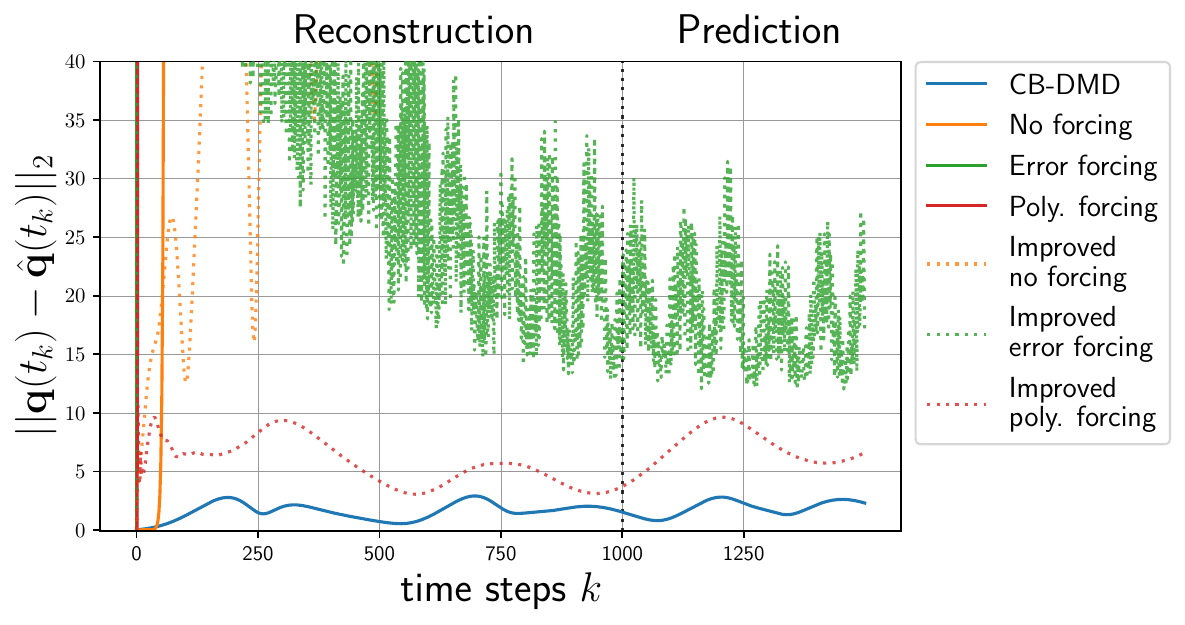}
    \includegraphics[height=5cm]{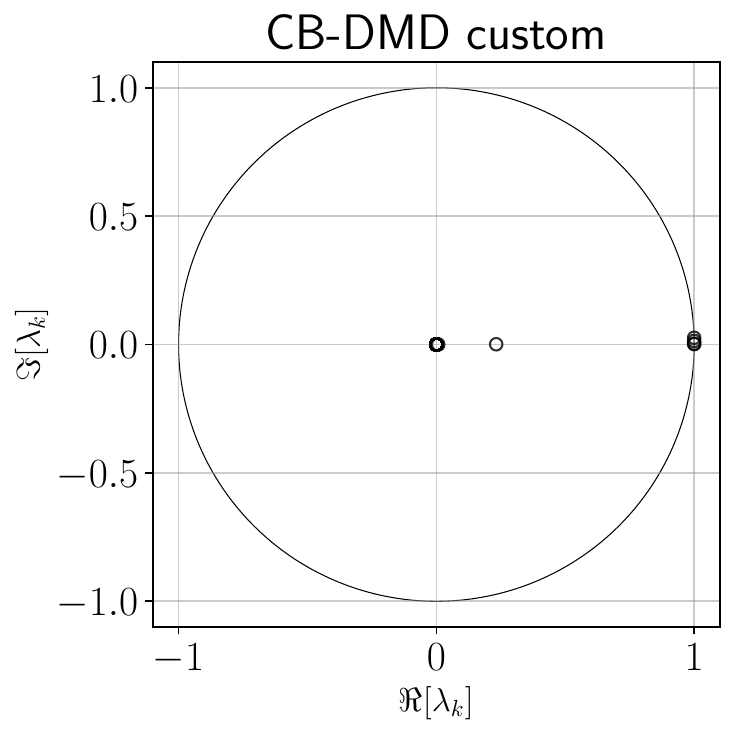}
\caption{Absolute error (left) and spectrum (right) of the CB-DMD operator for the NLS equation when the null space is computed with the custom threshold $\alpha=10^{-8}$. In the left panel, the vertical dashed black line marks the start of the prediction period, and the colored dashed lines show the errors of the improved versions of the corresponding operators. In the right panel, the circular line indicates the unit circle. Even with this custom threshold, CB-DMD remains the most effective method in terms of reconstruction and prediction accuracy.}
  \label{fig:NLS_app_error}
  \end{figure}

Fig.~\ref{fig:NLS_app_error} shows reconstruction and prediction errors obtained with the various approaches.
Here, CB-DMD was computed using $\alpha = 1e-8$. Compared to the default threshold, our error gets slightly worse. 
In this case, the operator learned by CB-DMD has only $6$ relevant eigenvalues, where relevant means $|\lambda|>10^{-3}$. Among them, $4$ lie close to the unit circle, in the sense that $0.99 \leq |\lambda| \leq 1.01$, and none has modulus greater than one. Visually, the resulting approximation is comparable to that obtained with the default threshold used to define the null space. Fig.~\ref{fig:NLS_app_error} shows that, even with the custom threshold, CB-DMD still yields the best reconstruction and prediction errors among the methods considered.

\subsection*{Cylinder wake flow}
\begin{figure}
  \centering
    \includegraphics[height=5cm]{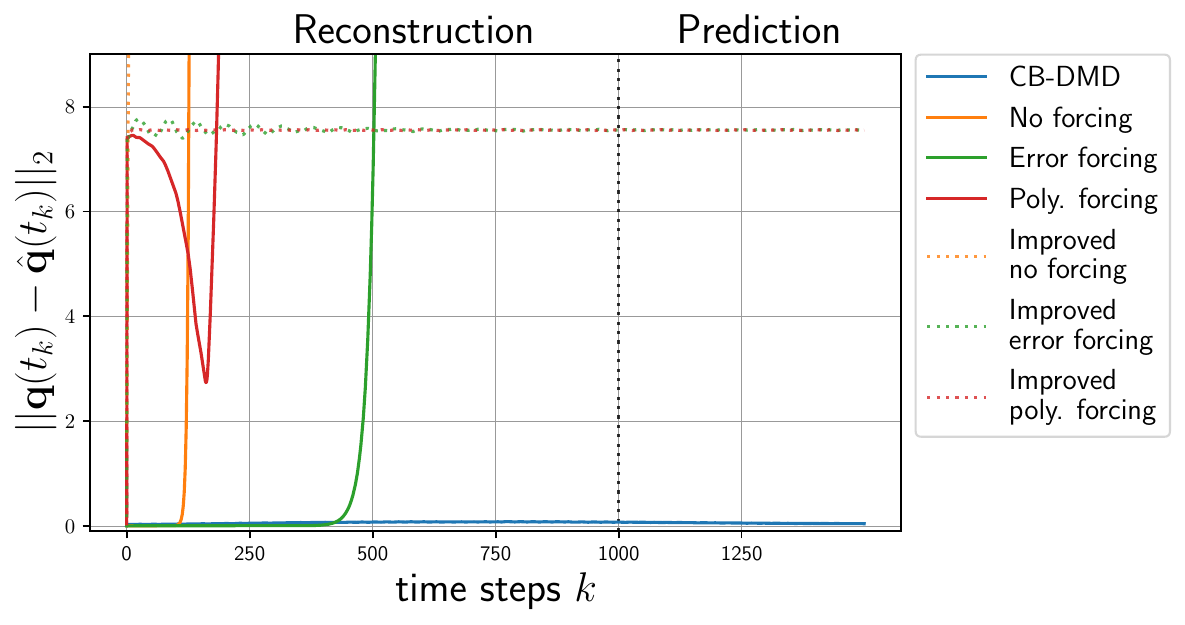}
    \includegraphics[height=5Cm]{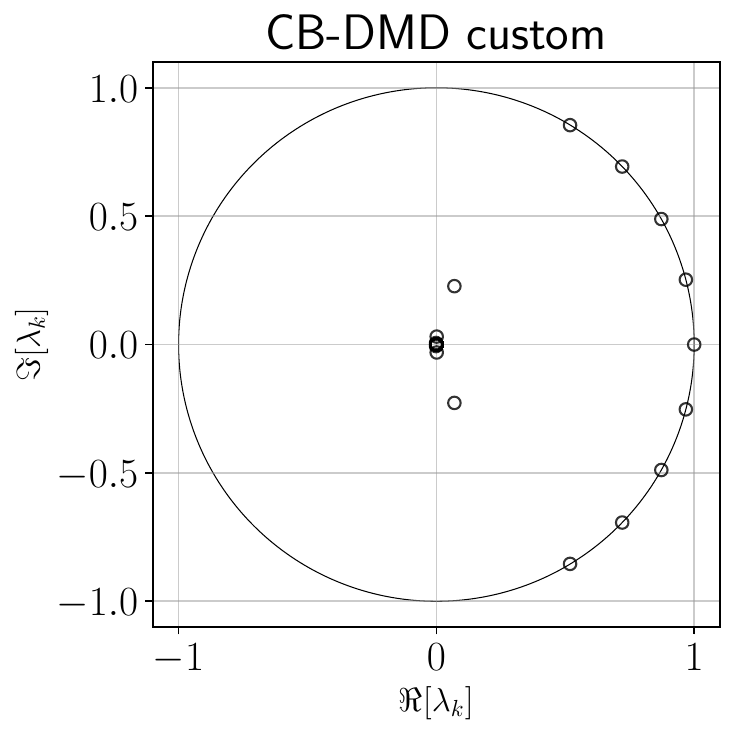}
  \caption{Absolute error (left) and spectrum (right) of the CB-DMD operator for the Cylinder wake system when the null space is computed with the custom threshold $\alpha=10^{-8}$. In the left panel, the vertical dashed black line marks the start of the prediction period, and the colored dashed lines show the errors of the improved versions of the corresponding operators. In the right panel, the circular line indicates the unit circle. Even with this custom threshold, CB-DMD remains the most effective method in terms of reconstruction and prediction accuracy.}
  \label{fig:Cylinder_app_error}
  \end{figure}

Fig.~\ref{fig:Cylinder_app_error} shows reconstruction and prediction errors obtained with the various approaches.
In these plots, the operator obtained with CB-DMD was computed using $\alpha = 1e-8$. Compared to the default threshold, our error gets slightly worse, similarly to what we observed for the NLS equation. 

  \subsection*{Cubic-quintic Ginzburg-Landau equation}
\begin{figure}
  \centering
    \includegraphics[height=5cm]{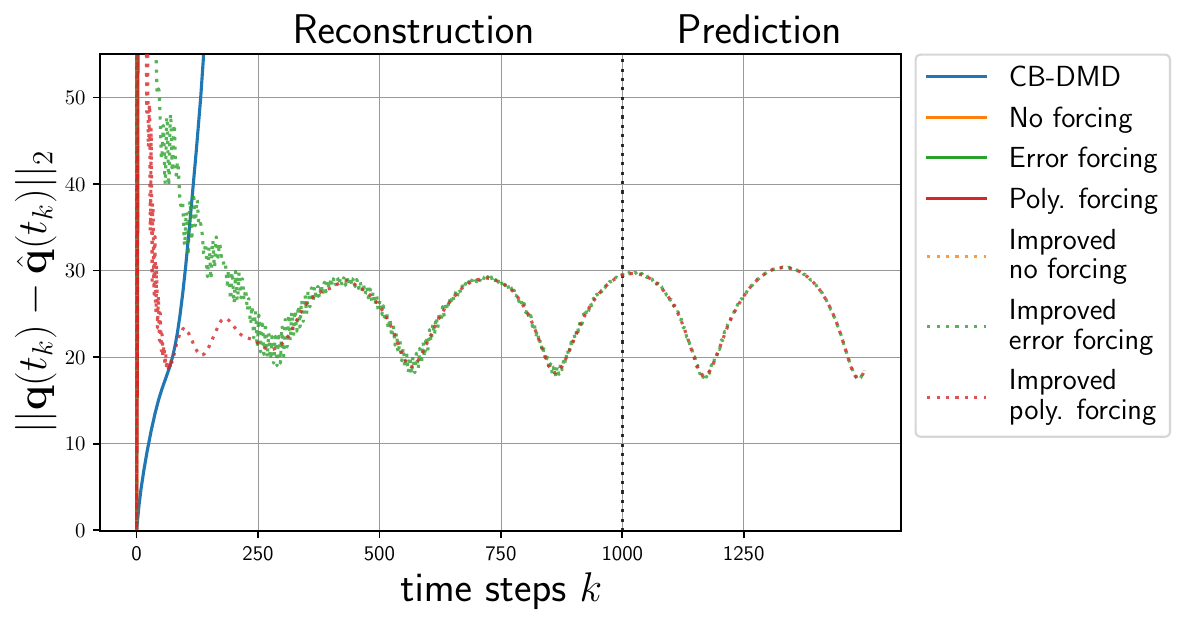}
    \includegraphics[height=5cm]{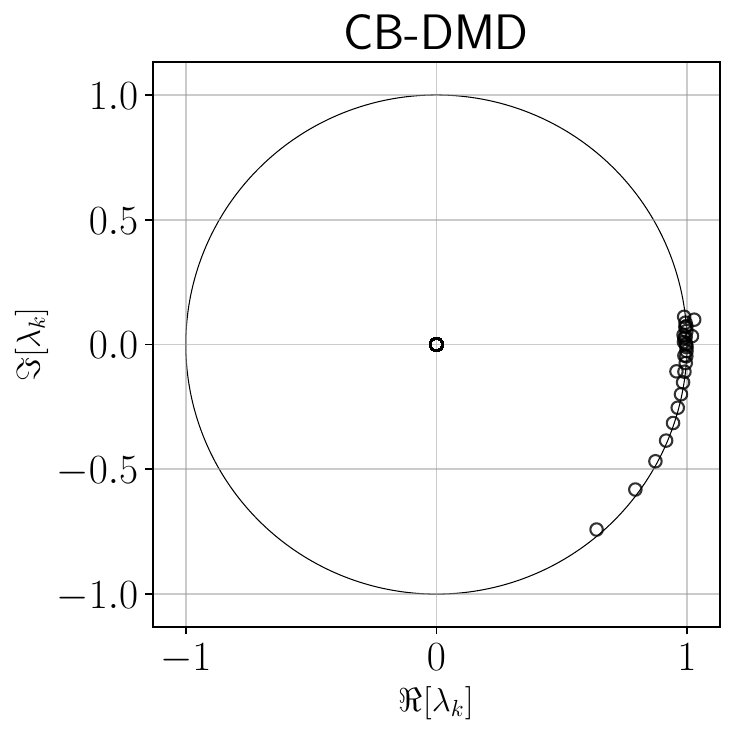}
\caption{Absolute error (left) and spectrum (right) of the CB-DMD operator for the CQGL equation when the null space is computed with the default threshold for $\alpha$. In the left panel, the vertical dashed black line marks the start of the prediction period, and the colored dashed lines show the errors of the improved versions of the corresponding operators. In the right panel, the circular line indicates the unit circle.. As reflected by the rapidly growing error, this choice of threshold leads to diverging dynamics for CB-DMD.}
  \label{fig:CQGL_app_error}
  \end{figure}

Fig.~\ref{fig:CQGL_app_error} shows reconstruction and prediction errors obtained with the various approaches.
These results are obtained using SciPy's \texttt{null\_space} method \cite{scipy}, that is, we set $\alpha$ as the default value given by machine precision. In this case, CB-DMD operator leads to diverging dynamics. 
In particular, the operator learned by CB-DMD has $28$ relevant eigenvalues, where relevant means $|\lambda|>10^{-3}$. Among them, $20$ lie close to the unit circle, in the sense that $0.99 \leq |\lambda| \leq 1.01$, while $3$ have modulus greater than one; the largest modulus is $1.0336$. As shown in Fig.~\ref{fig:CQGL_app_error}, the error of the resulting operator then grows rapidly.

\section*{Appendix D}
In the main text, we report improved reconstructions and predictions obtained by removing all modes associated with eigenvalues of modulus greater than one, since such modes can lead to diverging dynamics. This choice is, however, heuristic: modes with modulus only slightly above one may still improve reconstruction over a finite time horizon. To assess this effect, in this appendix, we consider the thresholds
$ \{0.99,\,0.999,\,0.9999,\,1,\,1.0001,$ $\,1.001,\,1.01\},$
together with the case of no thresholding, and compare the resulting mean absolute reconstruction errors. For each method, we then select the threshold that gives the best performance. 
For the baseline methods, the basis functions for the defect space have been determined according to the same procedures described in Section~\ref{results}. For CB-DMD, the basis functions are computed using the default choice of $\alpha$ for all systems. In the CQGLE results in Section~\ref{results}, we use a custom value of $\alpha$, since the default choice leads to diverging dynamics. Here, we instead consider the default choice to illustrate that, even when it introduces spurious modes, the main dynamical behavior can still be captured effectively once those modes are removed.

As in the main text, the improved linear error-forcing operator is constructed using the best reconstruction obtained from standard DMD. 
\subsection*{Nonlinear Schroedinger equation}
\begin{figure}
  \centering
      \includegraphics[height=4cm]{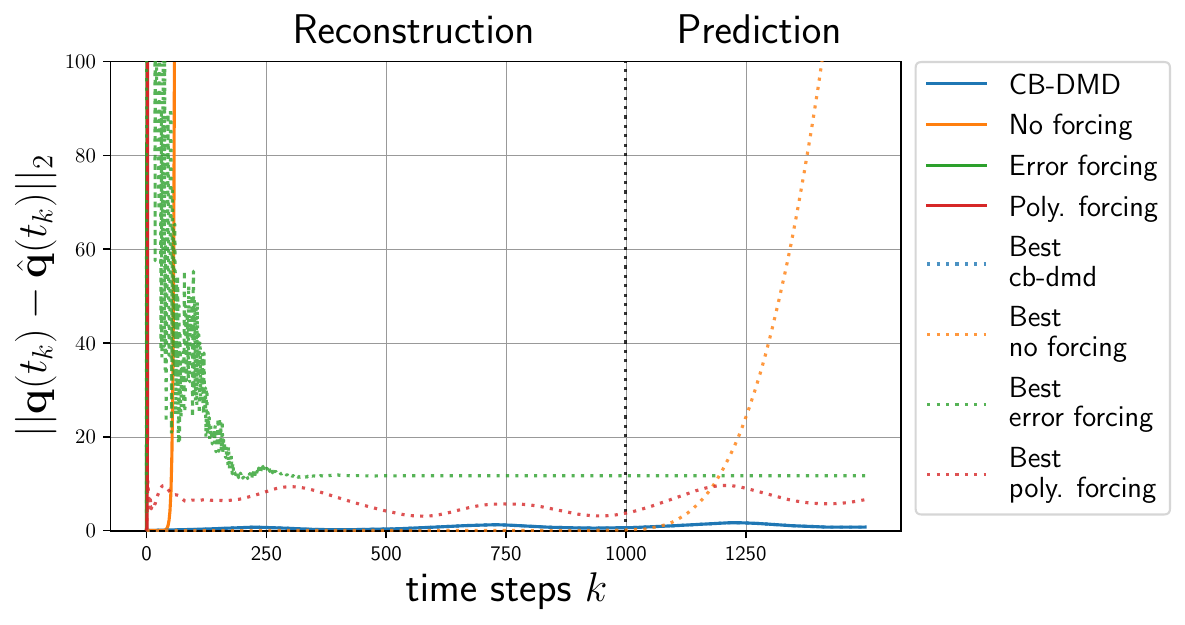}
      \includegraphics[height=4cm]{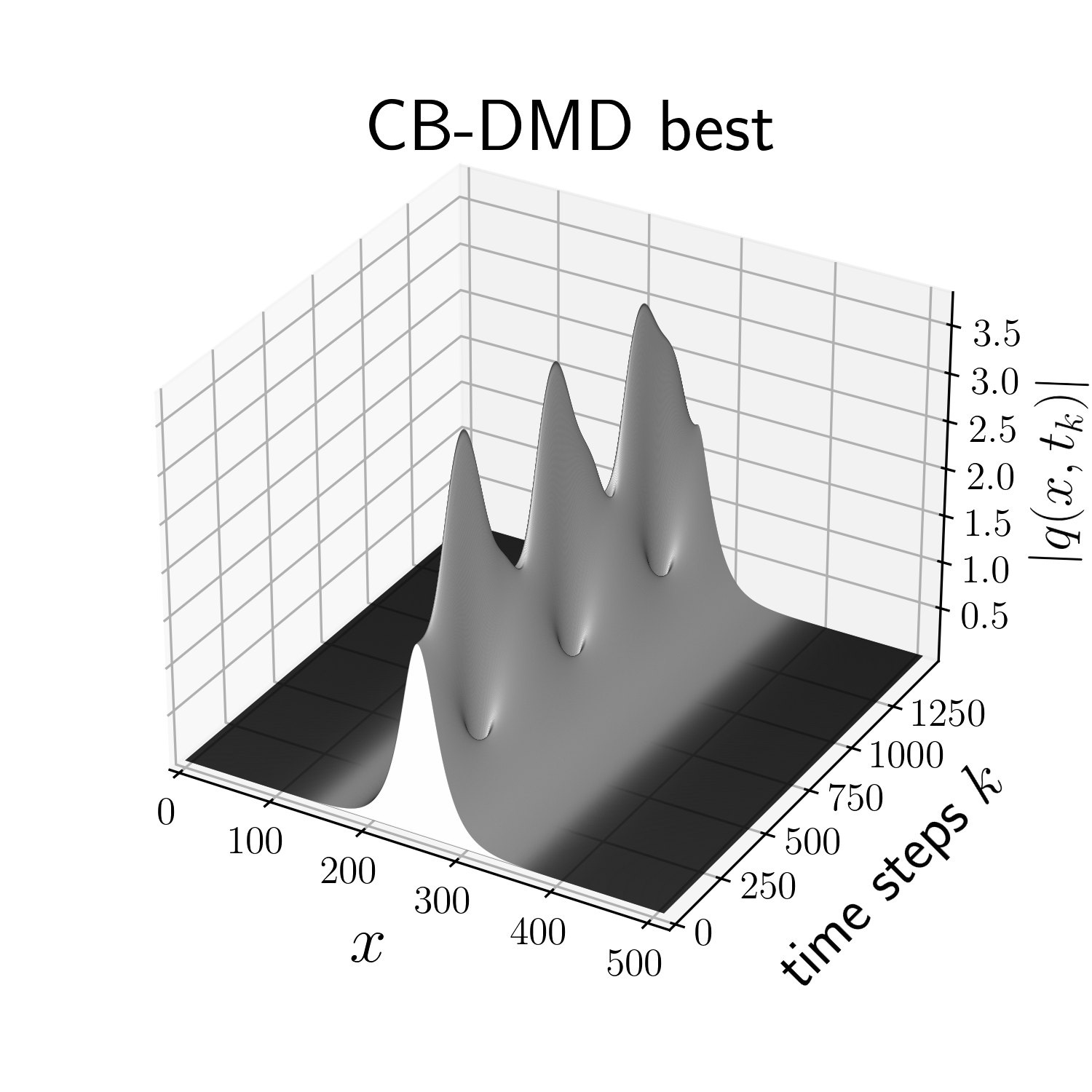}
\caption{Best eigenvalue-threshold comparison for the NLS equation. Left: absolute error of the reconstruction and prediction obtained by selecting, for each method, the eigenvalue-magnitude thres that gives the lowest mean absolute reconstruction error among the values tested. Dashed lines denote the corresponding best-threshold results, while solid lines denote the operators without eigenvalue thresholding. The vertical dashed black line marks the start of the prediction period. Right: reconstruction and prediction of CB-DMD for its best-performing threshold choice among the values tested. For CB-DMD, the best result is obtained without removing any modes, so it coincides with the result reported in Section~\ref{results}.}
  \label{fig:best_nls}
  \end{figure}

Fig.~\ref{fig:best_nls} shows the optimal reconstruction and prediction errors obtained with the various approaches in our setting.
For CB-DMD, the best performance is obtained without removing any modes, so the results reported in Section~\ref{results} are already optimal within this comparison. For the baseline methods, the best eigenvalue-magnitude thresholds are $1.01$ for standard DMD, $0.99$ for linear error forcing based on the best standard DMD reconstruction, and $1.0001$ for polynomial forcing.

\subsection*{Cylinder wake flow}
\begin{figure}
  \centering
      \includegraphics[height=4cm]{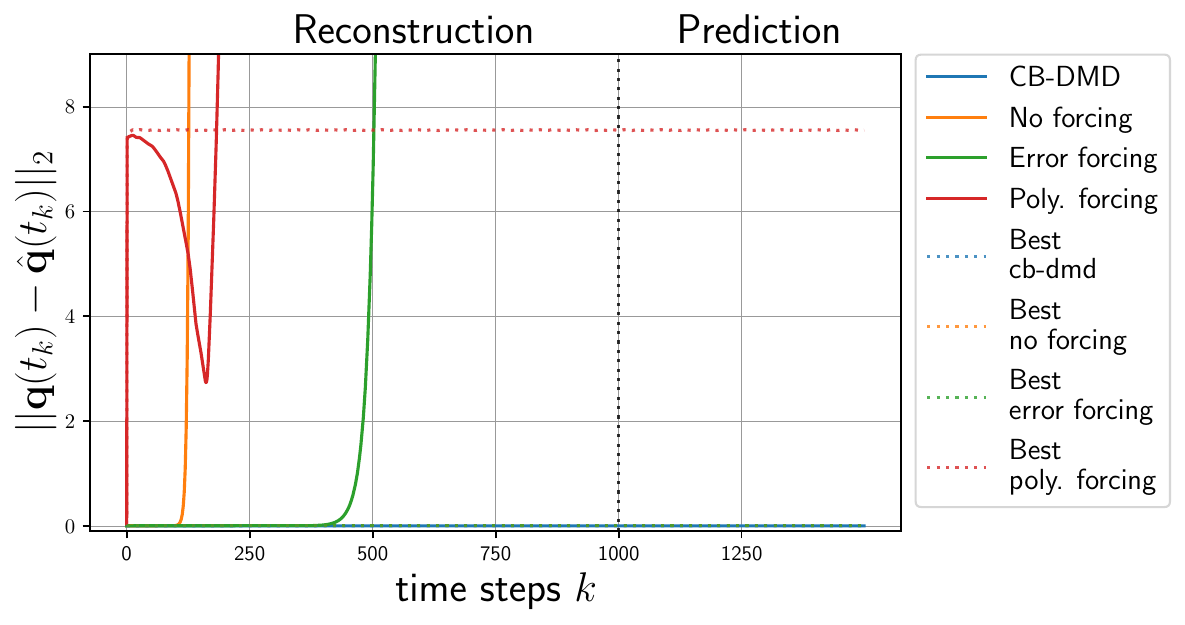}
  \caption{Best eigenvalue-threshold comparison for the cylinder wake dynamics. Left: absolute reconstruction and prediction errors obtained by selecting, for each method, the eigenvalue-magnitude threshold that minimizes the mean absolute reconstruction error among the tested values. 
  The vertical dashed black line marks the beginning of the prediction interval. With the optimal threshold, standard DMD and error forcing attain reconstruction and prediction errors comparable to those of CB-DMD}
  \label{fig:best_cylinder}
  \end{figure}

Fig.~\ref{fig:best_cylinder} shows the optimal reconstruction and prediction errors obtained with the various approaches in our setting.
The best result for the CB-DMD operator was obtained considering a threshold of 1.0001. 
The determined threshold for standard DMD is 1.0001, for the linear error forcing with respect to the best standard DMD reconstruction is 1.0001 and for the polynomial forcing is 1.. 
  \newpage

  \subsection*{Cubic-quintic Ginzburg-Landau equation}
\begin{figure}
  \centering
      \includegraphics[height=4cm]{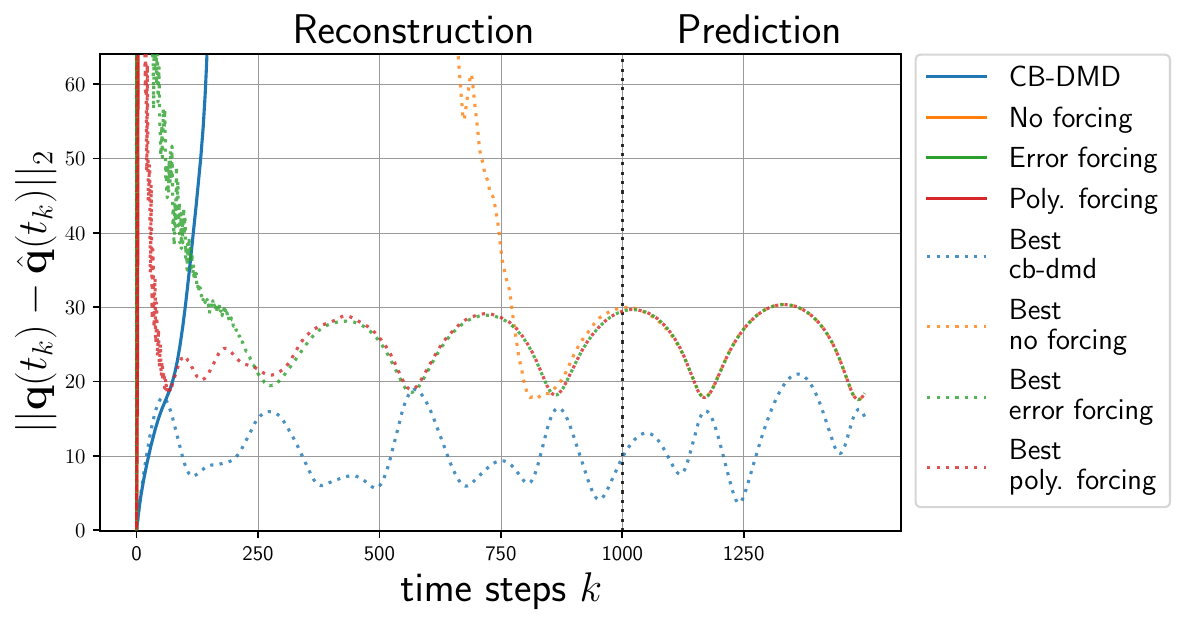}
      \includegraphics[height=4cm]{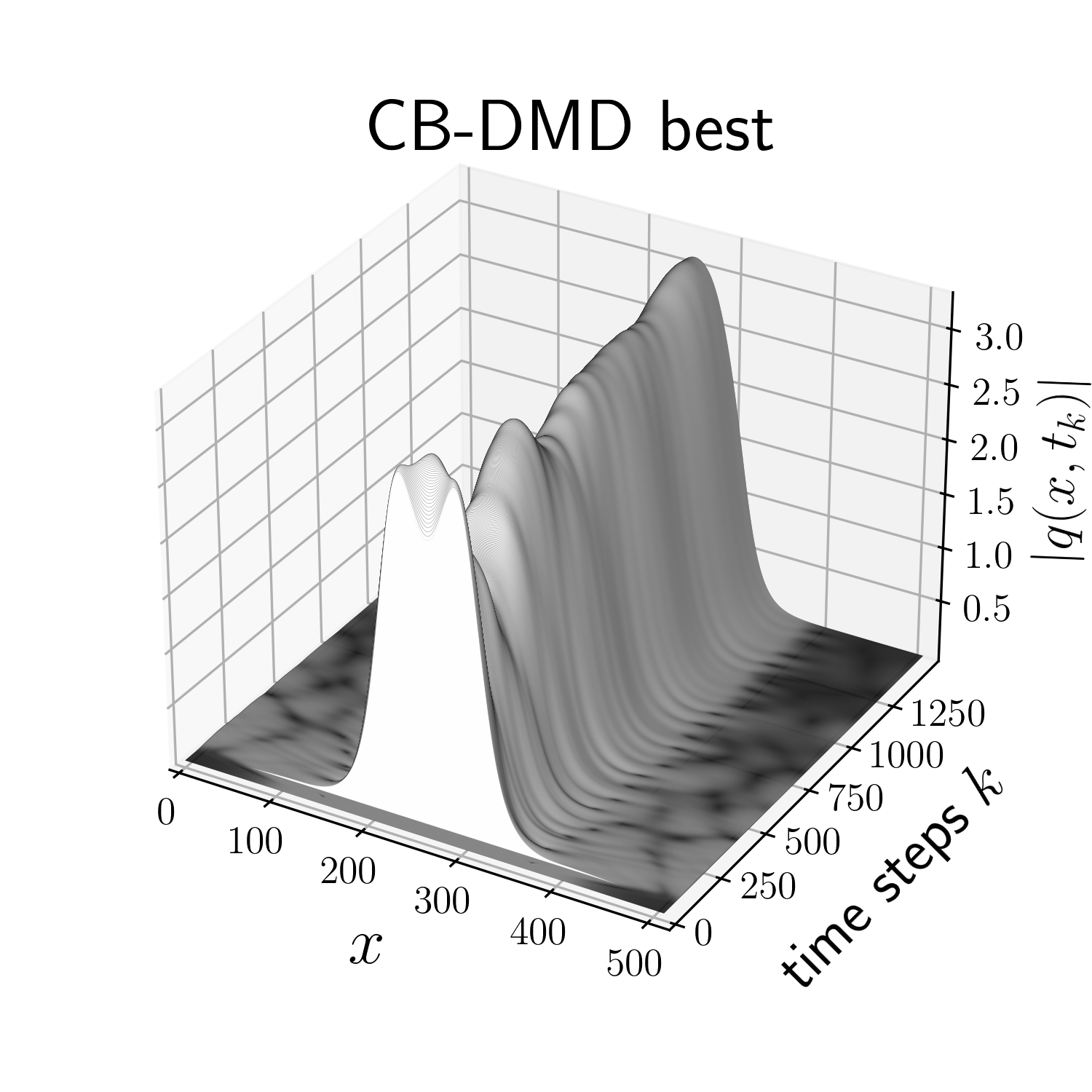}
\caption{Best eigenvalue-threshold comparison for the CQGL equation. Left: absolute error of the reconstruction and prediction obtained by selecting, for each method, the eigenvalue-magnitude threshold that gives the lowest mean absolute reconstruction error among the values tested. The vertical dashed black line marks the start of the prediction period. Right: reconstruction and prediction of CB-DMD corresponding to its best-performing threshold among the values tested.}
  \label{fig:best_cqgle}
  \end{figure}

Fig.~\ref{fig:best_cqgle} shows the optimal reconstruction and prediction errors obtained with the various approaches in our setting.
The selected eigenvalue-magnitude thresholds are $1.01$ for CB-DMD, $0.99$ for standard DMD, $1.0001$ for linear error forcing based on the best standard DMD reconstruction, and $1.0001$ for polynomial forcing. Fig. \ref{fig:best_cqgle} indicates that CB-DMD leads to the lowest reconstruction and prediction error.
\newpage
\section*{Appendix E}
This appendix includes the Dynamic Mode Decomposition Algorithm.
\begin{algorithm}[h]
  \caption{Dynamic Mode Decomposition (DMD)}\label{alg:DMD}
\begin{algorithmic}[1]
    \STATE Arrange the data into matrices $\bar{\mathbf{X}}$ and $\mathbf{X}$ as in
    $(\ref{snapshots})$. \STATE Compute the reduced singular value
    decomposition (SVD) of $\mathbf{X}$, i.e., $\mathbf{X}\approx \mathbf{U}_r \mathbf{\Sigma}_r \mathbf{V}^*_r$. Here
    $\mathbf{U}_r,$ $\mathbf{\Sigma}_r$ and $ \mathbf{V}^*_r$ are the rank-$r$ truncation of the
    matrices $\mathbf{U},\; \mathbf{\Sigma}$ and $ \mathbf{V}^*$ computed via SVD and such that  $\mathbf{X}= \mathbf{U}
    \mathbf{\Sigma} \mathbf{V}^*$. Further, $\mathbf{V}^*$ is the conjugate transpose of $\mathbf{V}$ and $r$ is
    a parameter to be chosen. \STATE Construct the matrix $\tilde{\mathbf{A}}:=
    \mathbf{U}^*_r\bar{\mathbf{X}}\mathbf{V}_r \mathbf{\Sigma}^{-1}_r$. \STATE Compute eigenvalues and eigenvectors
    of $\tilde{\mathbf{A}}$, solving $\tilde{\mathbf{A}}\mathbf{W}= \mathbf{\Lambda} \mathbf{W}$. With $\mathbf{\Lambda}$ a diagonal
    matrix of DMD eigenvalues and $\mathbf{W}$ a matrix where the columns are eigenvectors
    of $\tilde{\mathbf{A}}$. \STATE The DMD modes corresponding to DMD eigenvalues
    $\mathbf{\Lambda}$ are given by $\mathbf{\Phi}= \bar{\mathbf{X}}\mathbf{V}_r \mathbf{\Sigma}_r^{-1}\mathbf{W}$.

    \STATE The reconstruction and prediction is given by 
    \begin{equation}
    \label{prediction_wind}
    \hat{\mathbf{x}}_k= \mathbf{\Phi} \; \text{diag}( \text{exp}(\boldsymbol{\omega} k))\mathbf{b} \quad 
    i=0,\dots,T,
    \end{equation}
    where the amplitudes $\mathbf{b}$ are computed optimally as described in \cite{jovanovic14} or as $\mathbf{b}= \mathbf{\Phi}^{\dagger}\mathbf{x}_0$.

    The vector $\boldsymbol{\omega}$ contains the frequency values associated with each DMD mode obtained from the DMD eigenvalues $\lambda_k$:
\begin{equation}
  \omega_k= \log(\lambda_k)/\Delta t.
  \label{def:omega_frequency}
\end{equation} 
Here, $\Delta t$ represents
the time difference between two consecutive measurements.
\end{algorithmic}
\end{algorithm}

\end{document}

%% file: ex_shared.tex
\usepackage{lipsum}
\usepackage{amsfonts}
\usepackage{graphicx}
\usepackage{epstopdf}
\usepackage{algorithmic}
\ifpdf
  \DeclareGraphicsExtensions{.eps,.pdf,.png,.jpg}
\else
  \DeclareGraphicsExtensions{.eps}
\fi

\usepackage{enumitem}
\setlist[enumerate]{leftmargin=.5in}
\setlist[itemize]{leftmargin=.5in}

\newsiamremark{remark}{Remark}
\newsiamremark{hypothesis}{Hypothesis}
\crefname{hypothesis}{Hypothesis}{Hypotheses}
\newsiamthm{claim}{Claim}
\newsiamremark{fact}{Fact}
\crefname{fact}{Fact}{Facts}

\headers{A Koopman Framework}{P. Climaco, J. Garcke, and X. F. Gerloff}

\title{Forced-Term Modeling in the Koopman Framework for Linear Operator Learning in Nonlinear Dynamical Systems}

\author{Paolo Climaco\thanks{Department of Mathematics, University of California, Los Angeles, USA
  (\email{climaco@math.ucla.edu}).}
\and Jochen Garcke\thanks{Institut für Numerische Simulation, Universität Bonn, Germany 
  (\email{garcke@ins.uni-bonn.de}) and Fraunhofer SCAI, Sankt Augustin, Germany.}
\and Xenia F. Gerloff\thanks{Work done while at Institut für Numerische Simulation, Universität Bonn, Germany (\email{xeniagerloff@gmail.com}) and Fraunhofer SCAI, Sankt Augustin, Germany.}}

\usepackage{amsopn}
